 \documentclass[3p]{elsarticle}
\usepackage{graphicx}
\usepackage{amsmath,amssymb,mathtools,amsfonts}
\usepackage{hyperref}
\usepackage{multirow}
\usepackage{color}
\usepackage{psfrag}
\usepackage{empheq}
\usepackage{tabularx}
\usepackage{bm} 
\usepackage{url}
\usepackage{float}

\usepackage{pdflscape}
\usepackage{quiver}
\usepackage{tikz}
\usepackage{pgfplots}
\usetikzlibrary{arrows.meta}
\tikzstyle{block} = [draw, fill=white, rectangle, 
minimum height=2.5em, minimum width=5.em]
\tikzstyle{sum} = [draw, fill=white, circle, node distance=1cm]
\tikzstyle{input} = [coordinate]
\tikzstyle{output} = [coordinate]
\tikzstyle{pinstyle} = [pin edge={to-,thin,black}]
\newlength\fwidth
\newlength\fheight
\usetikzlibrary{arrows}
\usetikzlibrary{cd}
\usetikzlibrary{babel}
\usetikzlibrary{matrix}

\usepackage{algorithm}
\usepackage{algpseudocode}
\usepackage{comment}
\usepackage{caption, subcaption}
\usepackage{booktabs}

\newcommand{\R}{{\mathbb{R}}}
\newcommand{\dd}{{\mathrm{d}}}

\newcommand{\Nsys}{{n_{f}}} 
\newcommand{\Nstg}{{n_\mathrm{s}}} 
\newcommand{\NFE}{N_\mathrm{FE}} 
\newcommand{\Nsim}{{N_\mathrm{sim}}} 
\newcommand{\Nswitch}{{N_\mathrm{sw}}} 
\newcommand{\I}{{\mathcal{I}}}
\renewcommand{\I}{{\mathcal{I}}}
\newcommand{\C}{{\mathcal{C}}}
\newcommand{\J}{{\mathcal{J}}}
\newcommand{\K}{{\mathcal{K}}}
\newcommand{\G}{\mathcal{G}}
\newcommand{\ILift}{\mathcal{I}_{j}}
\newcommand{\lambdan}{\lambda^{\mathrm{n}}}
\newcommand{\lambdap}{\lambda^{\mathrm{p}}}

\newcommand{\Lambdan}{\Lambda^{\mathrm{n}}}
\newcommand{\Lambdap}{\Lambda^{\mathrm{p}}}

\newcommand{\irk}{{\mathrm{rk}}}

\newcommand{\fesd}{{\mathrm{fesd}}}

\newcommand{\Nctrl}{N}

\newcommand{\ts}{{t_{\mathrm{s}}}}
\newcommand{\tsn}{t_{\mathrm{s},n}}
\newcommand{\tsnn}{{t_{\mathrm{s},n+1}}}
\newcommand{\tsnhat}{{\hat{t}_{\mathrm{s},n}}}
\newcommand{\tsnnhat}{{\hat{t}_{\mathrm{s},n+1}}}
\newcommand{\hatts}{{\hat{t}_{\mathrm{s}}}}

\newcommand{\switchfun}{\psi}
\newcommand{\nswitchfun}{n_{\psi}}
\newcommand{\Ffilippov}{\mathcal{F}_{\mathrm{F}}}
\newcommand{\Fstep}{\mathcal{F}_{\mathrm{H}}}
\newcommand{\Fap}{\mathcal{F}_{\mathrm{AP}}}
\newcommand{\F}{\mathcal{F}}
\newcommand{\gdcseq}{g_{\mathrm{e}}}
\newcommand{\gdcscomp}{g_{\mathrm{c}}}

\newtheorem{theorem}{Theorem}
\newtheorem{example}{Example}
\newtheorem{remark}[theorem]{Remark}
\newtheorem{definition}[theorem]{Definition}
\newtheorem{assumption}[theorem]{Assumption}

\newtheorem{proposition}[theorem]{Proposition}

\newtheorem{lemma}[theorem]{Lemma}
\newcommand{\nosnoc}{{\texttt{nosnoc}}}

\newcommand{\secsub}[1]{\textcolor{blue}{#1}}

\DeclareMathSymbol{\shortminus}{\mathbin}{AMSa}{"39}

\usepackage[shortlabels]{enumitem} 

\journal{Preprenit}

\begin{document}

\begin{frontmatter}
	
\title{Finite Elements with Switch Detection for Numerical Optimal Control of Nonsmooth Dynamical Systems with Set-Valued Heaviside Step Functions}

\fntext[label1]{Corresponding author.} 	
\fntext[label2]{This research was supported by th DFG via Research Unit FOR 2401 and project 424107692, by the EU via ELO-X 953348, by the German Federal Ministry for Economic Affairs and Climate Action (BMWK.IIC6) via the project WOpS with the project number 03EN3054A.	
\textit{Email address}: \texttt{\{armin.nurkanovic,jonathan.frey,moritz.diehl\} @imtek.uni-freiburg.de, anton.pozharskiy@merkur.uni-freiburg.de}.}

\author[FreiburgImtek,label1]{Armin Nurkanovi\'c}
\author[FreiburgImtek]{Anton Pozharskiy}
\author[FreiburgImtek,FreiburgMath]{Jonathan Frey}
\author[FreiburgImtek,FreiburgMath]{Moritz Diehl}

\affiliation[FreiburgImtek]{organization={Department of Microsystems Engineering (IMTEK), University of Freiburg}, country={Germany}}
\affiliation[FreiburgMath]{organization={Department of Mathematics, University of Freiburg}, country={Germany}}

\begin{abstract}
This paper develops high-accuracy methods for numerically solving optimal control problems subject to nonsmooth differential equations with set-valued Heaviside step functions. 
An important subclass of these systems are Filippov systems.
Writing the Heaviside step function as the solution map of a linear program and using its optimality conditions, the initial nonsmooth system is rewritten into an equivalent dynamic complementarity system (DCS).
We extend the Finite Elements with Switch Detection (FESD) method~\cite{Nurkanovic2024a}, initially developed for Filippov systems transformed via Stewart's reformulation into DCS \cite{Stewart1990a}, to the above mentioned class of nonsmooth systems.
The key ideas are to start with a standard Runge-Kutta method for the DCS and to let the integration step sizes to be degrees of freedom. 
Next, we introduce additional conditions to enable implicit but exact switch detection and to remove possible spurious degrees of freedom if no switches occur.
The theoretical properties of the method are studied. 
Its favorable properties are illustrated on numerical simulation and optimal control examples. 
All methods introduced in this paper are implemented in the open-source software package \nosnoc~\cite{Nurkanovic2022b}.
\end{abstract}



\begin{keyword}
numerical simulation, Filippov systems, numerical optimal control, discontinuous differential equations, differential inclusions
\end{keyword}
\end{frontmatter}


\graphicspath{{image/}}
\section{Introduction}\label{sec:intro}
The set-valued version of the Heaviside step function, denoted as $\gamma : \R \to \mathcal{P}(\R)$, is defined as follows:
\begin{align}\label{eq:step_function}
	\gamma(y) &= \begin{cases}
		\{1\}, & y > 0, \\
		[0,1], & y = 0, \\
		\{0\}, & y < 0.
	\end{cases}
\end{align}
Here $\mathcal{P}(\R)$ represents the power set of $\R$. The set-valued Heaviside step function, often referred to simply as the step function, provides an intuitive way to model Boolean relationships within a dynamical system.

In the modeling process, \textit{switching functions} $\psi_i: \R^{n_x} \to \R$, for $i = 1,\ldots,\nswitchfun$, are commonly used as arguments of the step functions. We denote the concatenation of all scalar Heaviside step functions as $\Gamma(\switchfun(x)) \coloneqq (\gamma(\switchfun_1(x)),\dots,\gamma(\switchfun_{\nswitchfun}(x)))$, where $\switchfun(x) = (\switchfun_1(x),\dots,\switchfun_{\nswitchfun}(x))$.

In this paper, we study nonsmooth dynamical systems of the form:
\begin{align}\label{eq:basic_di}
	\dot{x} \in \F(x,u,\Gamma((\psi(x)))).
\end{align}
Here, $u \in \R^{n_u}$ represents a known control function that could be obtained, for example, by solving an optimal control problem. The equation \eqref{eq:basic_di} is an instance of a differential inclusion (DI) since the right-hand side is set-valued due to the function $\Gamma$, which can enter $\mathcal{F}$ in a nonlinear way. Thus, the set $\mathcal{F}(x,u,\Gamma((\psi(x))))$ can be convex or nonconvex.

The nonsmooth and set-valued nature of $\mathcal{F}(x,u,\Gamma((\psi(x))))$ in the DI \eqref{eq:basic_di} introduces complexities in its numerical treatment that we address in this paper.
\color{black}
In particular, we develop accurate and efficient numerical methods for \color{black} solving simulation and optimal control problems subject to such DIs. \color{black}
An important subclass of nonsmooth systems with Heaviside step functions are Filippov systems~\cite{Acary2014,Dieci2011,Filippov1988}.
Moreover, some classes of systems with state jumps can be reformulated into Filippov systems via the time-freezing reformulation \cite{Halm2023,Halm2023a,Nurkanovic2023a,Nurkanovic2022a,Nurkanovic2021}.
Therefore, this modeling approach enables us to treat a broad class of practical problems.

\color{black}
In a Piecewise Smooth System (PSS), the state space is split into nonempty regions where each of them is equipped with a different vector field.
Step functions are used in PSS modeling to determine in which region or on what boundaries the trajectory is.
Smoothed single-valued versions of the Heaviside step function are often used in the numerical treatment of a PSS~\cite{Bernardo2008,Guglielmi2022,Machina2011}.
\color{black}
A prominent application example of step functions are gene-regulatory networks~\cite{Acary2014,Machina2011}.
Another common application of step functions is to compute selections of Filippov sets in sliding modes~\cite{Dieci2009,Dieci2011}.

\color{black}
The event of $x(t)$ becoming nondifferentiable is called a \textit{switch}. 
In our case, this corresponds that one or more components of $\switchfun(x)$ become zero, or if they were zero, they become nonzero.
\color{black}
Many mature numerical simulation methods to treat ODEs with switches exist, and the corresponding theory is well-established~\cite{Acary2008}.
However, numerical methods for solving Optimal Control Problems (OCPs) subject to nonsmooth dynamical systems are not yet at such a mature stage.
The key ingredient to high-accuracy methods and efficient numerical optimal control is to detect the time points when the switches occur~\cite{Nurkanovic2024a}.
In the control community, a popular approach to deal with the nonsmoothness is to introduce integer variables to label all modes of the nonsmooth system~\cite{Bemporad1999b}.
This leads to mixed-integer optimization problems for solving OCPs.
However, they often become computationally intractable when nonconvexities appear or exact switching times need to be computed.

\color{black}
On the one hand, the computationally usually less favorable indirect methods for optimal control are not widespread since Pontryagin-like conditions are not established for many classes of nonsmooth systems, cf.~\cite{Guo2016,Bouali2023} for an overview.
\color{black}
On the other hand, the application of direct methods, i.e., in a first-discretize-then-optimize approach, is not as straightforward as for OCPs with smooth ODEs and can lead to spurious solutions and wrong conclusions.
This was first rigorously explained in the seminal paper of Stewart and Anitescu~\cite{Stewart2010}.
\color{black}
They show that the derivatives of an integration map obtained by time-stepping methods (e.g., the implicit Euler method) with respect to initial values and controls (called numerical sensitivities) do not converge to corresponding continuous-time values (called sensitivities), no matter how small the integrator step size is. 
In practice, this can result in termination at feasible, but nonoptimal points~\cite{Nurkanovic2020}.
\color{black}
Moreover, they show that the numerical sensitivities of the smoothed approximations of a nonsmooth system are only correct if the step size shrinks faster than the smoothing parameter.
This makes the smoothing approach impractical, as very small step sizes are needed even for moderate accuracy.

The limitations of direct methods based on time-stepping were recently overcome by the method of Finite Elements with Switch Detection (FESD)~\cite{Nurkanovic2024a}.
This method is based on standard Runge-Kutta (RK) discretizations of the Dynamic Complementarity systems (DCS), where the integrator step sizes are left as degrees of freedom as first proposed by Baumrucker and Biegler~\cite{Baumrucker2009}.
Additional constraints are introduced to have a well-defined system with exact switch detection.
The discretization yields Mathematical Programs with Complementarity Constraints (MPCC).
They are nonregular and nonsmooth Nonlinear Programs (NLP)~\cite{Anitescu2007,Scholtes2001}.
Still, with suitable reformulations and homotopy procedures, they can be solved efficiently using techniques from smooth optimization.

Initially, we developed FESD for Stewart's reformulation of Filippov systems into DCS~\cite{Stewart1990a}.
In this paper, we extend these ideas to nonsmooth systems of the form of \eqref{eq:basic_di} and corresponding OCPs. 
This enables us to cover a more general class of nonsmooth ODEs with a discontinuous r.h.s..

\paragraph*{Contributions}
\color{black}
We provide a detailed study of the transformation of a Filippov system into a dynamic complementarity system (DCS) via set-valued Heaviside step functions.
The key step in this transformation is to view the Heaviside step function as the solution of a parametric linear program and to use its optimality conditions.
Furthermore, if the active set in the DCS is fixed, we obtain locally a smooth ODE or Differential Algebraic Equation (DAE). 
We study the well-posedness of these systems.
In the theoretical analysis and algorithmic development, we focus on the Filippov systems.
\color{black}
Simple tutorial examples accompany all developments.
Most importantly, we present the extension of the FESD method to nonsmooth systems that are described with step functions. 
We also adapt the convergence and well-posedness results of~\cite{Nurkanovic2024a} to this case.
\color{black}
If Filippov systems are modeled via step functions, one ends up with multi-affine expressions, consisting of products of step functions, in the r.h.s. of the ODE.
We propose a lifting algorithm that introduces auxiliary variables and makes these expressions ``less nonlinear''. This can improve the convergence of the proposed method~\cite{Albersmeyer2010}.
\color{black}
The dynamical system \eqref{eq:basic_di} is more general than the Filippov system.
Therefore, the FESD method developed here applies to this broader class of systems.
The performance of the new method is compared to the original FESD~\cite{Nurkanovic2024a} and standard RK discretizations in terms of accuracy and computational time.
All methods are implemented in the open-source software package \nosnoc~\cite{Nurkanovic2022c}.

\paragraph*{Outline}
\color{black}
In Section \ref{sec:nonsmooth_systems_via_step}, we define and relate various classes of nonsmooth systems, including Dynamic Complementarity Systems (DCS), Filippov systems, and systems with set-valued Heaviside step functions. It is also demonstrated that the latter is equivalent to a DCS.
Section \ref{sec:dcs_properties} explores the properties of this DCS for a fixed active set and during active-set changes. 
In Section \ref{sec:fesd_step}, we introduce the FESD method for this class of problems. Section \ref{sec:fesd_step_theory} presents convergence and well-posedness results for the new method.
Additionally, in Section \ref{sec:compact_rep_and_lifting}, we demonstrate how to efficiently model piecewise smooth systems with step functions and introduce a lifting algorithm to reduce nonlinearity in the DCS.
\color{black}
In Section \ref{sec:numerical_examples}, we showcase the developments using several numerical examples.
\paragraph*{Notation}
The complementarity conditions for two vectors  $a,b \in \R^{n}$ \sloppy{read as ${0\leq a \perp b\geq 0}$, where $a \perp b$ means $a^{\top}b =0$}.
For two scalar variables $a,b$ the so-called C-functions \cite{Facchinei2003} have the property $\phi(a,b) = 0 \iff a\geq 0, b\geq 0, ab = 0$.
Examples are the natural residual functions $\phi_{\mathrm{NR}}(a,b)=\min(a,b)$ or the Fischer-Burmeister function $\phi_{\mathrm{FB}}(a,b) = a+b-\sqrt{a^2+b^2}$.
If $a,b \in \R^{n}$, we use $\phi(\cdot)$ component-wise and define $\Phi(a,b) = (\phi(a_1,b_1),\dots,\phi(a_{n},b_{n}))$.
All vector inequalities are to be understood element-wise, $\mathrm{diag}(x)\in \R^{n\times n}$ returns a diagonal matrix with $x \in \R^n$ containing the diagonal entries.
The concatenation of two column vectors $a\in \R^{n_a}$, $b\in \R^{n_b}$ is denoted by $(a,b)\coloneqq[a^\top,b^\top]^\top$, the concatenation of several column vectors is defined analogously.
A column vector with all ones is denoted by $e=(1,1,\dots,1) \in \R^n$ and its dimension is clear from the context.
The closure of a set $C$ is denoted by $\overline{C}$, its boundary as $ \partial C$.
\color{black}
For a given vector $a\in\R^n$ and set $\I \subseteq \{1,\ldots,n\}$, we define the projection matrix $P_\I \in \R^{|I| \times n}$, which has zeros or ones as entries. 
It selects all component $a_i, i\in \I$ from the vector $a$, i.e., $a_\I = P_\I a \in \R^{|\I|}$ and $a_\I = [a_i \mid i\in \I]$.
\color{black}
Given a matrix $M \in \R^{n \times m}$, its $i$-th row is denoted by $M_{i,\bullet}$ and its $j$-th column is denoted by $M_{\bullet,j}$.
For the left and right limits, we use the notation ${x(\ts^+)  = \lim\limits_{t\to \ts,\ t>\ts} x(t)}$ and ${x(\ts^-)  = \lim\limits_{t\to \ts,\  t<\ts}x(t)}$, respectively.

\renewcommand{\arraystretch}{1.2} 
\begin{table}
\centering
\caption{Key symbols used in this paper.}
\vspace{-0.35cm}
\label{tab:symbols_filippov}
\begin{tabular}{@{}ll@{}}
\specialrule{.15em}{1em}{0em} 
{Symbol} & {Description}\\
\hline
$x$ & differential state,  Eq.~\eqref{eq:pss}\\
$u$ & control function,  Eq.~\eqref{eq:pss}\\
$\theta$ & Filippov's convex multipliers\\
$S$ & sign matrix defining regions $R_i$,~\eqref{eq:sign_matrix}\\
$\alpha$ & selection of set-valued step function,~\eqref{eq:step_set}\\
$\lambdap$ & Lagrange multiplier in step DCS,~\eqref{eq:step_dcs}\\
$\lambdan$ & Lagrange multiplier in step DCS,~\eqref{eq:step_dcs}\\
$\beta$ & Lifting variable in the step DCS,  Sec.~\ref{sec:step_lifting_algortihm}\\
\hline
$\switchfun(x)$ & switching functions, Sec. \ref{sec:step_dcs}\\
$\gamma(\switchfun_j(x))$ &  set-valued Heaviside step function, Eq.~\eqref{eq:step_function}\\
$\Gamma(\switchfun(x))$ & concatenation of all regarded step functions , Eq.~\eqref{eq:step_function}\\
$f_i(x,u)$ & $i$-th mode of the total $\Nsys$ modes of the PSS system, Eq.~\eqref{eq:pss}\\
$F(x,u)$ & matrix collecting all PSS modes, F(x,u) = $[f_1(x,u),\dots, f_{\Nsys}(x,u)] \in \R^{n_x \times \Nsys}$,~\eqref{eq:step_dcs}\\ 
$G(x,\theta,\alpha,\lambdap,\lambdan)$ & algebraic eq. in the step DCS as nonsmooth DAE~\eqref{eq:dcs_step_lp}\\
$g(x)$ & Stewart's indicator function,~\eqref{eq:stewart_sets}\\
$W_{\mathcal{K},\mathcal{I}}(x,u)$ & auxiliary matrix used in the study of the DCS~\eqref{eq:step_dcs},  Sec.~\ref{sec:step_dcs_active_set_fixed}\\
$B_{\mathcal{K},\mathcal{I}}(x,u)$  & auxiliary matrix used in the study of the DCS~\eqref{eq:step_dcs},  Sec.~\ref{sec:step_dcs_active_set_fixed}\\
 \hline
$R_i$ & regions of the PSS, Eq.~\eqref{eq:pss}\\
$\tilde{R}_i$ & base sets, Def.~\ref{def:basis_sets}\\
$\mathcal{J}$ & index set of PSS modes, Eq.~\eqref{eq:pss}\\
\color{black}$\F(x,u,\Gamma(\switchfun(x)))$ & \color{black}r.h.s. of the nonsmooth ODE with set-valued Heaviside step functions, Eq.~\eqref{eq:basic_di}\\
\color{black}{$\Ffilippov(x,u)$ } & \color{black}r.h.s. of the Filippov differential inclusion, Eq.~\eqref{eq:filippov_di}, 
\eqref{eq:filippov_di_with_multiplers}\\
\color{black}{$\Fap(x,u)$ } &\color{black} r.h.s. of the Aizerman–Pyatnitskii differential inclusion, Eq.~\ref{eq:di_ap}\\
\color{black}{$\Fstep(x,u)$}  &\color{black} Filippov set obtained with set-valued Heaviside step functions, Eq.~\ref{eq:step_set}\\ 
$\I(\cdot)$ & active set for Filippov systems, Eq.~\eqref{eq:active_set_defintion}\\ 
$\mathcal{T}_{n}$  & $n-$th time interval with fixed active set $\mathcal{T}_{n}=(t_{\mathrm{s},n},t_{\mathrm{s},n+1})$ \\ 
$\I_{n}$  &  the fixed active set  $\mathcal{I}_{n} = \mathcal{I}(x(t)), t\in I_n$  \\ 
$\I_{n}^{0}$  & active set at $t_{\mathrm{s},n}$, i.e. $ \mathcal{I}_{n,0} = \mathcal{I}(t_{\mathrm{s},n})$ \\ 
$\C$ & index set of all switching functions \textcolor{black}{$\switchfun_i(x)$} \\ 
$\K$  & index set of all switching functions \textcolor{black}{$\switchfun_i(x)$} that are zero for a given $\I$, Sec.~\ref{sec:step_dcs_active_set_fixed}\\
\specialrule{.15em}{0.5em}{0em}
\end{tabular}
\end{table}

\section{Overview and relations between different classes of nonsmooth dynamical systems}\label{sec:nonsmooth_systems_via_step}
This section defines, reviews, and relates several classes of nonsmooth dynamical systems relevant to the algorithmic development in this paper.
The class of Dynamic Complementarity Systems (DCS) is defined in Section \ref{sec:dcs}.
In Section \ref{sec:DRHS_and_filippov}, generic ODEs with a Discontinuous Right-Hand Side (DRHS) and their Filippov extension are considered.
In Section \ref{sec:pss_and_filippov}, these results are applied to a more structured ODE with DRHS, namely to piecewise smooth systems.
This is followed by Section~\ref{sec:ap_inclusion}, where another class of structured ODEs with DRHS, i.e., nonsmooth systems where Heaviside step functions appear in the dynamics. 
This class is the focus of this paper, but it is closely related to other classes that we review.
Section \ref{sec:step_functions} relates these systems to Filippov and DCSs.
The key tool is to rewrite the set-valued Heaviside step function in terms of the optimality conditions of an appropriate linear program. 
For comparison, we review Stewart's way of rewriting a Filippov system as a DCS in Section \ref{sec:stewarts_representation}.
We conclude with Section \ref{sec:nonsmooth_systems_summary}, where we summarize the relations between all system classes regarded in this section.

This section assumes that a continuous control function $u: [0,T]\to \R^{n_u}$ is given. 
Discontinuous $u(t)$ can be considered by partitioning the considered time interval $[0, T]$ into pieces where $u(t)$ is locally continuous, and considering a sequence of ODEs with continuous control functions.
\color{black}
\subsection{Dynamic complementarity systems}\label{sec:dcs}
A dynamic complementarity system (DCS) \cite{Brogliato2020,Pang2008} is the problem:
\begin{subequations}\label{eq:dcs_general}
	\begin{align}
		\dot{x}(t) &= f(x(t),y(t),z(t),u(t)), \label{eq:dcs_general_ode}\\
		0&= \gdcseq(x(t),y(t),z(t)), \label{eq:dcs_general_alg}\\
		0 &\leq y(t) \perp \gdcscomp(x(t),y(t),z(t)) \label{eq:dcs_general_cc}\ \text{for almost all } t,
	\end{align}
\end{subequations}
where $x \in \R^{n_x}$  are the differential states, $y \in  \R^{n_y}$ and $z \in \R^{n_z}$ are the algebraic states of the DCS.
The DCS consists of an ODE (defined the function $f: \R^{n_x}\times \R^{n_y} \times \R^{n_z}\times \R^{n_u} \to \R^{n_x}$) coupled with algebraic equality constraints~\eqref{eq:dcs_general_alg} (defined by $\gdcseq:\R^{n_x}\times \R^{n_y} \times \R^{n_z} \to \R^{n_z}$) and complementarity constraints~\eqref{eq:dcs_general_cc} (defined by $y:[0,T] \to \R^{n_y}$ and $\gdcscomp:\R^{n_x}\times \R^{n_y} \times \R^{n_z} \to \R^{n_y}$).
We assume that the problem functions $f,\gdcseq$ and $\gdcscomp$ are at least twice continuously differentiable in all arguments.
In the complementarity conditions \eqref{eq:dcs_general_cc}, both components $y_i$ and $g_{\mathrm{c}, i}(x,y,z)$ have to be non-negative for $i = 1,\ldots,n_y$, but only one of them is allowed to be strictly positive at a time $t$.
These conditions encode combinatorial structure in the problem and make it nonsmooth.
In contrast to many other nonsmooth system formulations, complementarity conditions can be efficiently treated with Newton-type methods, which makes them attractive from a computational point of view~\cite{Facchinei2003}. 

The formulation in Eq. \eqref{eq:dcs_general} is very general and further assumptions on $f, \gdcseq$ and $\gdcscomp$ determine the existence, uniqueness, and qualitative properties of $x$, $y$ and $z$.
Several important models fit into the form of \eqref{eq:dcs_general}, most prominently rigid bodies with friction and impact \cite{Brogliato2020,Stewart2011} and electric circuits with electronic devices \cite{Acary2010}. 
We refer the reader to \cite{Brogliato2020,Pang2008,Stewart2011} for an extensive collection of theoretical results and 
application examples.

If one writes \eqref{eq:dcs_general_cc} equivalently via a C-function $\Phi$, i.e., $\Phi(y,\gdcscomp(x,y,z)) = 0$, the DCS \eqref{eq:dcs_general} can be seen as a nonsmooth Differential Algebraic Equation (DAE).
The DCS is an instance of a more general class of problems called Differential Variational Inequalities (DVI), formally introduced by Pang and Stewart~\cite{Pang2008}. 
One way to classify DCS (and DVIs) is their \textit{index} \cite{Pang2008}, that is, how many times $\gdcscomp(x,y,z)=0$ must be differentiated with respect to time $t$ to find $z$ and $y$ as a function of $x$. 
A similar concept is the \textit{relative degree}, cf. \cite[Appendix C]{Brogliato2020}.
Index zero DVIs (and DCS) have continuously differentiable solutions~\cite[Proposition 5.1]{Pang2008}.
The existence and uniqueness of solutions of index zero systems are proven in \cite{Pang2008}.
Index one DVIs have in general absolutely continuous solutions~\cite{Stewart2011}.
Uniqueness for some index one DVIs is proven by Stewart \cite{Stewart2008}.
Index two systems are systems with state jumps \cite[Chapter 6]{Stewart2011}.
For index two DVIs only existence can be proven \cite{Stewart1997}. 
In general, for higher index DVIs solutions fail to exist \cite{Stewart2011}, and if they exist they are distributions~\cite{Brogliato2020}. 
In this paper, we are interested in DCS where $x(t)$ and $y(t)$ are continuous functions of time, and $z(t)$ is allowed to be discontinuous.


\subsection{Discontinuous ODEs and Filippov systems}\label{sec:DRHS_and_filippov}
Consider an ODE with a discontinuous right-hand side (DRHS):
\begin{align}\label{eq:discontinuous_ode}
	\dot{x}(t) = f(x(t),u(t)),
\end{align}
with a given initial value $x(0) = x_0$.
More precisely, $f$ is discontinuous in $x$ but continuous in $u$.

Classical solutions $x(t)$ to \eqref{eq:discontinuous_ode} on an interval $[0,T]$ are continuously differentiable, which is impossible with a discontinuous $f$.
Carath\'eodory solutions of an ODE satisfy the ODE in integral form, i.e., $x(t) = x(0) + \int_{0}^{t} f(x(t),u(t))\dd t$ for $t >0$, where the integral is the Lebesgue integral. 
This allows discontinuities of $f$ in $x$ and in some situations, a Carath\'eodory solution might be of use, as $x$ is absolutely continuous, and it requires \eqref{eq:discontinuous_ode} to be satisfied for almost all $t\in [0,T]$. 
However, this relaxation is not always sufficient as can be seen from the following example~\cite[Example 5]{Cortes2008a}:
\begin{align}\label{ex:cortes}
	\dot{x} &=  \begin{cases}
		1, & x<0,\\
		-1, & x\geq0,
	\end{cases}
\end{align}
with $x(0) = x_0$. 
For $x_0 > 0$, there exist the solution $x(t) = x(0)-t$ for $t\in \left[0, x_0 \right)$. 
Similarly, for For $x_0 < 0$, there exist the solution $x(t) = x(0)+t$ for $t\in \left[0, -x_0 \right)$.
For $t$ larger than $|x_0|$ in both cases, each solution reaches the point $x(t) = 0$ and cannot leave it as the vector fields from both sides push towards it. 
However, since $\dot{x} = 0 \neq -1$, we have no solution in the classical or Carath\'eodory sense, and a more general concept is necessary.

Deciding on which more general notion of solution to use is a modeling step.
A popular option in the control community is Fillipov's extension, which suggests embedding the ODE \eqref{eq:discontinuous_ode} into the following differential inclusion~\cite{Filippov1964}:
\begin{align}\label{eq:filippov_di}
	\dot{x}(t) \in \Ffilippov(x(t),u(t))\coloneqq \bigcap_{\delta>0} \bigcap_{\mu(N) = 0} \overline{\mathrm{conv}} f(x(t)+\delta \mathcal{B}(x(t))\setminus N,u(t)).
\end{align}
The right-hand side is now a closed convex set.
Given a point $x(t)$ and the corresponding $u(t)$, the idea behind this definition is to regard the closed convex hull of all neighboring values in a ball $x(t)\in \delta B(x(t))$ instead of only $f(x(t),u(t))$ and thereby ignoring all values of $f(x(t),u(t))$ on sets of measure zero. 
In practice, these sets are surfaces on which $f$ becomes discontinuous or points/surfaces where $f$ is not even defined. 
If $f(\cdot)$ is continuous at $x(t)$ we obtain $\Ffilippov(x(t),u(t)) = \{f(x(t),u(t))\}$. 
This DI \eqref{eq:filippov_di} is a so-called Filippov system.
An absolutely continuous function ${x}$ is said to be a Filippov solution if it is a solution to the initial value problem defined \eqref{eq:filippov_di}.
It is guaranteed that at least one solution exists to the Filippov DI \eqref{eq:filippov_di}~\cite{Cortes2008a,Filippov1988}.
For the uniqueness of solutions, more assumptions must be made, for explicit statements see~\cite{Brogliato2020,Cortes2008a,Filippov1988}.

For illustration, we revisit the example in \eqref{ex:cortes} and apply Filippov's extension to it. 
We obtain
\begin{align}\label{ex:cortes_filippov}
	\dot{x} &\in  \begin{cases}
		\{1\}, & x<0,\\
		[-1,1], & x=0,\\
		\{-1\}, & x\geq0,
	\end{cases}
\end{align}
It can be seen that this DI has a solution for any $x(0)$ and for all $t \in \left[0,\infty\right)$, since for $t>|x_0|$ we have that $\dot{x} = 0 \in [-1,1]$.

\subsection{Piecewise smooth systems as Filippov systems}\label{sec:pss_and_filippov}
In the modeling of physical systems and control applications the discontinuities in \eqref{eq:discontinuous_ode} usually appear in a very structured way.
Such an example are piecewise smooth systems, which read as:
\color{black}
\begin{align} \label{eq:pss}
	\dot{x}(t) = f_i(x(t),u(t)),\ &\mathrm{if} \  x(t) \in R_i \subset \R^{n_x}, \ i \in \mathcal{J} \coloneqq \{ 1,\dots,\Nsys  \},
\end{align}
where $R_i$ are disjoint, connected, and open sets.
They are assumed to be nonempty and to have piecewise-smooth boundaries $\partial R_i$.
We assume that $\overline{\bigcup\limits_{i\in \mathcal{J}} R_i} = \R^{n_x}$ and that $\R^{n_x} \setminus \bigcup\limits_{i\in\mathcal{J}} R_i$ is a set of measure zero.
The functions $f_i(\cdot)$ are assumed to be at least twice continuously differentiable and Lipschitz continuous on an open neighborhood of $\overline{R}_i$.
\color{black}
We regard PSS where the regions $R_i$ are defined via a finite number of switching functions $\switchfun_j(x)$, $j = 1,\ldots,\nswitchfun$. 
It is assumed that the functions $\switchfun_j(x)$ are Lipschitz continuous and at least twice continuously differentiable.
\color{black}
\begin{definition}[Base regions]\label{def:basis_sets}
	Given $\nswitchfun$ scalar switching functions $\switchfun_j(x),\ j \in \mathcal{C} \coloneqq \{1,\ldots,\nswitchfun\}$, we define the \textcolor{black}{$n_f = 2^{\nswitchfun}$} base regions:
	\begin{align*}
		\tilde{R}_1 &= \{x \in \R^{n_x}  \mid \switchfun_1(x) >0,\switchfun_2(x) >0, \ldots, \switchfun_{\nswitchfun-1}(x) >0 , \switchfun_{\nswitchfun}(x) >0 \},\\
		\tilde{R}_2 &= \{x \in \R^{n_x}  \mid \switchfun_1(x) >0,\switchfun_2(x) >0, \ldots, \switchfun_{\nswitchfun-1}(x) >0 , \switchfun_{\nswitchfun}(x) <0 \},\\
		\vdots\\
		\tilde{R}_{n_f} &= \{x \in \R^{n_x}  \mid \switchfun_1(x)<0,\switchfun_2(x) <0, \ldots, \switchfun_{\nswitchfun-1}(x) <0 , \switchfun_{\nswitchfun}(x) <0 \},
	\end{align*}
	such that $\R^{n_x} = \cup_{i=1}^{n_f} \tilde{R}_i$.
	These definitions are compactly expressed via a dense sign matrix $S \in \R^{n_f \times \nswitchfun}$:
	\begin{align}\label{eq:sign_matrix}
		S = \begin{bmatrix}
			1 & 1 & \dots &1 & 1\\
			1 & 1 & \dots & 1& -1\\
			\vdots &  \vdots & \cdots &  \vdots &  \vdots\\
			-1 & -1 & \dots & -1& -1\\
		\end{bmatrix}.
	\end{align}
	The matrix $S$ has no repeating rows and no zero entries. 
	The sets $\tilde{R}_i$ are defined using the rows of the matrix $S$:
	\begin{align}\label{eq:standard_sets_step_dense}
		\tilde{R}_i &= \{ x\in \R^{n_x} \mid S_{i,j} \switchfun_j(x)>0, j \in \mathcal{C} \},\ i = 1,\ldots, n_f.
	\end{align}
\end{definition}
\color{black}
Note that the boundaries of the regions $\partial \tilde{R}_i$ \textcolor{black}{are unions of subsets} of the zero-level sets of corresponding functions $\switchfun_j(x)$.
\color{black}
For notational convenience and ease of exposition, we assume that the PSS regions are equal to the base regions, i.e., $n_f = 2^{\nswitchfun}$ and $R_i = \tilde{R_i}$ for $i = 1,\ldots,n_f$. 
Other cases are discussed later in Section \ref{sec:compact_rep_and_lifting}.

Note that the ODE \eqref{eq:pss} is not defined on the region boundaries $\partial R_i$. 
However, if the vector fields from both sides of $\partial R_i$ push towards it, the solution needs to evolve on $\partial R_i$.
The cases when $x(t)$ must evolve on $\partial R_i$ \cite{Cortes2008a,Filippov1988} are so-called \textit{sliding modes}. 
To define the sliding mode dynamics, one can apply Filippov's extension to the PSS~\eqref{eq:pss}.
Furthermore, the special structure of the PSS allows a more explicit definition of \eqref{eq:filippov_di} via a finite number of multipliers $\theta_i$~\cite{Filippov1988,Stewart1990a}. 
The corresponding DI reads as
\color{black}
\begin{align}\label{eq:filippov_di_with_multiplers}
	\dot{x}(t) \in \Ffilippov(x,u) = \Big\{ &\sum_{i\in \mathcal{J}}
	f_i(x,u) \, \theta_i  \mid \sum_{i\in \mathcal{J}}\theta_i = 1, \theta_i \geq 0,
	\ \theta_i = 0 \  \mathrm{if} \;  x \notin \overline{R}_i,
	\forall  i  \in \mathcal{J} \Big\}. 
\end{align}
An important notion for Filppov PSS \eqref{eq:filippov_di_with_multiplers} is the active set, which is defined as the set
\begin{align}\label{eq:active_set_defintion}
	\mathcal{I}(x(t)) = \{ i \in \mathcal{J} \mid \theta_i(t) >0 \}.
\end{align}
Note that if $x(t)$ is in the interior of a region, then $\I$ is a singleton, and for sliding modes, it consists of the indices of regions neighboring the current sliding surface.
\color{black}
\begin{remark}\label{rem:time_freezing}
	A solution of the Filippov PSS~\eqref{eq:filippov_di_with_multiplers} is an absolutely continuous function.
	However, systems with state jumps (and thus discontinuous $x(t)$) are often of practical interest.
	One way to transform such systems into a PSS, and thus to apply the methods developed in this paper, is to use the time-freezing reformulation~\cite{Halm2023,Halm2023a,Nurkanovic2022a,Nurkanovic2021,VanRoy2023}.
	In several classes of systems with state jumps some parts of the state space are forbidden (e.g., a hooping robot is not allowed to penetrate the ground). 
	The main idea in the time-freezing is to define auxiliary regions $R_i$ in the forbidden region and corresponding auxiliary dynamics $f_i(x)$. 
	By construction, the endpoints of the trajectory parts evolving in the auxiliary regions satisfy the state jump law.
	Furthermore, a clock state $t(\tau)$ is introduced, which evolves in the feasible region of the initial system ($\frac{\dd t}{\dd \tau} >0$) and is frozen in the initially forbidden region ($\frac{\dd t}{\dd \tau} =0$). 
	Here $\tau$ is the time of the new relaxed system. 
	Now, by taking only the parts of the trajectory where the clock state $t(\tau)$ evolved, we recover the discontinuous solution of the original system.
	An example of such a system is treated in Section \ref{sec:robot_ocp}.
\end{remark}

\subsection{Aizerman–Pyatnitskii differential inclusions}\label{sec:ap_inclusion}
So far we have introduced dynamic complementarity and Filippov systems. 
Now we introduce another special case of the discontinuous ODE~\eqref{eq:discontinuous_ode} and relate it to the others later with the help of step functions.
Regard and ODE with DRHS of the form of:
\begin{align}\label{eq:ap_ode}
	\dot{x} = f(x,u,v(x)),
\end{align}
where $f:\R^{n_x} \times \R^{n_u} \in \R^{n_v}$ is continuous in all arguments, but $v(x)$ is discontinuous, e.g. it could be the usual single-valued Heaviside step function.
Such systems appear commonly in sliding mode control or the modeling of gene-regulatory networks~\cite{Acary2014}.
It is assumed that at each point of discontinuity, for the components $v_i(x)$ a closed convex set $V_i(x)$ is given, out of which the corresponding arguments of \eqref{eq:ap_ode} take their values. 
Hence, one can define the following differential inclusion:
\begin{align}\label{eq:di_ap}
	\dot{x} \in \mathcal{F}_{\mathrm{AP}}(x,u) \coloneqq \Big\{f(x,u,v)\mid v_i \in V_i(x), i =1,\ldots, n_v \Big\}.
\end{align}
The set $\mathcal{F}_{\mathrm{AP}}(x,u)$ is in general nonconvex, except in some special cases. 
For example, $\mathcal{F}_{\mathrm{AP}}(x,u)$ is convex if $v(x)$ enters the r.h.s. of \eqref{eq:ap_ode} linearly and $V(x)$ is closed convex.
The DI \eqref{eq:di_ap} does not have as rich a theory as Filippov DIs, and there are fewer results available on the existence of solutions and convergence of numerical methods~\cite{Acary2014}.
These systems are called Aizerman–Pyatnitskii DIs, cf. \cite{Acary2014} and \cite[page 55, Definition c]{Filippov1988}.
Time-stepping methods for such systems were developed in~\cite{Acary2014}.

In this paper, we focus on the case where the functions $V_i(x)$ are given by set-valued Heaviside step functions. 
In particular, we regard the DI \eqref{eq:basic_di} introduced at the beginning of the paper, where $V(x) = \Gamma(\switchfun(x))$.

\subsection{Rewriting Filippov PSS as DCS via set-valued Heaviside step functions}
\label{sec:step_functions}
\color{black}
The next question we address is: How can we find an (implicit) function of $x$ for computing the Filippov multipliers $\theta$ in \eqref{eq:filippov_di_with_multiplers}? 
To derive such a function, we will make use of the set-valued step function and the algebraic representations of the regions $R_i$ in \eqref{eq:standard_sets_step_dense}.
Moreover, expressing the step function via the Karush–Kuhn–Tucker (KKT) conditions of a suitable Linear Program (LP), we can write the Filippov DI \eqref{eq:filippov_di_with_multiplers} as an equivalent DCS. 
This DCS will be the main formulation used in the development of the FESD method. 

\subsubsection{Heaviside step functions via linear programming}
We start the described procedure with a closer look at the step functions.
\color{black}
Let us denote by $\alpha \in \R^{\nswitchfun}$ a selection $\alpha \in \Gamma(\switchfun(x))$.
A well-known way to express the function $\Gamma(\switchfun(x))$~\cite{Acary2014,Baumrucker2009} is the use of the solution map of the parametric linear program:
\begin{subequations}\label{eq:step_parametric_lp}
	\begin{align}
		\Gamma(\switchfun(x))  = \arg \min_{\alpha\in \R^{\nswitchfun}} \ &-\switchfun(x)^\top \alpha \\
		\quad \mathrm{s.t.} \quad & 0 \leq \alpha_i \leq 1, \ i=1,\ldots,\nswitchfun \label{eq:step_parametric_lp_ineq}.
	\end{align}
\end{subequations}
Note that all components of $\alpha$ are decoupled in this LP, i.e., every $\alpha_i$ can be obtained by solving a one-dimensional LP with the objective $-\switchfun_i(x)\alpha_i$ and the feasible set $0\leq \alpha_i \leq 1$.
Let $\lambdan,\lambdap \in \R^{\nswitchfun}$ be the Lagrange multipliers for the lower and upper bound on $\alpha$ in \eqref{eq:step_parametric_lp_ineq}, respectively.
The KKT conditions of \eqref{eq:step_parametric_lp} read as
\begin{subequations}\label{eq:step_parametric_lp_kkt}
	\begin{align}
		&\switchfun(x) = \lambdap - \lambdan, \label{eq:step_parametric_lp_kkt_lagrangian}\\
		&0 \leq  \lambdan \perp \alpha \geq 0,\\
		&0 \leq  \lambdap \perp e-\alpha \geq 0,
	\end{align}
\end{subequations}
\textcolor{black}{Now we have a purely algebraic representation of the set-valued Heaviside step function.}
Let us look at a single component $\alpha_j$ and the associated functions $\switchfun_j(x)$.
From the LP \eqref{eq:step_parametric_lp} and its KKT conditions, one can see that for $\switchfun_j(x) > 0$, we have $\alpha_j=1$.
Since the upper bound is active, we have that $\lambdan_j = 0$ and from \eqref{eq:step_parametric_lp_kkt_lagrangian} that $\lambda_{\mathrm{p},j} = \switchfun_j(x) > 0$.
Likewise, for $\switchfun_j(x) < 0$, we obtain $\alpha_j=0$, $\lambdap_j = 0$ and $\lambdan_j = -\switchfun_j(x) > 0$.
On the other hand, $\switchfun_j(x) = 0$ implies that $\alpha_j \in [0,1]$ and $\lambdap_j = \lambdan_j = 0$.
From this discussion, it can be seen that $\switchfun(x)$, $\lambdan$ and $\lambdap$ are related by the following expressions:
\begin{align}\label{eq:step_continuity_lambda}
	\lambdap &= \max(\switchfun(x),0),\	\lambdan = -\min(\switchfun(x),0).
\end{align}
That is, $\lambdap$ collects the positive parts of $\switchfun(x)$ and $\lambdan$ the absolute value of the negative parts of $\switchfun(x)$.
From this relation, we can immediately conclude the following:
\begin{lemma}\label{lem:cont_lambda}
	Let $\switchfun(x(t))$ be a continuous function of time, and then the functions $\lambdap(t)$ and $\lambdan(t)$ are continuous in time.
\end{lemma}
\color{black}
The continuity of the Lagrange multipliers $\lambdap(t)$ and $\lambdan(t)$ plays a crucial role in the switch detection in the FESD method. 
The exact switch detection is necessary for high-order accuracy and the correct computation of numerical sensitivities.
\color{black}
\subsubsection{Filippov system as a dynamic complementarity systems}\label{sec:step_dcs}
\color{black}
For the sake of clarity, we start by illustrating what we want to achieve on a simple example and give in the sequel the general expression.	
\color{black}
\begin{example}[Step representation]\label{ex:sparse_sets}
We regard a PSS with four regions defined via two scalar switching functions $\switchfun_1(x)$ and $\switchfun_2(x)$.
The regions are equal to the base sets from Definition \ref{def:basis_sets} and read as
${R}_1 = \{ x\in \R^{n_x} \mid \switchfun_1(x) > 0, \switchfun_2(x)>0\},\; 
{R}_2 = \{ x\in \R^{n_x} \mid \switchfun_1(x) > 0, \switchfun_2(x)<0\},\;
{R}_3 = \{ x\in \R^{n_x} \mid \switchfun_1(x) < 0, \switchfun_2(x)>0\}$ and
${R}_4 = \{ x\in \R^{n_x} \mid \switchfun_1(x) < 0, \switchfun_2(x)<0\}$.
\color{black}
The corresponding sing matrix $S \in \R^{4\times 2}$ read as
\color{black}
\begin{align*}
	{S} = \begin{bmatrix}
			1& 1\\
			1& -1\\
			-1& 1\\
			-1&-1
 		\end{bmatrix}.
\end{align*}
The corresponding Filippov system (defined by \eqref{eq:filippov_di_with_multiplers} reads as:
\color{black}
\begin{align}\label{eq:filippov_example}
	\dot{x} \in \{ \sum_{i=1}^{4}\theta_i f_i(x) \mid \theta \geq 0, \sum_{i=1}^{4}\theta_i = 1,\ \theta_i =0, \textit{if} \ x\notin \overline{R_i} \}.
\end{align}
Let $x \in R_1$, then $\alpha_1 \in \gamma(\psi_1(x)) = \{1\}$ and $\alpha_2 \in \gamma(\psi_2(x)) = \{1\}$, thus $\theta_1 = \alpha_1\alpha_2 = 1$. On the other hand, by direct evaluation of the step functions, one can see that $\theta_1 = \alpha_1\alpha_2 = 0$ if $x \notin \overline{R_1}$, since at least one of the selections $\alpha_1$ or $\alpha_2$ is zero. 
Similarly, for $x \in R_2$, we observe that $\alpha_1 \in \gamma(\psi_1(x)) = \{1\}$ and $\alpha_2 \in \gamma(\psi_2(x)) = \{0\}$ and we can set $\theta_2 = \alpha_1 (1-\alpha_2) =1$. 
By direct evaluation of we can see that $\theta_2 > 0$ if $x\in \overline{R_2}$ and $\theta_2 = 0$, otherwise.
Following the same pattern, we conclude that $\theta_3 = (1-\alpha_1)\alpha_2$ and $\theta_4 = (1-\alpha_1)(1-\alpha_2)$. 
Thus, we can define the system
\begin{align}\label{eq:ap_example}
	\begin{split}
		\dot{x} \in \Big\{ &\alpha_1 \alpha_2 f_1(x) + 
						\alpha_1 (1-\alpha_2)f_2(x) + 
						(1-\alpha_1)\alpha_2 f_3(x) + 
						(1-\alpha_1)(1-\alpha_2) f_4(x) \\
						&\mid \alpha_1 \in \gamma(\switchfun_1(x)), \alpha_2 \in \gamma(\switchfun_2(x)) \Big\}.
		\end{split}
\end{align}
Since $\alpha_1, \alpha_2 \in [0,1]$ it is clear that $\theta_i \in [0,1], i\in \{1,\ldots,4\}$.
Moreover, direct calculation shows that $\sum_{i=1}^{4} \theta_i = 1$.
Therefore, we conclude that the sets in the r.h.s. of \eqref{eq:filippov_example} and \eqref{eq:ap_example} are the same sets, i.e., we can express the Filippov set via set-valued Heaviside step functions.

Furthermore, we can observe how the sign pattern in ${S}$ determines how $\alpha_j$ enters the expression for $\theta_i$.
For ${S}_{i,j} =1$ we have $\alpha_j$, for ${S}_{i,j} =-1$ we have $(1-\alpha_j)$.
In summary, the definition of $\theta_i$ consists of products of of $\alpha_j$ and $(1-\alpha_k)$, i.e., it is multi-affine in the selections $\alpha_j$, $j = 1,\ldots,\nswitchfun$.
\color{black}
\end{example}
We generalize the patterns observed in the previous example and define the set
	\begin{align}\label{eq:step_set}
		\mathcal{F}_{\mathrm{H}}(x,u) &\coloneqq
		\Big\{ \sum_{i=1}^{\textcolor{black}{n_f}} \prod_{j =1}^{\nswitchfun} 	\Big( \frac{1-{S}_{i,j}}{2}+{S}_{i,j}\alpha_i \Big) f_i(x,u) \mid \alpha \in \Gamma(\switchfun(x)) \Big\}.
	\end{align}
	Note that we have
	\begin{align*}
		\frac{1-{S}_{i,j}}{2}+{S}_{i,j}\alpha_i	&= \begin{cases}
			\alpha_i,\ & \textrm{ if } {S}_{i,j} = 1, \\
			1-\alpha_i,\ & \textrm{ if }  {S}_{i,j} = -1.
		\end{cases}
	\end{align*}
	Similar definitions of $\Fstep(x,u)$ as in \eqref{eq:step_set} can be found in \cite[Section 4.2]{Dieci2011} and \cite[Section 2.1]{Guglielmi2022}.
	Observe that this set has the same form as the r.h.s. $\Fap(x,u)$ in \eqref{eq:di_ap}.
	Next, we show that $\Fstep(x,u)$ is indeed the same set as $\Ffilippov(x,u)$, i.e., the set in the r.h.s. of \eqref{eq:filippov_di_with_multiplers}.
	
		\begin{lemma}[Lemma 1.5 in \cite{Dieci2011}]\label{lem:dieci}
			Let $a_1,a_2,\ldots, a_m \in \R$. 
			Consider the $2^m$ non-repeated products of the form
			$p_i = (1\pm a_1)(1\pm a_2)\cdots(1\pm a_m)$, then it holds that $\sum_{i=1}^{2^m} p_i = 2^m$.
		\end{lemma}
		\begin{proposition}\label{prop:filippov_ap}
			Let 
			\begin{align}\label{eq:theta_via_step}
				\theta_i &= \prod_{j =1}^{\nswitchfun}   \Big(\frac{1-{S}_{i,j}}{2}+{S}_{i,j}\alpha_j \Big),\ \textrm{for all}\ i\in \mathcal{J} = \{1,\ldots,\Nsys\},
			\end{align}
			then it holds that $\Ffilippov(x,u) =\Fstep(x,u)$.
		\end{proposition}
		\textit{Proof.}
		We only need to show that $\theta_i \geq 0$ for all $i \in \J$ and $\sum_{i\in \J} \theta_i  =1$.
		It is easy to see that $\theta_i\in [0,1]$ since it consists of a product of terms that takes value in $[0,1]$.
		
		\color{black}
		Without loss of generality, regard $\theta_1$ and suppose that $x \notin \overline{R_1}$.
		This means that $x\in R_i, i \neq 1$ and that at least one $\switchfun_j(x) <0, j\in\mathcal{C}$, which implies that $\alpha_j = 0$. From \eqref{eq:theta_via_step} it follows that $\theta_1 = 0$ if $x \notin \overline{R_1}$.
		By similar arguments it follows that $\theta_i = 0$ if $x \notin \overline{R_i}$ for $i = 1,\ldots,n_f$.
		\color{black}

		Next we show that $\sum_{i\in \mathcal{J}} \theta_i  =1$.
		We introduce the change of variables:
		\begin{align*}
			\frac{1+b_j}{2} &= \alpha_j,\; \frac{1-b_j}{2}= 1-\alpha_j.
		\end{align*}
		Then all $\theta_i$ are of the form
		\begin{align*}
			\theta_i = 2^{-\nswitchfun} \prod_{j =1}^{\nswitchfun} (1\pm b_j).
		\end{align*}

		By applying Lemma \ref{lem:dieci}, we conclude that $\sum_{i\in \mathcal{J}}\theta_i=1$ and the proof is complete.
		\qed
		
		To pass from the definition in Eq. \eqref{eq:step_set} to a dynamic complementarity system, we state the KKT conditions of \eqref{eq:step_parametric_lp} to obtain an algebraic expression for $\Gamma(\switchfun(x))$.
		Combining this with the definition of the Filippov set in \eqref{eq:step_set} and the expression for $\theta_i$ in \eqref{eq:theta_via_step}, we obtain the following DCS:
			\begin{subequations}\label{eq:step_dcs}
			\begin{align}
				&\dot{x} = F(x,u)\; \theta, \\
				&0=\theta_i - \prod_{j=1}^{\nswitchfun} \Big(  \frac{1-{S}_{i,j}}{2}+{S}_{i,j}\alpha_j\Big),\; \textrm{for all } i\in \J, \label{eq:step_dcs_theta}\\
				&0=\switchfun(x) - \lambdap + \lambdan,\label{eq:step_dcs_lp_lagrangian}\\
				&0 \leq  \lambdan \perp \alpha \geq 0, \label{eq:step_dcs_cc1}\\
				&0 \leq  \lambdap \perp e-\alpha \geq 0, \label{eq:step_dcs_cc2}
			\end{align}
		\end{subequations}
		where $F(x,u) = [f_1(x,u),\dots, f_{\Nsys}(x,u)] \in \R^{n_x \times \Nsys}$, $\theta =(\theta_1,\dots,\theta_{\Nsys}) \in \R^{\Nsys}$ and $\lambdap,\lambdan,\alpha \in \R^{\nswitchfun}$.
		We group all algebraic equations into a single function and use a C-function $\Psi(\cdot,\cdot)$ for the complementarity condition to obtain a more compact expression:
		\begin{align}\label{eq:dcs_step_lp}
			G(x,\theta,\alpha,\lambdap,\lambdan) \coloneqq
			\begin{bmatrix}
				\theta_1 - \prod_{j=1}^{\nswitchfun}   \Big(\frac{1-{S}_{1,j}}{2}+{S}_{1,j}\alpha_j\Big)\\
				\vdots\\
				\theta_{n_f} - \prod_{j=1}^{\nswitchfun}   \Big(\frac{1-{S}_{n_f,j}}{2}+{S}_{n_f,j}\alpha_j\Big)\\
				\switchfun(x) - \lambdap + \lambdan\\
				\Psi( \lambdan,\alpha)\\
				\Psi( \lambdap,e-\alpha)
			\end{bmatrix}.
		\end{align}
		Finally, we obtain a compact representation of \eqref{eq:step_dcs} in the form of a nonsmooth DAE:
		\begin{subequations}\label{eq:dcs_step_2}
			\begin{align}
				\dot{x} & = F(x,u)\theta,\\
				0&= G(x,\theta,\alpha,\lambdap,\lambdan) \label{eq:dcs_step_2_alg}.
			\end{align}
		\end{subequations}
		
		\color{black}
		We can see that \eqref{eq:step_dcs} is an instance of the generic DCS \eqref{eq:dcs_general}.
		Set $z = (\theta,\alpha)$ and $y = (\lambdan,\lambdap)$. 
		It follows that $f(x,y,z)  = F(x,u) \theta$, the function $\gdcseq(x,y,z)$ is define by the expressions in \eqref{eq:step_dcs_theta} and \eqref{eq:step_dcs_lp_lagrangian}, and $\gdcscomp(x,z,y)  = (\alpha, e-\alpha)$.
		
		\begin{remark}
			In this section, we have focused on PSS, their Filippov extension, and their relation to Aizerman-Pyatnitskii DIs \eqref{eq:di_ap}.
			However, we can also reformulate generic Aizerman-Pyatnitskii DIs with Heaviside step functions \eqref{eq:basic_di} into DCS by replacing the step function by~\eqref{eq:step_parametric_lp_kkt}.
			There is no need for equivalence to Filippov systems to apply this procedure, and the numerical method developed in this paper can also be applied directly to such DCS. 
			Such a more general example is treated in Section \ref{sec:numerical_simulation_gen_reg}.
		\end{remark}
		\color{black}
		
		\subsection{Stewart's representation}\label{sec:stewarts_representation}
		\color{black}
		The FESD method~\cite{Nurkanovic2024a} was originally developed for Stewart's representation~\cite{Stewart1990a} of Filippov systems.
		Stewart's representation assumes a specific definition of the regions $R_i$ and uses a LP to transform the Filippov system into an equivalent DCS, which is not the same DCS as \eqref{eq:step_dcs}.
		For comparison and completeness, we briefly introduced Stewart's DCS.
		\color{black}
		
		In Stewart's representation~\cite{Stewart1990a}, the regions $R_i$ are defined via so-called indicator functions $g_i(x)$ for all $i \in \mathcal{J} = \{1,\ldots,\Nsys\}$. 
		The definition reads as 
		\begin{align}\label{eq:stewart_sets}
			R_i = \{ x \in \R^{n_x} \mid g_i(x) < \min_{j\in \mathcal{J}\setminus\{i\}} g_j(x)\}.
		\end{align}
		\color{black}
		Arguably, this definition of the regions might be less intuitive than \eqref{eq:standard_sets_step_dense}.
		\color{black}
		However, if the regions $R_i$ match the base sets $\tilde{R}_i$ from Definition \ref{def:basis_sets}, 
		it was shown in \cite[Proposition 2]{Nurkanovic2024a} that the function $g(x) = (g_1(x),\ldots,g_{\Nsys}(x))$ can be obtained as:
		\begin{align}\label{eq:stewart_indicators}
			g(x) = -S \switchfun(x).
		\end{align}
		The multiplier vector $\theta$ is expressed as the solution of an LP parameterized by $x$:
		\begin{subequations}\label{eq:stewart_lp}
			\begin{align}
				\quad \theta(x)  \in\arg\min_{\tilde{\theta} \in \R^{\Nsys}} \quad & g(x)^\top \, \tilde{\theta} \\
				\textrm{s.t.} \quad & e^\top \tilde{\theta} = 1 \label{eq:stewart_lp_eq}
				\\
				&\tilde{\theta}\geq 0 \label{eq:stewart_lp_ineq}.
			\end{align}
		\end{subequations}
		Using its KKT condition, one can obtain a DCS equivalent to \eqref{eq:dcs_step_2}, which reads as:
		\begin{subequations}\label{eq:stewart_dcs}
			\begin{align}
				&\dot{x} = F(x,u) \theta,\\
				&0=g(x)  - \lambda + \mu e,\\
				&0=e^\top  {\theta} -1,\\
				&0 \leq \lambda \perp \theta \geq 0,
			\end{align}
		\end{subequations}
		where $\mu \in \R$  and $\lambda \in \R^{\Nsys}$  are the Lagrange multipliers for the constraints \eqref{eq:stewart_lp_eq} and \eqref{eq:stewart_lp_ineq}, respectively.
		
		\color{black}
		The DCS \eqref{eq:stewart_dcs} is also an instance of the generic DCS \eqref{eq:dcs_general}.
		Set $z = (\theta,\mu)$ and $y = \lambda$. 
		It follows that $f(x,y,z)  = F(x,u) \theta$, $\gdcseq(x,y,z) = (g(x)  - \lambda + \mu e, e^\top  {\theta} -1)$ and $\gdcscomp(x,z,y)  = \theta$.
		\color{black}
		\color{black}
		\subsection{Summary and relations between different formalisms}\label{sec:nonsmooth_systems_summary}
		The diagram in Figure \ref{fig:di_dcs_relations} summarizes the relationships between the nonsmooth systems studied in this paper.
		The first column consists of the different ODEs with DRHS that we have treated.
		After treating a generic ODE with DRHS in \eqref{eq:discontinuous_ode}, we specialize in two structured cases: PSS in \eqref{eq:pss} and ODEs \eqref{eq:ap_ode}, where their dynamics contain some discontinuous expressions of $x$, e.g. as Heaviside step functions.
		
		These ODEs may not have a classical or Carath\'eodory solution in some cases, so we embed them in differential inclusions and obtain more general notions, such as Filippov solutions.
		These concepts are summarized in the second column.
		A generic ODE with DRHS is generalized to Filippov DIs~\eqref{eq:filippov_di}. 
		If the ODE is more structured as a piecewise smooth system, then its Filippov extensions become more explicit~\eqref{eq:filippov_di_with_multiplers}.
		In both cases, the r.h.s. is now a convex and compact set.
		If the discontinuous function in the structured ODE \eqref{eq:ap_ode} is replaced by set-valued extensions, we obtain an Aizerman–Pyatnitskii DI \eqref{eq:di_ap}. 
		The set on the r.h.s. may even be nonconvex.
		In Proposition \ref{prop:filippov_ap} we show that when we use PSS and Heaviside step functions, the Aizerman-Pyatnitskii DI and the Filippov extension for PSS are equivalent.

		The third column consists of dynamic complementarity systems obtained from the DIs, which are useful representations for numerical computations.
		The Heaviside step DCS \eqref{eq:step_dcs} and the Stewart DCS \eqref{eq:stewart_dcs} are both instances of a generic DCS \eqref{eq:dcs_general}. 
		They are derived from the corresponding differential inclusions using the KKT conditions of parametric LPs, as shown in Sections \ref{sec:step_dcs} and \ref{sec:stewarts_representation}, respectively.
		The numerical methods developed in \cite{Nurkanovic2024a} and in this paper exploit the continuity properties of the Lagrange multipliers in the KKT conditions of these LPs.
		We conclude this section by illustrating the different formalisms with an example. 
		\color{black}
		\begin{figure}[t]
			\centering
			\[\begin{tikzcd}
				{\text{ODE with DRHS}}\ \eqref{eq:discontinuous_ode}& {\text{Filippov DI}} \ \eqref{eq:filippov_di} & {\text{DCS} }\ \eqref{eq:dcs_general} \\
				{\text{PSS}}\ \eqref{eq:pss} & {\text{Filippov PSS}}  \ \eqref{eq:filippov_di_with_multiplers}& {\text{Heaviside DCS}} \eqref{eq:step_dcs}\\
				{\text{Structured ODE with DRHS}}\ \eqref{eq:ap_ode} & {\text{Aizerman–Pyatnitskii DI}}\ \eqref{eq:di_ap} & {\text{Stewart DCS}} \eqref{eq:stewart_dcs}
				\arrow[Rightarrow, from=1-1, to=2-1]
				\arrow[shift right=5, curve={height=15pt}, Rightarrow, from=1-1, to=3-1]
				\arrow[Rightarrow, from=1-1, to=1-2]
				\arrow[shorten <=30pt, Rightarrow, from=2-1, to=2-2]
				\arrow[Rightarrow, from=3-1, to=3-2]
				\arrow[Rightarrow, from=1-2, to=2-2]
				\arrow[Rightarrow, 2tail reversed, from=2-2, to=3-3]
				\arrow[Rightarrow, from=2-2, to=2-3]
				\arrow[Rightarrow, from=1-3, to=2-3]
				\arrow[curve={height=-54pt}, Rightarrow, from=1-3, to=3-3]
				\arrow[Rightarrow, from=3-2, to=2-2]
			\end{tikzcd}\]
			\caption{Summary of relations between nonsmooth systems treated in this paper.}
			\label{fig:di_dcs_relations}
		\end{figure}
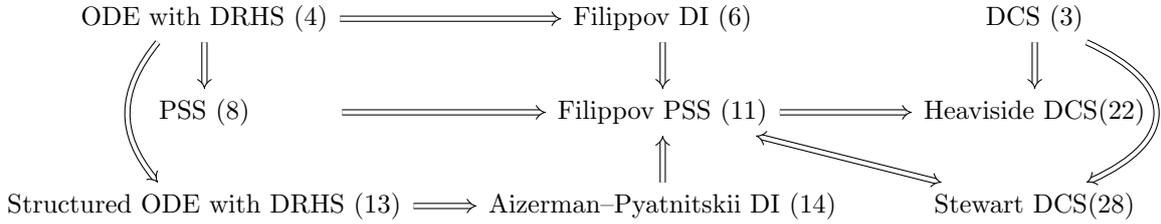
		\begin{example}\label{ex:illustrating_different_formalisms}
			Let us illustrate the different classes of nonsmooth systems on the discontinuous ODE (which is in the form of \eqref{eq:discontinuous_ode}):
			\begin{align*}
				\dot{x} = \begin{cases}
					1, & x >0,\\
					3, & x <0.
				\end{cases}
			\end{align*}
			This is a PSS (as in \eqref{eq:pss}) with the switching function  $\switchfun(x) = x$, with the regions $R_1 = \{ x \mid x > 0\}$ and 
			$R_2 = \{ x \mid x < 0\}$. 
			The Filippov extensions of this PSS (Eq. \eqref{eq:filippov_di_with_multiplers}) reads as:
			$\dot{x} \in \{ \theta_1 + 3\theta_2 \mid \theta \geq 0, \theta_1+\theta_2 = 1\}$, with $\theta = (\theta_1,\theta_2)$.  
			In the form of an Aizerman–Pyatnitskii DI \eqref{eq:basic_di}, i.e. \eqref{eq:di_ap}, the system reads as $\dot{x} \in 3-2\gamma(x)$. 	
			We can write also as a DCS of the form of \eqref{eq:step_dcs}:
			\begin{align*}
				&\dot{x} = \begin{bmatrix}
					3 & 1 
				\end{bmatrix}\, \theta,\\
				&\theta_1 = \alpha,\ \theta_2 =  1-\alpha,\ x = \lambdap-\lambdan,\\ 
				&0 \leq \lambdap \perp \alpha \geq 0,\ 0 \leq \lambdan \perp 1-\alpha \geq 0.
			\end{align*}
			Similarly, by defining the indicator function $g(x) = (-x,x)$, we can state Stewart DCS \eqref{eq:stewart_dcs}:
			\begin{align*}
				&\dot{x} = \begin{bmatrix}
					3 & 1 
				\end{bmatrix}\, \theta,\\
				&-x = \lambda_1 -\mu,\ x = \lambda_2 -\mu,\ \theta_1 +\theta_2 = 1,\\ 
				&0 \leq \theta_1 \perp \lambda_1 \geq 0,\ 0 \leq \lambda_2 \perp \theta_2 \geq 0.
			\end{align*}
		   \end{example}

		\color{black}
\section{Properties of the step representation DCS}\label{sec:dcs_properties}
In this section, we study some properties of the DCS obtained via Heaviside step functions \eqref{eq:step_dcs} for a fixed active set and at active-set changes.
These properties will be useful for algorithmic development in the subsequent sections.
\textcolor{black}{For a fixed active set, i.e. $\mathcal{I}(x(t)) = \{ i\mid \theta_i(t) >0 \} $ being constant on some time interval, there are no switches and the dynamics have locally no discontinuities. 
On the other hand, at an active-set change there is a switch and a discontinuity in the dynamics.}

\subsection{Active-set changes and continuity of $\lambdap$ and $\lambdan$}
\label{sec:step_dcs_active_set_changes}
Active-set changes are paired with discontinuities in some of the algebraic variables.
We have seen \textcolor{black}{in Lemma \ref{lem:cont_lambda}} that $\lambdap$ and $\lambdan$ are continuous functions of time.

\color{black}
At an active-set change, we at least one of the switching functions $\switchfun_j(x)$ either becomes zero, or if it was zero it becomes nonzero.
It follows from $\psi_j(x(t)) = \lambdap_j(t)- \lambdan_j(t)$ (in \eqref{eq:step_continuity_lambda}), that also both $\lambdap_j(t)$ and $\lambdan_j(t)$ must be zero at an active-set change.
\color{black}
We use now the DCS formulation via step functions in \eqref{eq:step_dcs} to illustrate
the different switching cases that arise in Filippov systems.
\begin{figure}[t]
	\centering
	{\includegraphics[scale=0.63]{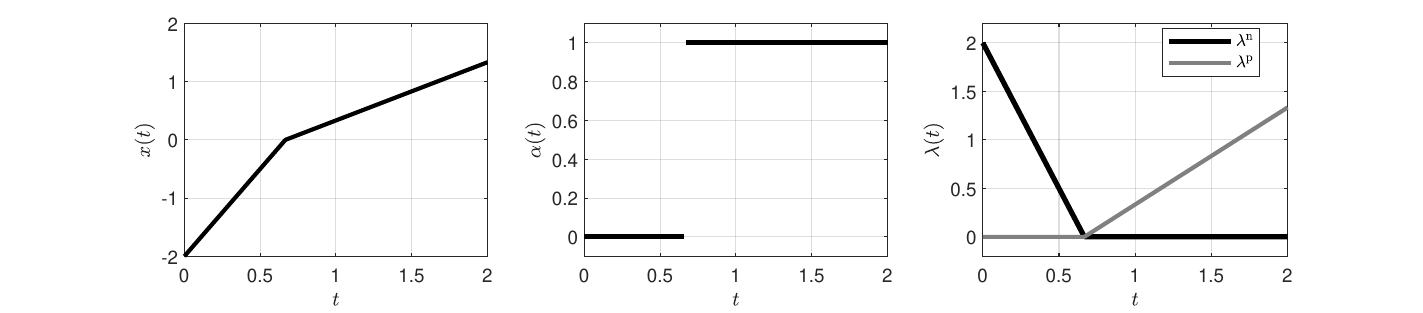}}
	{\includegraphics[scale=0.63]{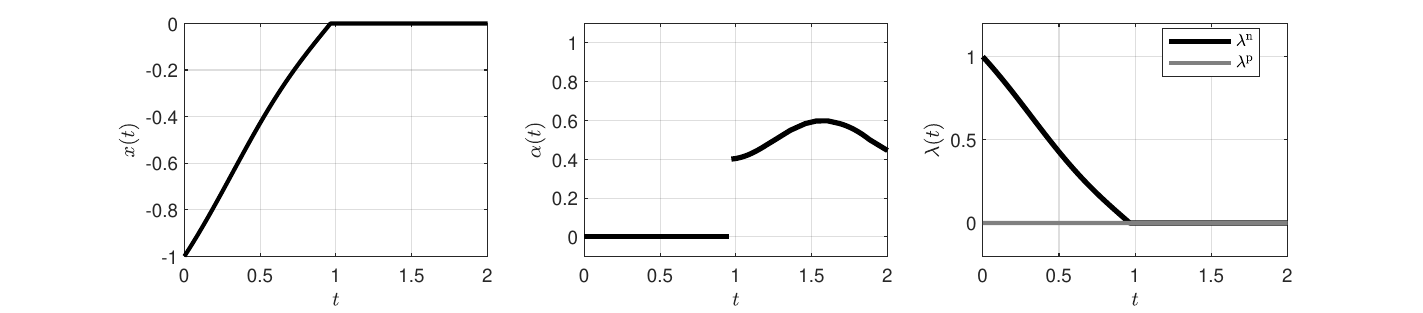}}
	{\includegraphics[scale=0.63]{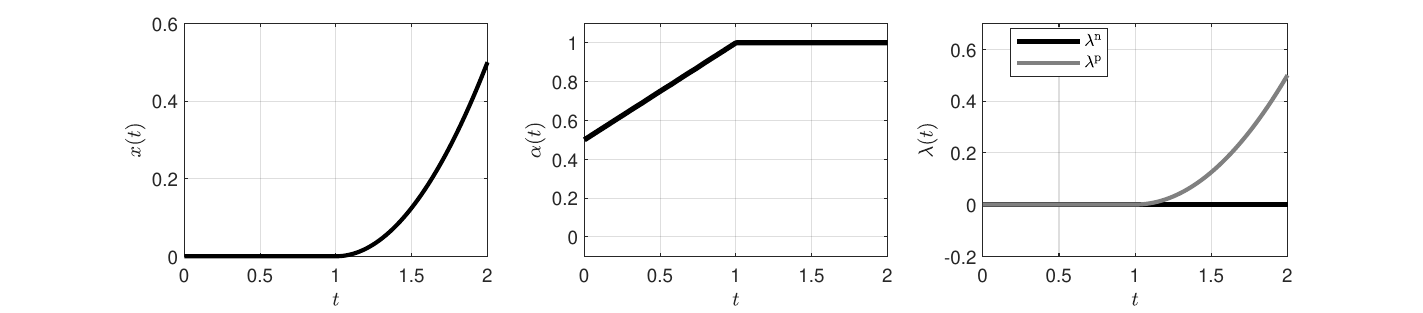}}
	{\includegraphics[scale=0.63]{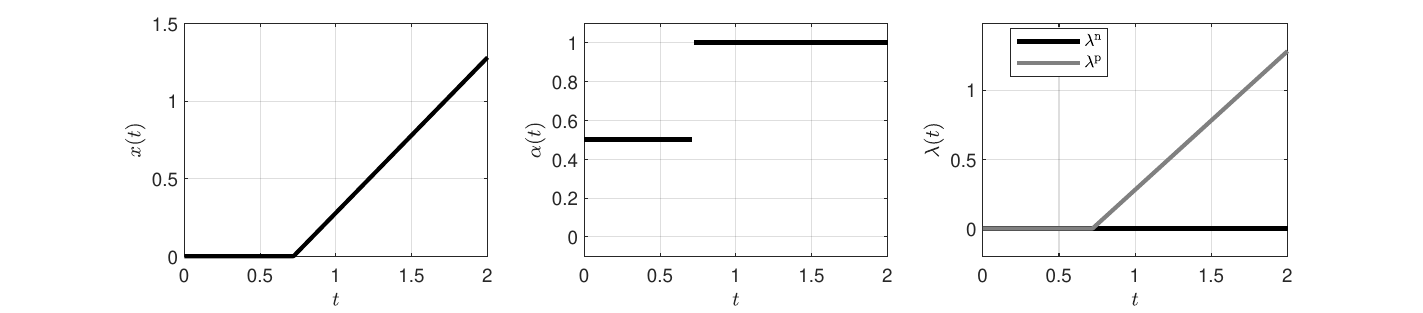}}
	\caption{
		Illustration of example solution trajectories for different switching cases.
		The rows from top to bottom show $x(t)$, $\alpha(t)$, $\lambdap(t)$ and  $\lambdan(t)$ for the cases (a)-(d) in Example~\ref{ex:switching_cases}, respectively.
	}
	\label{fig:switcing_cases_step}
	\vspace{-0.2cm}
\end{figure}
\begin{example}\label{ex:switching_cases}
	There are four possible switching cases which we illustrate with the following examples:
	\begin{enumerate}[(a)]
		\item crossing a surface of discontinuity, $\dot{x}(t) \in 2-\mathrm{sign}(x(t))$ (same system as in Example \ref{ex:illustrating_different_formalisms}),
		\item entering a sliding mode, $\dot{x}(t) \in -\mathrm{sign}(x(t))+0.2\sin(5t)$,
		\item leaving a sliding mode $\dot{x}(t) \in -\mathrm{sign}(x(t))+t$,
		\item spontaneous switch, $\dot{x}(t) \in \mathrm{sign}(x(t))$.
	\end{enumerate}
	In case (a), for $x(0) < 0$ the trajectory reaches $x = 0$ and crosses it.
	In example (b), for any finite $x(0)$, the trajectory reaches $x=0$ and stays there.
	On the other hand, in example (c), for $x(0) =0$, the DI has a unique solution and leaves $x=0$ at $t=1$.
	In the last example, the DI has infinitely many solutions for $x(0) =0$, and $x(t)$ can spontaneously leave $x = 0$ at any $t\geq0$.
	The trajectories are illustrated in Figure~\ref{fig:switcing_cases_step}.
	\color{black}
	Whenever $x(t)$ has a kink, which corresponds to a switch and discontinuity in the dynamics, both the Lagrange multipliers $\lambdap(t)$ and $\lambdan(t)$ are zero at that time.
	\color{black}
\end{example}

\subsection{Fixed active set in the step formulation}\label{sec:step_dcs_active_set_fixed}
We study the properties of the DCS \eqref{eq:step_dcs} for a fixed active set $\I(x(t))$.
Without loss of generality, the corresponding time interval is $\mathcal{T} = (0,T)$.
For a fixed active set, the DCS \eqref{eq:step_dcs} reduces either to an ODE or to a Differential Algebraic Equation (DAE).

We start with the simpler ODE case.
Let $\switchfun_j(x) \neq 0$ for all $j \in \C \coloneqq \{1,\ldots,\nswitchfun\}$, then $x(t)$ is in the interior of some region $R_i$.
It can be seen from the LP \eqref{eq:step_parametric_lp} that $\alpha_j \in \{0,1\}$ for all $j \in \C $.
This implies that $\theta_i=1$ and $\theta_k = 0, k \neq i$.
\color{black}
Therefore, $\I(x(t))  = \{i\}$ and the Filippov DI reduces to $\dot{x} \in F_{\mathrm{F}}(x) = \{f_i(x)\}$, i.e., we have locally an ODE.

Next, we regard the case when $\I(x(t))$ is not a singleton, i.e., the trajectory evolves at the boundary of two or more regions.
\color{black}
Consequently, we have at least one $\switchfun_j(x) = 0$.
Let us associate with $\I(x(t))$ the index set $\K(x(t))  = \{ j \in \C \mid  \switchfun_j(x) = 0\}$, i.e., the set of indices of all switching functions that are zero for a given active set~$\I(x(t))$.

In the sequel, we make use of the following notation. For a given vector $a\in\R^n$ and set $\I \subseteq \{1,\ldots,n\}$, we define the projection matrix $P_\I \in \R^{|I| \times n}$, which has zeros or ones as entries. 
It selects all component $a_i, i\in \I$ from the vector $a$, i.e., $a_\I = P_\I a \in \R^{|\I|}$ and $a_\I = [a_i \mid i\in \I]$.

Following the discussion from the previous section, for all nonzero $\switchfun_j(x)$, i.e., $j \in \C \setminus \K$, we can compute $\alpha_j \in \{0,1\}$ via the LP \eqref{eq:step_parametric_lp} and $\lambda_{\mathrm{p},j}, \lambda_{\mathrm{n},j}$ via  \eqref{eq:step_continuity_lambda}.
Next, we have that $\lambda_{\mathrm{p},j} = \lambda_{\mathrm{n},j} = 0$ for all $j \in \K$.
It is left to determine $\alpha_j$ for all $j \in \K$ and thus implicitly all $\theta_i$, for all $i \in \I$.
Recall that $\theta_i = 0$ for all $i \notin \I$.
By fixing the already known variables in \eqref{eq:step_dcs} we obtain the DAE:
\begin{subequations}\label{eq:step_dae}
	\begin{align}
		&\dot{x} = F_{\mathcal{I}}(x,u)\; \theta_{\mathcal{I}}, \\
		&\theta_i - \prod_{j =1}^{\nswitchfun}  \Big(\frac{1-{S}_{i,j}}{2}+{S}_{i,j}\alpha_j\Big) = 0,\; i\in\I, \label{eq:step_dae_theta}\\
		&\switchfun_j(x) = 0,\;  j\in \K \label{eq:step_dae_c},
	\end{align}
\end{subequations}
where we define \color{black} $F_{\I}(x,u) \coloneqq F(x,u) P_{\I}^\top \in \R^{n_x \times |\I|}$ \color{black}, i.e., we select only the columns of $F(x,u)$ with the index $i \in \I$.
Note that  $\alpha_j$ for all $j \in \mathcal{C}\setminus\mathcal{K}$ are known and thus no degrees of freedom.
We keep them for ease of notation.
Thus we have a DAE with $|\I| + |\K|$ unknowns, namely $\theta_{\I} \in \R^{|\I|}$ and $\alpha_{\K} \in \R^{|\K|}$, and $|\I| + |\K|$ algebraic equations in \eqref{eq:step_dae_theta} and \eqref{eq:step_dae_c}.

Next, we investigate conditions under which the DAE \eqref{eq:step_dae} is well-posed.
For this purpose, we define the matrix
\color{black}
\begin{align*}
	W_{\K,\I}(x,u) &= \nabla \switchfun_{\K}(x)^\top F_{\I}(x,u) \in \R^{|\K| \times |\I|},
\end{align*}
where $\nabla \switchfun_{\K}(x) = \begin{bmatrix}
	\nabla \switchfun_j(x) \mid j \in \K
\end{bmatrix} \in \R^{n_x \times |\K|}$  is a matrix, whose columns are the gradients of the switching functions that are zero for the given active set $\I$. 
\color{black}
Moreover, we define a compact notation for the partial Jacobian $B_{\K,\I}(\alpha)  \in \R^{|\I| \times |\mathcal{K}|}$  of \eqref{eq:step_dae_theta} w.r.t. to $\alpha_{\K}$, with the elements:
\begin{align*}
	&B_{i,j}(\alpha) \coloneqq 
	\frac{\partial}{\partial \alpha_{j}} \prod_{l \in \mathcal{C}}  \frac{1-{S}_{i,l}}{2}+{S}_{i,l}\alpha_l,\ i \in \I, j \in \mathcal{K}.
\end{align*}
\color{black}
\begin{assumption}\label{ass:solution_existence_step}
	Given a fixed active set $\I(x(t)) = \I$ for $t \in  \mathcal{T}$, it holds that the matrix functions $W_{\mathcal{K},\I}(x,u)$ and $B_{\mathcal{K},\I}(\alpha)$ are Lipschitz continuous, and that $W_{\mathcal{K},\I}(x,u) B_{\mathcal{K},\I}(\alpha)$ has rank~$|\K|$, i.e. it is full rank, for all $t\in \mathcal{T}$.
\end{assumption}
\color{black}
\begin{proposition}\label{eq:step_dae_fixed_active_set}
	Suppose that Assumption\ref{ass:solution_existence_step} holds.
	Given an initial value $x(0)$, the DAE \eqref{eq:step_dae} has a unique solution for all $t\in \mathcal{T}$.
\end{proposition}
\textit{Proof.}
First, we differentiate \eqref{eq:step_dae_c} 
\color{black}
with respect to $t$, such that algebraic variables appear explicitly in the algebraic equations (this correspond to so-called index reduction in the theory of DAEs, cf. \cite{Hairer1991}):
\color{black}
\begin{subequations}\label{eq:step_dae_index1}
	\begin{align}
		&\dot{x} = F_{\mathcal{I}}(x,u)\; \theta_{\mathcal{I}}, \label{eq:step_dae_index1_ode}\\
		&\theta_i - \prod_{j \in \C}  \frac{1-{S}_{i,j}}{2}+{S}_{i,j}\alpha_j = 0,\; i\in\I, \label{eq:step_dae_index1_theta}\\
		&W_{\K,\I}(x,u) \theta_\I = 0 \label{eq:step_dae_index1_c}.
	\end{align}
\end{subequations}
\color{black}
Next, we prove that the partial Jacobian of \eqref{eq:step_dae_index1_theta}-\eqref{eq:step_dae_index1_c} w.r.t. to the algebraic variables $(\theta_{\I},\alpha_{\K})$ has rank $|\I| + |\K|$, i.e., it is an invertible matrix.
\color{black}
We omit the dependencies on $\alpha$ and $x$ for brevity.
The Jacobian of \eqref{eq:step_dae_index1_theta}-\eqref{eq:step_dae_index1_c} w.r.t. to $(\theta_{\I},\alpha_{\K})$
has the form
\begin{align*}
	A =	\begin{bmatrix}
		I_{|\I|} & -B_{\K,\I}\\
		W_{\K,\I} & \mathbf{0}
	\end{bmatrix}.
\end{align*}
\textcolor{black}{To prove that this matrix has full rank, we show that the only solution $(v,w)\in \R^{|\I|} \times \R^{|\K|}$, to the following linear system is the zero vector:}
\begin{align}\label{eq:linear_sys_step_dae}
	\begin{bmatrix}
		I_{|\I|} & -B_{\I,\K}\\
		W_{\K,\I} & \mathbf{0}
	\end{bmatrix}
	\begin{bmatrix}
		v \\w
	\end{bmatrix} = 0.
\end{align}
From the first line we have $v = B_{\I,\K}w$ and substituting this into the second line we have
$W_{\I,\K}B_{\I,\K}w = 0$. 
Since the matrix $W_{\I,\K}B_{\I,\K}\in \R^{|\K|\times |\K|}$ has rank $|\K|$, the only solution to \eqref{eq:linear_sys_step_dae} is $w =0$, and $v = Bw = 0$. 
Hence, $A$ has full rank. 

\color{black}
Now we can apply the implicit function theorem \cite[Theorem 1B.1]{Dontchev2014} to \eqref{eq:step_dae_index1_theta}-\eqref{eq:step_dae_index1_c}, which guarantees the existence of continuously differentiable functions $\theta_{\I}(x)$ and $\alpha_{\K}(x)$.  
Since the function $\theta_{\I}(x)$ is continuously differentiable, it is also Lipschitz continuous for a fixed $\I(x(t))$ on $t\in\mathcal{T}$.
By substituting $\theta_{\I}(x)$ into \eqref{eq:step_dae_index1_ode}, we have a product of two Lipschitz continuous functions (all columns of $F_{\I}(x,u)$, are Lipschitz by assumption), and the DAE \eqref{eq:step_dae_index1} reduces to an ODE with a Lipschitz continuous r.h.s.
\color{black}
This enables us to apply the Picard-Lindel\"of Theorem to obtain the assertion of the proposition.
\qed

	\color{black}
	We make a few comments on Assumption~\ref{ass:solution_existence_step}.
	The rank condition can be checked explicitly, since one can simply compute the matrix $W_{\mathcal{K},\I}(x,u) B_{\mathcal{K},\I}(\alpha)$. 
	We have already assumed Lipschitz continuity of all columns of $f_i(x,u)$ and of the gradients $\nabla \switchfun_j(x)$.  Here we additionally assume it for the matrix $W_{\K,\I}(x,u)$, whose entries are computed as inner products on these vectors. 
	The entries of the matrix $B_{\mathcal{K},\I}(\alpha)$ are multi-affine terms, which are also Lipschitz, at least on the bounded domains we consider here.
	In Stewart's reformulation, we consider a square matrix with entries $\nabla g_i(x)^\top f_j(x,u)$ (which is structurally similar to $W_{\K,\I}(x,u)$). 
	For well-posedness with a fixed active set, the invertibility of this matrix is assumed \cite{Nurkanovic2024a,Stewart1990a}.	
	\color{black}
	
	In \cite{Dieci2011}, the authors make some assumptions on the signs of the entries of $W_{\K,\I}(x)$ and prove the existence, but not uniqueness, of solutions with a fixed-point argument.
	For the case of $|\mathcal{K}| \leq 2$, i.e., sliding modes with co-dimension one or two, and with additional assumptions on the signs of the entries of $W_{\K,\I}(x)$ they even prove the uniqueness of solutions.
	
	Observe that for a given $x(t)$, there might be several $\I(x(t))$ that yield meaningful DAEs of the form of Eq.~\eqref{eq:step_dae}.
	This may happen when the Filippov DI does not have unique solutions, such as in Example \ref{ex:switching_cases} case~(d). 
	\color{black}
	The trajectory may stay in sliding mode or leave it any point of time. 
	The trajectory pieces in either scenario are well-posed, even thought the overall trajectory is not unique.
	\color{black}

\section{Finite Elements with Switch Detection (FESD) for the step representation}\label{sec:fesd_step}
\color{black}
\begin{figure}[th]
	\centering
	{\includegraphics[width= 0.32\textwidth]{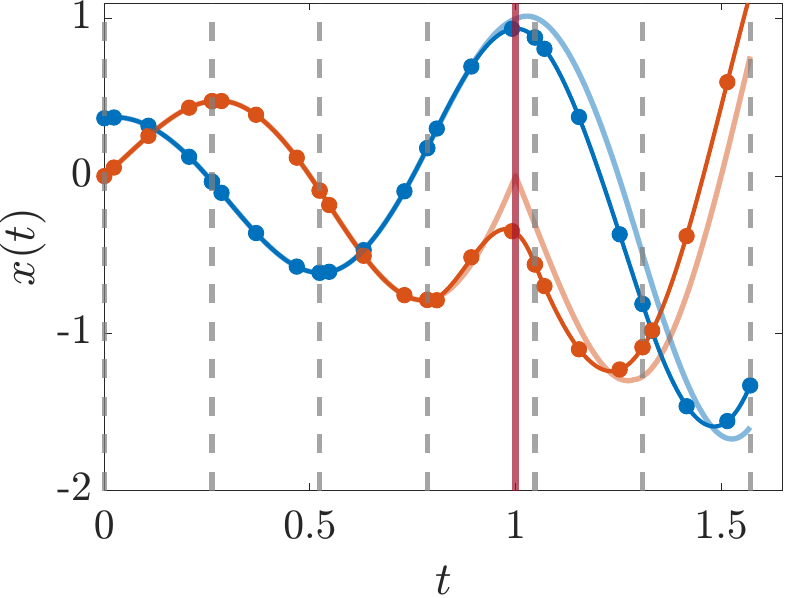}}
	{\includegraphics[width= 0.32\textwidth]{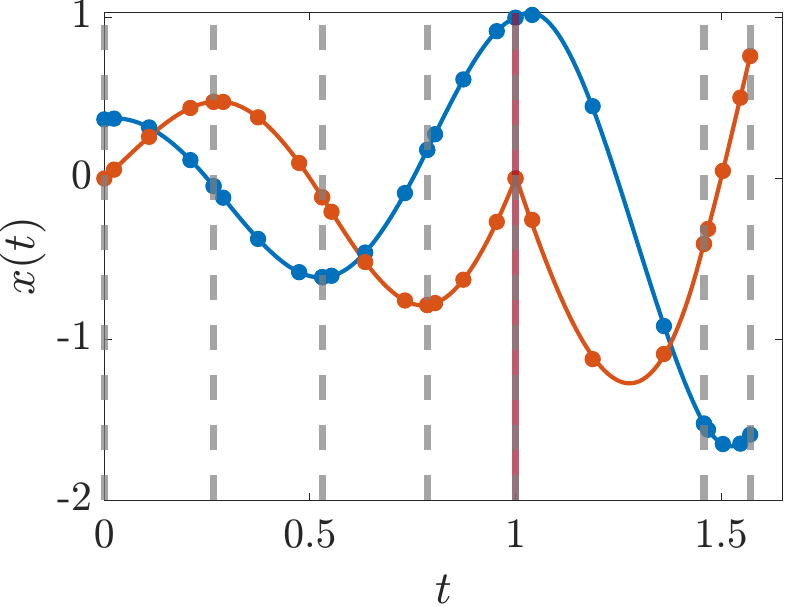}}
	{\includegraphics[width= 0.32\textwidth]{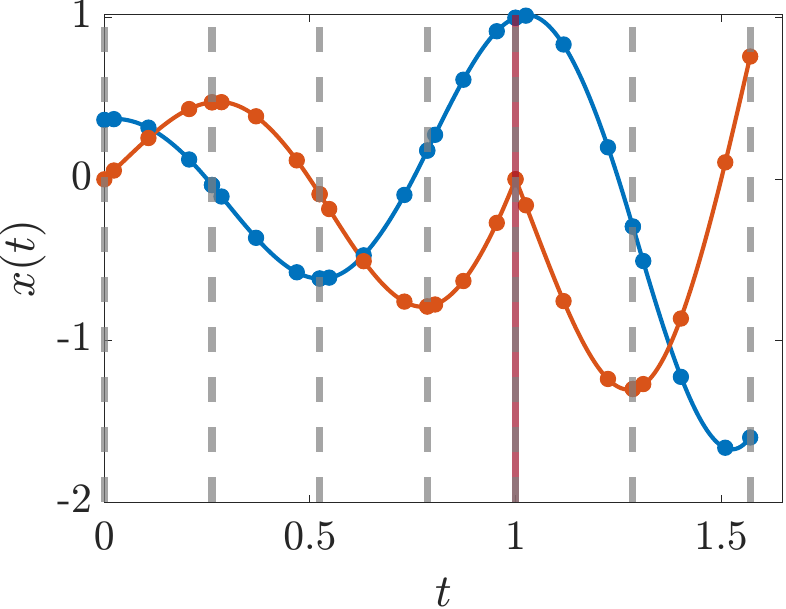}}
	\caption{
		\color{black}
		Example trajectories of an ODE with DRHS, with a switch at $t = 1$. 
		The solid lines correspond to numerical approximations, the transparent ones to analytic solutions. 
		The circle markers correspond to the stage values of an underlying Runge-Kutta (RK) method, the vertical dashed lines to the boundaries of the integration intervals. 
		The left plots show an approximation obtained with a time-stepping RK method, the middle one with FESD but without step equilibration, and the right one with FESD including step equilibration.
	}
	\label{fig:fesd_motivation}
	\vspace{-0.35cm}
\end{figure}
\subsection{Main ideas}\label{sec:fesd_intro}
We outline the main ideas in the derivation of the FESD discretization, and in the following sections we explain each step in more detail. 
The starting point is a standard time-stepping Runge-Kutta (RK) discretization of the DCS~\eqref{eq:step_dcs}.
In this case, a fixed integration step size is assumed. 
If the switch occurs within an integration interval, the high accuracy properties of these methods are lost. 
This is shown in the left plot of Figure \ref{fig:fesd_motivation}. 
The switch occurred at an RK stage point inside the integration interval, and thereafter there is a large discrepancy between the numerical approximation and the exact solution. 
Nevertheless, the standard RK methods are a starting point for the development of the FESD method and we recall them in Section~\ref{sec:step_dcs_rk}.

If the switches always coincided with the integration interval boundaries, the RK method would not lose its accuracy.
This can be achieved by detecting the switches and adjusting the step size.
To achieve this, in the FESD method~\cite{Nurkanovic2024a}, the integration step sizes are left as degrees of freedom. 
This idea was first proposed by Baumrucker and Biegler in \cite{Baumrucker2009} and extended and theoretically analyzed in \cite{Nurkanovic2024a}.
Moreover, additional complementarity conditions, called cross complementarities, prohibit switches within an integration interval (finite element), which leads to their isolation at the boundaries and recovery of high accuracy. 
An illustration is given in the middle plot of Figure \ref{fig:fesd_motivation}. 
In contrast to the fixed step size discretization, the numerical and exact solutions are indistinguishable.
We study these extensions in Section \ref{sec:fesd_cross_comp_step_reformulation}.

However, when there are no switches, the step sizes are not uniquely determined. 
As can be seen in the middle plot of Figure \ref{fig:fesd_motivation}, the integration intervals all have different and somewhat random lengths.
In Section \ref{sec:fesd_step_equilibration_step_reformulation} we introduce the step equilibration conditions that make neighboring finite elements have the same length if no switches occur, In Section \ref{sec:fesd_step_equilibration_step_reformulation} we introduce the step equilibration conditions that make neighboring finite elements have the same length if no switches occur, and thus determine them uniquely.
The right plot in Figure \ref{fig:fesd_motivation} illustrates this effect, where now before and after the switch have equidistant grids.
This is also a key ingredient for having locally unique solutions, as we will show later.

In Section \ref{sec:fesd_all_fesd_eq_step}, we summarize the developments of the previous the sections and summarize the full FESD discretization in a compact way.
We emphasize that in the FESD method, switch detection is fully implicit and based only on algebraic conditions.
This makes it suitable for the discretization of optimal control problems.
We conclude with Section \ref{sec:fesd_ocp}, where we show how to apply the FESD method to discretize an optimal control problem subject to a nonsmooth dynamical system with Heaviside step functions.
\color{black}
\subsection{Standard implicit Runge-Kutta discretization}\label{sec:step_dcs_rk}
As a starting point in our analysis, we regard a standard Runge-Kutta (RK) discretization of the DCS~\eqref{eq:step_dcs}.
For ease of notation, we work with the nonsmooth DAE formulation of the DCS \eqref{eq:dcs_step_2}, which we restate for convenience
\begin{subequations}\label{eq:dcs_step_2_restate}
	\begin{align}
		&\dot{x}  = F(x,u)\theta,\\
		&0= G(x,\theta,\alpha,\lambdap,\lambdan) \label{eq:dcs_step_2_restate_alg}
	\end{align}
\end{subequations}
In the sequel, one should keep in mind that \eqref{eq:dcs_step_2_restate_alg} collects all algebraic equations including the complementarity conditions \eqref{eq:step_dcs_cc1} and \eqref{eq:step_dcs_cc2}.

\color{black}
In a discretized Optimal Control Problem (OCP) we usually have several control intervals. 
In our derivations here, it is sufficient to consider a single control interval $[0,T]$ with a constant control input $q \in \R^{n_u}$, i.e., we set $u(t) = q$ for $t\in [0,T]$. 
The extension to more elaborate control parameterizations is straightforward~\cite{Rawlings2017}.
In Section \ref{sec:fesd_ocp} we show how to discretize and OCPs with multiple control intervals.
\color{black}

Let $x(0) = s_0$ be the given initial value.
The control interval $[0,T]$ is divided into $\NFE$ {finite elements} (i.e., integration intervals) $[t_{n},t_{n+1}]$ via the grid points $0= t_0 < t_1 < \ldots <t_{\NFE} = T$.
On each of the finite elements we regard an $\Nstg$-stage Runge-Kutta method which is characterized by the Butcher tableau entries $a_{i,j} ,b_i$ and $c_i$ with $i,j\in\{1,\ldots,\Nstg\}$~\cite{Hairer1991}.
The step-sizes read as $h_{n} = t_{n+1} - t_{n},\; n = 0, \ldots,\NFE-1$.
The approximation of the differential state at the grid points $t_n$ is denoted by $x_n \approx x(t_n)$.

We regard the so-called differential representation of the Runge-Kutta method~\cite{Hairer1991}.
Thus, the derivatives of the states at the RK stage points $t_{n,i} \coloneqq t_n + c_i h_n,\; i = 1,\ldots, \Nstg$, are the degrees of freedom.
For a single finite element, we summarize them in the vector $V_n \coloneqq (v_{n,1}, \ldots, v_{n,\Nstg}) \in \R^{\Nstg n_x}$.
The stage values for the algebraic variables are collected in the vectors:
$\Theta_n \coloneqq (\theta_{n,1}, \ldots, \theta_{n,\Nstg} )\in \R^{\Nstg \cdot \Nsys}$,
${A}_n \coloneqq (\alpha_{n,1}, \ldots, \alpha_{n,\Nstg} )\in \R^{\Nstg \cdot n_c}$,
$\Lambdap_n \coloneqq (\lambdap_{n,1}, \ldots, \lambdap_{n,\Nstg} )\in \R^{\Nstg \cdot n_c}$ and
$\Lambdan_n \coloneqq (\lambdan_{n,1}, \ldots, \lambdan_{n,\Nstg} )\in \R^{\Nstg \cdot n_c}$.

We collect all \textit{internal} variables in the vector
$Z_n =(x_n,\Theta_n,A_n,\Lambdap_n,\Lambdan_n,V_n)$.
The vector $x_n^{\mathrm{next}}$ denotes the state value at $t_{n+1}$, which is obtained after a single integration step.
Now, we can state the RK equations for the DCS \eqref{eq:dcs_step_2_restate} for a single finite element as
\begin{align}\label{eq:dcs_irk_single_step_reformulation}
	&0 = G_{\irk}(x_n^{\mathrm{next}}\!,Z_n,h_n,q)\!\coloneqq\! 
	\begin{bmatrix}
		\! v_{n,1}\! -\!  F(x_n +h_n \sum_{j=1}^{\Nstg} a_{1,j} v_{n,j},q)\theta_{n,1}\\
		\vdots\\
		v_{n,\Nstg} \! -\! F(x_n +h_n \sum_{j=1}^{\Nstg} a_{\Nstg,j} v_{n,j},q)\theta_{n,\Nstg}\\
		G(x_n + h_n\sum_{j=1}^{\Nstg} a_{1,j} v_{n,j},\theta_{n,1},\alpha_{n,1},\lambdap_{n,1},\lambdan_{n,1})\\
		\vdots\\
		G(x_n \!+\! h_n\sum_{j=1}^{\Nstg}\!a_{\Nstg,j} v_{n,j},\theta_{n,\Nstg},\alpha_{n,\Nstg},\lambdap_{n,\Nstg},\lambdan_{n,\Nstg})\\
		x_n^{\mathrm{next}} - x_n - h_n \sum_{i=1}^{\Nstg} b_i v_{n,i}
	\end{bmatrix}.
\end{align}
Next, we summarize the equations for all $\NFE$ finite elements over the entire interval $[0, T]$ in a discrete-time system format.
To make it more manageable, we use some additional shorthand notation and group all variables of all finite elements for a single control interval into the following vectors:
$\mathbf{x}= (x_0,\ldots,x_{\NFE}) \in \R^{(\NFE+1) n_x}$,
$\mathbf{V} = (V_0,\ldots,V_{\NFE-1}) \in \R^{\NFE \Nstg n_x}$ and
$\mathbf{h}\coloneqq (h_0,\ldots,h_{\NFE-1})\in \R^{\NFE}$.
Recall that the simple continuity condition $x_{n+1} = x_{n}^{\mathrm{next}}$ holds.
We collect the stage values of the Filippov multipliers in the vector
$\mathbf{\Theta} = ({\Theta}_0,\ldots,\Theta_{\NFE-1})\in \R^{n_{{\theta}}}$ and
$n_{{\theta}}= \NFE\Nstg\Nsys$.
Similarly, we group the stage values of the algebraic variables in the vectors
$\mathbf{A},\mathbf{\Lambda}^{\mathrm{p}}, \mathbf{\Lambda}^{\mathrm{n}}\in \R^{n_{\alpha}}$, where $n_{\alpha} = \NFE\Nstg n_c$.
Finally, we collect all internal variables in the vector
$\mathbf{Z} = (\textbf{x},\mathbf{V},\mathbf{\Theta},\mathbf{A},\mathbf{\Lambda}^{\mathrm{p}},\mathbf{\Lambda}^{\mathrm{n}})\in \R^{n_{\mathbf{Z}}}$, where $n_{\mathbf{Z}} = (\NFE+1)n_x + \NFE\Nstg n_x + n_{\theta}+3n_{\alpha}$.

All computations over a single control interval of the {standard discretization} (denoted by the subscript $\textrm{std}$ in the corresponding functions) are summarized in the following equations:
\begin{subequations}\label{eq:dcs_irk_step_reformulation}
	\begin{align}
		{s}_1 \! &= \!F_{\mathrm{std}}(\textbf{Z}),\\
		0\! &=\! {G}_{\mathrm{std}}(\mathbf{Z},\mathbf{h},s_0,q),
	\end{align}
\end{subequations}
where $s_1\in \R^{n_x}$ is the approximation of $x(T)$ and
\begin{align*}
	F_{\mathrm{std}}(\textbf{Z})  &= x_{\NFE},\\
	G_{\mathrm{std}}(\mathbf{Z},\mathbf{h},s_0,q)
	\coloneqq	&
	\begin{bmatrix}
		x_0- s_0\\
		G_{\irk}(x_1,Z_0,h_0,q)\\
		\vdots\\
		G_{\irk}(x_{\NFE},Z_{\NFE-1},h_{\NFE-1},q)
	\end{bmatrix}.
\end{align*}
In \eqref{eq:dcs_irk_step_reformulation}, $\mathbf{h}$ is a given parameter and implicitly fixes the discretization grid.
In contrast to standard RK discretizations, we will now proceed by letting $\mathbf{h}$ be degrees of freedom and introduce the cross-complementarity conditions.

\subsection{Cross complementarity}\label{sec:fesd_cross_comp_step_reformulation}
For ease of exposition, suppose that the underlying RK scheme satisfies $c_{\Nstg} =1$. 
This means that the right boundary point of a finite element is a stage point, since $t_{n+1} = t_n+c_{\Nstg} h_n$.
For example, this assumptions is satisfied for Radau and Lobatto schemes \cite{Hairer1991}.
We provide extensions for $c_{\Nstg} \neq 1$ at the end of the section.

The goal is to derive additional constraints that will allow active-set changes only at the boundary of a finite element, compare left and middle plots in Figure \ref{fig:fesd_motivation}.
Moreover, in this case, the step size $h_n$ should adapt such that the switch is detected exactly.
Recall that for the step reformulation at every stage point we have the complementarity conditions:
\begin{subequations}\label{eq:fesd_standard_comp_step_reformulation}
\begin{align}
	0 \leq &\lambdan_{n,m} \perp \alpha_{n,m} \geq  0,\;&\; n= 1,\ldots,\NFE m = 1,\ldots,\Nstg,\\
	0 \leq &\lambdap_{n,m} \perp e-\alpha_{n,m} \geq  0,\; \
	&\textrm{for all}\; n = 1,\ldots,\NFE\}, m = 1,\ldots,\Nstg,
\end{align}
\end{subequations}
\color{black}
where $n$ is the index of the finite elements (integration interval) and $m$ the index of the RK-stage.
\color{black}
We exploit the continuity of the Lagrange multipliers $\lambdap$ and $\lambdan$, cf. Lemma \ref{lem:cont_lambda}.
We regard the boundary values of the approximation of $\lambdap$ and $\lambdan$ on an interval $[t_n,t_{n+1}]$.
They are denoted by $\lambdap_{n,0},\; \lambdan_{n,0}$ (which we define \secsub{below}) at $t_n$ and $\lambdap_{n,\Nstg},\; \lambdan_{n,\Nstg}$ at $t_{n+1}$.

Next, we impose a continuity condition for the discrete-time versions of $\lambdap$ and $\lambdan$ for all ${n \in \{0,\ldots, \NFE-1\}}$:
\begin{align}\label{eq:continuity_of_lambda_step_representation}
	\lambdap_{n+1,0} = \lambdap_{n,\Nstg},\;
	\lambdan_{n+1,0}= \lambdan_{n,\Nstg}.
\end{align}
Note that $\lambdap_{0,0}$ and $\lambdan_{0,0}$ are not defined via~\eqref{eq:continuity_of_lambda_step_representation}, as we do not have a preceding finite element for $n=0$.
Nevertheless, they are crucial for determining the active set in the first finite element.
They are not degrees of freedom but parameters determined by a given $s_0$.
Using \eqref{eq:step_continuity_lambda} we obtain $\lambdap_{0,0} = \max(\switchfun(s_0),0)$ and $\lambdan_{0,0}= -\min(\switchfun(s_0),0)$.

We have seen in Section \ref{sec:step_dcs_active_set_changes} that, due to continuity, $\lambdap_i(t)$ and $\lambdan_i(t)$ must be zero at an active set change, see also Figure \ref{fig:switcing_cases_step}.
Moreover, on an interval $t\in(t_n,t_{n+1})$ with a fixed active set, the components of these multipliers are either zero or positive on the whole interval.
The discrete-time counterparts, i.e., the stage values $\lambdap_{n,m}$ and $\lambdan_{n,m}$ should satisfy these properties as well.
We achieve these goals via the cross complementarity conditions, which read as, for all  $n \in \{0,\ldots,\NFE\!-\!1\}$:
\begin{subequations}\label{eq:cross_cc_true_step_reformulation}
	\begin{align}
		0&  = \mathrm{diag}(\lambdan_{n,m'})\alpha_{n,m},& m = 1,\ldots, \Nstg,\, m' = 0,\ldots, \Nstg, \, m  \neq m',
		\\
		0&  = \mathrm{diag}(\lambdap_{n,m'})(e-\alpha_{n,m}), & m = 1,\ldots, \Nstg,\, m' = 0,\ldots, \Nstg, \, m  \neq m'.
	\end{align}
\end{subequations}
In contrast to ~\eqref{eq:continuity_of_lambda_step_representation}, here we have conditions relating variables corresponding to different RK stages within a finite element. 
Equation \eqref{eq:cross_cc_true_step_reformulation} extends the complementarity conditions for the same RK-stage, i.e., for $m  = m'$, which are part of the standard RK equations, cf.~\eqref{eq:fesd_standard_comp_step_reformulation}.

Some of the claims about the constraints \eqref{eq:cross_cc_true_step_reformulation} are formalized by the next lemma.
Recall that in our notation $\alpha_{n,m,j}$ is the $j$-th component of the vector $\alpha_{n,m}$.
\begin{lemma}\label{lem:cross_cc_statemnt_step_representation}
Regard a fixed $n \in \{0,\ldots,\NFE\!-\!1\}$ and a fixed $j \in \C$.
If any $\alpha_{n,m,j}$ with $m \in \{1,\ldots, \Nstg\}$ is positive, then all $\lambdan_{n,m',j}$ with $m'\in \{0,\ldots, \Nstg\}$ must be zero.
Conversely, if any $\lambdan_{n,m',j}$ is positive, then all $\alpha_{n,m,j}$ are zero.
\end{lemma}
\textit{Proof.} Let $\alpha_{n,m,i}$ be positive, and suppose $\lambdan_{n,j,i} = 0 $ and $\lambdan_{n,k,i}>0$ for some $k,j\in \{0,\ldots, \Nstg\}, k\neq j$, then $\alpha_{n,m,i}\lambdan_{n,k,i} >0$, which violates \eqref{eq:cross_cc_true_step_reformulation}, thus all $\lambdan_{n,m',i}=0,\ m'\in \{0,\ldots, \Nstg\}$.
The converse is proven similarly. \qed

An according statement holds for $\lambdap_{n,m}$ and $(e-\alpha_{n,m})$.
\color{black}
\begin{lemma}\label{lem:cross_cc_statemnt_step_representation2}
	Regard a fixed $n \in \{0,\ldots,\NFE\!-\!1\}$ and a fixed $j \in \C$.
	If any $1-\alpha_{n,m,j}$ with $m \in \{1,\ldots, \Nstg\}$ is positive, then all $\lambdap_{n,m',j}$ with $m'\in \{0,\ldots, \Nstg\}$ must be zero.
	Conversely, if any $\lambdap_{n,m',j}$ is positive, then all $1-\alpha_{n,m,j}$ are zero.
\end{lemma}
\color{black}


It is now left to discuss why the boundary points $\lambdap_{n+1,0} = \lambdap_{n,\Nstg}$ and  $\lambdan_{n+1,0} = \lambdan_{n,\Nstg}$ of the previous finite element are included in the cross complementarity conditions \eqref{eq:cross_cc_true_step_reformulation}.
It turns out, they are the key to switch detection.
A consequence of Lemmata \ref{lem:cross_cc_statemnt_step_representation} and \ref{lem:cross_cc_statemnt_step_representation2} is that, if the active set changes in the $j$-th component between the $n$-th and $(n+1)$-th finite element, then it must hold that $\lambdap_{n,\Nstg,j} =   \lambdap_{n+1,0,j} = 0$ and $\lambdan_{n,\Nstg,j} =   \lambdan_{n+1,0,j} = 0$.
Since $x_n^{\mathrm{next}} = x_{n+1}$, we have from \eqref{eq:dcs_irk_single_step_reformulation} the condition
\begin{align*}
	\switchfun_j(x_{n+1}) = 0,
\end{align*}
which defines exactly the switching surface between two regions.
Therefore, we have implicitly a constraint that forces $h_n$ to adapt such that the switch is detected exactly.

\color{black}
Given $a,b \in \R^p$, the complementarity conditions $0\leq a \perp b \geq 0$ mean that $a,b \geq 0$ and $a_i b_i = 0, i = 1,\ldots, p$. 
Due to the non-negativity of $a$ and $b$, the last conditions can be replaced by $a^\top b = 0$. 
Similar aggregations can be made with the cross complementarity conditions.
For clarity, we stated \eqref{eq:cross_cc_true_step_reformulation} in their most sparse form, without any aggregation.
\color{black}
However, the nonnegativity of $\alpha_{n,m},\lambdap_{n,m}$ and $\lambdan_{n,m}$ allows more compact forms. 
In the sequel, we use a formulation such that, together with the constraint $\sum_{n=0}^{\NFE-1} h_n = T$,  we have the same number of new equations as new degrees of freedom by varying $h_n$.
Thus, we combine constraints of two neighboring finite elements and have the compact formulation
\begin{align}\label{eq:cross_comp_step_formulation}
		&G_{\mathrm{cross}}(\mathbf{A},\mathbf{\Lambda^{\mathrm{p}}},\mathbf{\Lambda^{\mathrm{n}}}) = 0,
\end{align}
whose entries are for all $n \in \{0,\ldots, \NFE-2\}$ given by
\begin{align*}
	G_{\mathrm{cross},n}&(\mathbf{A},\mathbf{\Lambda^{\mathrm{p}}},\mathbf{\Lambda^{\mathrm{n}}})=	
	\sum_{k=n}^{n+1} \!
	\Big(
	\sum_{m = 1}^{\Nstg}\!\sum_{\substack{m'=0,\\m'\neq m}}^{\Nstg}\alpha_{k,m}^\top \lambdan_{k,m'}
	+(e-\alpha_{k,m})^\top \lambdap_{k,m'}		\Big)\!.
\end{align*}
We remind the reader that we use this seemingly complicated form to obtain a square system of equations.
This simplifies the study of the well-posedness of the FESD equations later.
However, in an implementation one can use any of the equivalent more sparse, or dense, formulations.
Many possible variants are implemented in \nosnoc~\cite{Nurkanovic2022b}, and the user can control the sparsity.

\subsection{Step size equilibration}\label{sec:fesd_step_equilibration_step_reformulation}
To complete the derivation of the FESD method for the Heaviside step reformulation \eqref{eq:dcs_step_2_restate}, we need to derive the step equilibration conditions.
Here, \textit{step} refers to the integration step size and should not to be confused with the set-valued Heaviside step function.

If no active-set changes happen, the cross complementarity constraints \eqref{eq:cross_cc_true_step_reformulation} are implied by the standard complementarity conditions~\eqref{eq:fesd_standard_comp_step_reformulation}.
\color{black}
This can easily be verified by looking at the stage point in the middle plot of Figure \ref{fig:fesd_motivation}.
\color{black}
Therefore, we end up with a system of equations with more degrees of freedom than conditions.
The step equilibration constraints aim to remove the degrees of freedom in the appropriate $h_n$ if no switches happen.
This results in a piecewise uniform discretization grid for the differential and algebraic states on the considered time interval.

We achieve the goals outlined above via the equation:
\begin{align}\label{eq:step_eq_step_reformulation}
	0&= G_{\mathrm{eq}}(\mathbf{h},\mathbf{A},\mathbf{\Lambda^{\mathrm{p}}},\mathbf{\Lambda^{\mathrm{n}}})
	\coloneqq 
	\begin{bmatrix}
		(h_{1}-h_{0})\eta_1(\mathbf{A},\mathbf{\Lambda^{\mathrm{p}}},\mathbf{\Lambda^{\mathrm{n}}}) \\
		\vdots\\
		(h_{\NFE\!-\!1}-h_{\NFE\!-\!2})\eta_{\NFE\!-\!1}(\mathbf{A},\mathbf{\Lambda^{\mathrm{p}}},\mathbf{\Lambda^{\mathrm{n}}})
	\end{bmatrix},
\end{align}
where $\eta_{n}$ is an indicator function that is zero only if a switch occurs, otherwise its value is strictly positive.
This provides a condition that removes the spurious degrees of freedom.
In the remainder of this section, we derive a possible expression for $\eta_{n}$.

The derivations below are motivated by the following facts.
Let $t_n$ be a switching point with $\switchfun_j(x(t_n)) =0$ for some $j \in \C$.
Consequently, it holds that $\lambdap_j(t_n) = \lambdan_j(t_n) =0$.
If, for example, a switch occurs at $t_n$ such that $\switchfun(x(t_n^-))<0$ and $\switchfun(x(t_n^+))>0$, we have that
$\dot{\lambda}^{\mathrm{n}}_j(t_n^-)<0$, $\dot{\lambda}^{\mathrm{n}}_j(t_n^-)=0$ and that $\dot{\lambda}^{\mathrm{p}}_j(t_n^+)=0$, $\dot{\lambda}^{\mathrm{p}}_j(t_n^+)>0$.
This can be verified by looking at the plots in Figure \ref{fig:switcing_cases_step}.
The symmetric case is possible as well.
The absolute values of these directional derivatives help us to encode the switching logic.

Now, instead of looking at the time derivatives, in the discrete-time case, we exploit the non-negativity of  $\lambdap_{n,m}$, $\lambdan_{n,m}$, and the fact no switches occur within a finite element.
For $n \in \{1,\ldots,\NFE-1\}$, we define the following backward and forward sums of the stage values over the neighboring finite elements $[t_{n-1},t_n]$ and $[t_{n},t_{n+1}]$:
\begin{align*}
	\sigma_{n}^{\lambdap,\mathrm{B}} &=  \sum_{m=0}^{\Nstg}  \lambdap_{n-1,m},\;\;
	\sigma_{n}^{\lambdap,\mathrm{F}} =  \sum_{m=0}^{\Nstg}  \lambdap_{n,m},\\
	\sigma_{n}^{\lambdan,\mathrm{B}} &=  \sum_{m=0}^{\Nstg}  \lambdan_{n-1,m},\;\;
	\sigma_{n}^{\lambdan,\mathrm{F}} =  \sum_{m=0}^{\Nstg}  \lambdan_{n,m}.
\end{align*}
They are zero if the left and right time derivatives are zero, respectively.
Likewise, they are positive when the left and right time derivatives are nonzero.


\begin{table}[t]	
	\centering
	\caption{Overview of switching cases for the step size equilibration}
	\centering
		\begin{tabular}{@{}l|l l|l l | l l | l @{}}	
			\specialrule{.15em}{0em}{0.0em} 
			\textbf{Switching case} & $	\sigma_{n}^{\lambdan,\mathrm{B}} $ &$	\sigma_{n}^{\lambdan,\mathrm{F}} $ & $	\sigma_{n}^{\lambdap,\mathrm{B}} $ & $	\sigma_{n}^{\lambdap,\mathrm{F}}$ &$	\pi_{n}^{\lambdan} $ & $\pi_{n}^{\lambdap}$  & $\upsilon_n$\\
			\hline
			No switch &  1 & 1 & 0 & 0 & 1 & 0 &1 \\
			Crossing &  1 & 0 & 0 & 1 & 0 & 0 &0 \\
			Entering sliding mode &  1 & 0 & 0 & 0 & 0 & 0 &0 \\
			Leaving sliding mode &  0 & 0 & 0 & 1 & 0 & 0 &0 \\
			Spontaneous switch &  0 & 0 & 0 & 1 & 0 & 0 &0\\
			\specialrule{.15em}{0em}{0.0em} 
		\end{tabular}
		\label{tab:step_size_equilibration}
	\end{table}
Moreover, for all $n \in \{ 1,\dots,\NFE-1$\} we define the following variables to summarize the logical dependencies:
\begin{align*}
	\pi_{n}^{\lambdan}  &=
	\mathrm{diag}(\sigma_{n}^{\lambdan,\mathrm{B}})	\sigma_{n}^{\lambdan,\mathrm{F}} \in \R^{\nswitchfun},\\
	\pi_{n}^{\lambdap}  &=
	\mathrm{diag}(\sigma_{n}^{\lambdap,\mathrm{B}})	\sigma_{n}^{\lambdap,\mathrm{F}} \in \R^{\nswitchfun},
\end{align*}
and 
\begin{align*}
	\upsilon_n  = \pi_{n}^{\lambdan} +\pi_{n}^{\lambdap} \in \R^{\nswitchfun}.
\end{align*}	
The switching cases and sign logic is summarized in Table \ref{tab:step_size_equilibration}. 
For readability, we put in the table a one if a variable is positive and a zero if it is zero. 
Let us discuss how the variables above encode the switching logic, and for this purpose, we regard the $j$-th switching functions $\switchfun_j(x)$.
If no switch occurs, and for example, we have that $\switchfun_j(x(t))<0$ during the regard time interval, it follows that $\lambdan_j(t) >0$ and $\lambdap_j(t) = 0$ during this time interval.
In the discrete time setting, we have $\sigma_{n,j}^{\lambdan,\mathrm{B}}, \sigma_{n,j}^{\lambdan,\mathrm{F}} >0$ and 
$\sigma_{n,j}^{\lambdap,\mathrm{B}} = \sigma_{n,j}^{\lambdap,\mathrm{F}} = 0$.
This means that $\pi_{n,j}^{\lambdan} > 0$, $\pi_{n,j}^{\lambdan} = 0$ and $\upsilon_{n,j} >0$.
It can be seen that the symmetric case with $\switchfun_j(x(t))>0$ leads also to $\upsilon_{n,j} >0$, hence we do not enumerate all symmetric cases in Table \ref{tab:step_size_equilibration}.

On the other hand, if we have a switch of the crossing type (cf. top plots in Figure \ref{fig:switcing_cases_step} with $\switchfun_j(x(t)) <0$ for $t<\tsn$ and  $\switchfun_j(x(t)) >0$ for $t>\tsn$, it follows that
$\lambdan_j(t) >0, \lambdap_j(t) =0$ for $t<\tsn$ and $\lambdan_j(t) =0, \lambdap_j(t) >0$ for $t>\tsn$.
In the discrete-time setting we obtain the sign pattern as in the second row of Table \ref{tab:step_size_equilibration}, with $\upsilon_{n,j} = 0$. 

In general, if there is an active-set change in the $j$-th complementarity pair, then at most one of the $j$-th components of $\sigma_{n}^{\lambdap,\mathrm{B}}$ and $\sigma_{n}^{\lambdap,\mathrm{F}}$, or
$\sigma_{n}^{\lambdan,\mathrm{B}}$ and $\sigma_{n}^{\lambdan,\mathrm{F}}$ is nonzero. 
In these cases, we obtain that $\upsilon_{n,j} = 0$, and if now switch happens we have that $\upsilon_{n,j} >0$.

In other words, $\upsilon_{n,j}$ is only zero if there is an active-set change in the $j$-th complementarity pair at~$t_n$, otherwise, it is strictly positive.
We summarize all logical relations for all switching functions into a single scalar expression and define
\begin{align*}
	\eta_n(\mathbf{A},\mathbf{\Lambda^{\mathrm{p}}},\mathbf{\Lambda^{\mathrm{n}}}) &\coloneqq \prod_{i=1}^{\nswitchfun} (\upsilon_n)_i.
\end{align*}
It is zero only if an active-set change happens at the boundary point $t_{n}$, otherwise, it is strictly positive.
\subsection{The FESD discretization}\label{sec:fesd_all_fesd_eq_step}
We have now introduced all extensions needed to pass from a standard RK discretization \eqref{eq:dcs_irk_step_reformulation} to the FESD discretization for the step reformulation.
With a slight abuse of notation, we collect all equations in a discrete-time system form:
\begin{subequations} \label{eq:fesd_compact_step_representation}
	\begin{align}
		s_{1} \!&=\! F_{\fesd}(\mathbf{Z}), \label{eq:fesd_compact_state_transition_step_representation}\\
		0 \!&= \!G_{\fesd}(\mathbf{Z},\mathbf{h},s_0, q , T),
	\end{align}
\end{subequations}
where $F_{\fesd}(\mathbf{x})\!=x_{\NFE}$ is the state transition map and  $G_{\fesd}(\mathbf{x},\mathbf{h},\mathbf{Z},q, T)$ collects all other internal computations including all RK steps within the regarded time interval:
\begin{align}\label{eq:fesd_compact_algebraic_step_representation}
	&G_{\fesd}(\mathbf{Z},\mathbf{h},s_0,q, T)\coloneqq
	\begin{bmatrix}
		{G}_{\mathrm{std}}(\mathbf{Z},\mathbf{h},s_0,q,T)\\
		G_{\mathrm{cross}}(\mathbf{A},\mathbf{\Lambda^{\mathrm{p}}},\mathbf{\Lambda^{\mathrm{n}}})\\
		G_{\mathrm{eq}}(\mathbf{h},\mathbf{A},\mathbf{\Lambda^{\mathrm{p}}},\mathbf{\Lambda^{\mathrm{n}}})\\
		\sum_{n=0}^{\NFE-1} h_n - T
	\end{bmatrix}.
\end{align}
Here, the control variable $q$, the horizon length $T$, and the initial value $s_0$ are given parameters.
\paragraph*{Remark on RK methods with $c_{\Nstg}\neq 1$}
The extension for the case of an RK method with $c_{\Nstg}\neq1$ follows similar lines as in Stewart's formulation~\cite{Nurkanovic2024a}.
We have that $t_n+c_{\Nstg} h_n < t_{n+1}$.
Hence, the variables $\lambdap_{n,\Nstg}$ and $\lambdan_{n,\Nstg}$ do not correspond to the boundary values $\lambdan(t_{n+1})$ and $\lambdap(t_{n+1})$ anymore.
We denote the true boundary points by $\lambdap_{n,\Nstg+1}$ and $\lambdan_{n,\Nstg+1}$.
They can be computed from the KKT conditions of the step reformulation LP~\eqref{eq:step_parametric_lp}.
For all $n \in \{1,\ldots,\NFE-2 \}$ we have
\begin{align}\label{eq:lp_boundary_points_step_representatio}
	\begin{bmatrix}
		\switchfun(x_{n+1}) - \lambdap_{n,\Nstg+1} + \lambdan_{n,\Nstg+1}\\
		\Psi(\lambdan_{n,\Nstg+1},\alpha_{n,\Nstg+1})\\
		\Psi(\lambdap_{n,\Nstg+1},e-\alpha_{n,\Nstg+1})
	\end{bmatrix} = 0.
\end{align}
These equations are appended to the FESD equation in~\eqref{eq:fesd_compact_algebraic_step_representation}.

However, to make the switch detection work, we must update the continuity conditions for the discrete-time versions of the Lagrange multipliers and adapt the cross-complementarity conditions accordingly.
For all $n = \{0,\ldots, \NFE-1\}$, \eqref{eq:continuity_of_lambda_step_representation} is replaced by:
\begin{align}\label{eq:continuity_of_lambda_2_step_representation}
	\lambdap_{n,\Nstg+1}= \lambdap_{n+1,0},\; 	\lambdan_{n,\Nstg+1}= \lambdan_{n+1,0}.
\end{align}
We append to the vectors $\mathbf{A},\mathbf{\Lambda^{\mathrm{p}}}$ and $\mathbf{\Lambda^{\mathrm{n}}}$ the new variables $\alpha_{n,\Nstg+1},\lambdap_{n,\Nstg+1}$ and $\lambdan_{n,\Nstg+1}$ accordingly.
For the whole control interval, we have in total $3(\NFE-1)n_c$ new variables.
It is only left to state the modified cross complementarity conditions, including the expressions' $(\Nstg+1)$-th point.
More explicitly, the $n$-th component of \eqref{eq:cross_comp_step_formulation}
reads now for all $n \in \{0,\ldots, \NFE-2\}$ as
\begin{align*}
	G_{\mathrm{cross},n}&(\mathbf{A},\mathbf{\Lambda^{\mathrm{p}}},\mathbf{\Lambda^{\mathrm{n}}})=\ \sum_{k=n}^{n+1}
	\sum_{m = 1}^{\Nstg}\! \sum_{\substack{m'=0,\\m'\neq m}}^{\Nstg+1}\alpha_{k,m}^\top \lambdan_{k,m'}
	+(e-\alpha_{k,m})^\top \lambdap_{k,m'}.
\end{align*}
\subsection{Discretizing optimal control problems with FESD}\label{sec:fesd_ocp}
Regard a optimal control problem subject to a nonsmooth dynamical system \eqref{eq:basic_di}:
 of the following form:
\begin{subequations} \label{eq:ocp}
	\begin{align}
		\min_{x(\cdot),u(\cdot)} \quad & \int_{0}^{T} L(x(t),u(t))\dd t +  R(x(T)) \label{eq:ocp_cost}\\
		\textrm{s.t.} \quad  x(0) &= s_0, \\
		\dot{x}(t) &\! \in \! \F(x(t),u(t),\Gamma(\switchfun(x(t)))),& \textrm{for a.a. } t \in [0,\!T], \label{eq:ocp_pss} \\
		0&\geq G_{\mathrm{path}}(x(t),u(t)),& t \in [0,T],\\
		0&\geq G_{\mathrm{terminal}}(x(T)),
	\end{align}
\end{subequations}
where $L: \R^{n_x} \times \R^{n_u} \to \R$ is the running cost and $R:\R^{n_x}\to \R$ is the terminal cost, $s_0\in\R^{n_x}$ is a given initial value.
The path and terminal constraints are grouped into the functions $G_{\mathrm{path}} : \R^{n_x}  \times \R^{n_u} \to \R^{n_{G_{\mathrm{p}}}}$ and $G_{\mathrm{terminal}} : \R^{n_x}  \to \R^{n_{G_{\mathrm{t}}}}$, respectively.

\color{black}
In this paper, we consider a direct approach~\cite{Rawlings2017}, i.e., we first discretize the continuous-time OCP \eqref{eq:ocp} and then solve a finite-dimensional nonlinear program (NLP).
Here we discretize the OCP using the FESD method.
\color{black}
First, for the discretization of the control function, we consider $\Nctrl\geq 1$ control intervals of equal length, indexed by $k$.
We use a piecewise constant control discretization, where the control variables are collected 
$\mathbf{q} = (q_0,\ldots,q_{\Nctrl-1})\in \R^{\Nctrl n_u}$.
Such a discretization is typically used in feedback control, but extensions to more sophisticated control are straightforward.
All internal variables are additionally indexed by $k$.

\color{black}
Second, we discretize the cost function \eqref{eq:ocp_cost} and the dynamics \eqref{eq:ocp_pss}.
To apply the FESD method, the differential inclusion in \eqref{eq:ocp_pss} is transformed into an equivalent DCS, as described in Section \ref{sec:step_functions}.
The DCS, now of the from of \eqref{eq:step_dcs}, is discretized with the FESD method summarized in the previous section.
For each control interval $k$ we use \eqref{eq:fesd_compact_step_representation} with $N_{\mathrm{FE}}$ internal finite elements.
\color{black}
The state values at the control interval boundaries are grouped in the vector  $\mathbf{s} = (s_0,\ldots,s_{\Nctrl})\in\R^{(\Nctrl+1)n_x}$.
In ${\mathcal{Z}} = (\mathbf{{Z}}_0,\ldots,\mathbf{{Z}}_{\Nctrl-1})$ we collect all internal variables and in $\mathcal{H} = (\mathbf{h}_0,\ldots,\mathbf{h}_{\Nctrl-1})$ all step sizes.
\color{black}
For the cost discretization, one can derive a quadrature formula~\cite[Chapter 8]{Rawlings2017}, or introduce a scalar quadrature state $\dot{\ell}(t) = L(x(t),u(t)), \ell(0) = 0$, integrate it together with the dynamics equations, and use $\ell(T)$ in the objective, which approximates the integral term in~\eqref{eq:ocp_cost}.
\color{black}

\color{black}
Third, the path constraints are relaxed and evaluated only at the control interval boundary points, i.e., at the variables $s_k$ and $q_k, k = 0,\ldots,\Nctrl-1$.
If necessary, the path constraint can be evaluated on a finer grid using the values computed in the internal integration intervals within a control interval or at the RK stage points.
The terminal constraint and cost are simply evaluated at $s_N$, which is an approximation of $x(T)$. 
In summary, we obtain a discrete-time variant of \eqref{eq:ocp}:
\color{black}
\begin{subequations}\label{eq:ocp_discrete_time}
	\begin{align}
		\min_{\mathbf{s},\mathbf{q},\mathcal{Z},\mathcal{H}} \quad & \sum_{k=0}^{\Nctrl-1} \hat{L}(s_k,\mathbf{x}_k,q_k)+ {R}(s_{\Nctrl}) \\
		\textrm{s.t.} \quad  &s_{0} = \bar{x}_0,\\
		&{s}_{k+1}  = F_{\fesd}(\mathbf{x}_k),&k = 0,\ldots,\Nctrl-1,\\
		&0 = G_{\fesd}(\mathbf{x}_k,\mathbf{Z}_k,q_k), &k = 0,\ldots,\Nctrl-1,\\
		&0 \geq G_{\mathrm{ineq}}(s_k,q_k), &k = 0,\ldots,\Nctrl-1,\\
		&0 \geq G_{\mathrm{terminal}}(s_{\Nctrl}),
	\end{align}
\end{subequations}
where $\hat{L}:\R^{n_x}\times \R^{(\NFE+1)\Nstg n_x} \times \R^{n_u}\to \R$ is the discretized running cost.

Due to the complementarity constraints in the FESD discretization, this NLP is a mathematical program with complementarity constraints.
This are degenerate NLPs since the complementarity constraints lead to the violation of standard constraint qualifications at all feasible points.
In practice, they can often be efficiently solved by solving a sequence of related and relaxed NLPs within a homotopy approach.
Such an approach with some of the standard reformulations \cite{Scholtes2001,Anitescu2007,Ralph2004} is implemented in \nosnoc.
\color{black}
A survey and comparison of state-of-the-art methods for solving NLPs in nonsmooth optimal control problems is given in \cite{Nurkanovic2024b}.
All these homotopy solution methods are implemented in~\nosnoc~\cite{Nurkanovic2022b}. 
In practice, Scholtes' relaxations~\cite{Scholtes2001} together with \texttt{IPOPT}~\cite{Waechter2006}, called via its \texttt{CasADi} interface~\cite{Andersson2019}, often work very well.
We use this method in the numerical experiments in the paper, when solving problems like \eqref{eq:ocp_discrete_time}.
\color{black}
\section{Convergence theory of FESD for the step representation}\label{sec:fesd_step_theory}
In this section, we present the main convergence result of the FESD method for the step representation outlined in~\eqref{eq:fesd_compact_step_representation}.
Specifically, we show that: 
(1) the solutions to the FESD problem are locally isolated; 
(2) both the solution approximations and 
(3) the numerical sensitivities obtained via the FESD method converge with the same order of accuracy as the underlying RK method. 
The proofs are similar to those used in Stewart's case in \cite{Nurkanovic2024a}, hence, we will not repeat them here.
The main difference is in the assumptions we make. 
\subsection{Main assumptions}\label{sec:fesd_assumptions_step}
We start by stating all assumptions.
The first assumption relates to the underlying RK methods:
\begin{assumption}(Runge-Kutta method)\label{ass:irk_scheme_step}
	A Butcher tableau with the entries $a_{i,j} ,b_i$ and $c_i$, $i,j\in\{1,\ldots,\Nstg\}$ related to an $\Nstg$-stage Runge-Kutta (RK) method 
	is used in the FESD \eqref{eq:fesd_compact_step_representation}. 
	Moreover, we assume that:
	\begin{enumerate}[(a)]
		\item If the same RK method is applied to the differential algebraic equation \eqref{eq:step_dae} on an interval $[t_a,t_b]$, with a fixed active set, it has a global accuracy of $O(h^p)$ for the differential states.
		\item The RK equations applied to \eqref{eq:step_dae} have a locally isolated solution for a sufficiently small $h_n>0$. 
	\end{enumerate}
\end{assumption}
\color{black}
The purpose of this assumption is to describe the properties of an underlying RK method used in FESD when applied to a smooth ODE or DAE. 
Both requirements (a) and (b) are standard and are satisfied by many RK methods used in practice, e.g. Radau IIA or Gauss-Legendre methods, cf. \cite{Hairer1991}.
\color{black}
The second assumption concerns the existence of solutions to the FESD problem outlined in \eqref{eq:fesd_compact_step_representation}.
\begin{assumption}(Solution existence)\label{ass:solution_existence_fesd_step}
	For given parameters $s_0,q$ and $T$, there exists a solution to the FESD problem \eqref{eq:fesd_compact_step_representation}, such that for all $n \in \{0,\ldots,\NFE-1\}$ it holds that ${h}_n>0$.
\end{assumption}
\color{black}
The problem \eqref{eq:fesd_compact_step_representation} is a nonlinear complementarity problem \cite{Facchinei2003}.
The proof of existence of solutions for the standard RK equations \eqref{eq:dcs_irk_step_reformulation} can probably be done with standard tools from complementarity theory \cite{Facchinei2003}, as it was done for example in \cite[Proposition 15]{Acary2014} for the implicit Euler method.
With the additional cross complementarity and step equilibration conditions appearing in \eqref{eq:fesd_compact_step_representation}, the same proof technique is no longer applicable.
The existence of solutions is still an open problem. 
Motivated by empirical observations, we assume the existence of solutions in this paper.
\color{black}

The next assumption is slightly more technical and relates to the regularity of the problem under consideration.
\begin{assumption}(Regularity)	\label{ass:regularity_step}
	Given the complementarity pairs $\Psi(\alpha_{n,m},\lambdan_{n,m})=0$ and $\Psi(e-\alpha_{n,m},\lambdap_{n,m})=0$, for all $n \in \{ 0,\ldots\NFE-1\}$ there exists an $m \in \{1,\dots,\Nstg\}$ and $i \in \{1,\ldots,\Nsys\}$, such that the strict complementarity property holds, i.e., $\alpha_{n,m,i}+\lambdan_{n,m,i}>0$ and $e-\alpha_{n,m,i}+\lambdap_{n,m,i}>0$. 
	Moreover, for the RK equations \eqref{eq:dcs_irk_single_step_reformulation}, it holds for all $n \in \{0,\ldots\NFE-1\}$, that at least one entry of the vector $\nabla_{h_n} G_{\irk}(x_{n+1},Z_n,h_n,q)$ is nonzero.
\end{assumption}
\color{black}
This assumption is made to ensure the correct rank of partial Jacobians of \eqref{eq:fesd_compact_step_representation} (with a fixed active set). 
It is used to prove the local uniqueness of solutions of the FESD problem, and similar assumptions are made in \cite{Nurkanovic2024a}.
The first part of the assumption requires that at least one complementarity pair on stage points within an integration interval (finite element) satisfies strict complementarity. 
Looking at the common switching cases in Figure \ref{fig:switcing_cases_step}, one can see that this assumption always holds and is thus not restrictive.
The second part requires that at least one term, $G_{\irk}(x_{n+1},Z_n,h_n,q)$ multiplied by $h_n$, is nonzero.
Both assumptions can be checked computationally once a candidate solution has been computed.
\color{black}

Before starting the final assumption, let us introduce some notation. 
We use a suitable interpolation scheme to construct a continuous-time approximation of the solution based on the stage values obtained by solving equation \eqref{eq:fesd_compact_step_representation}.
\color{black}
We denote the continuous-time approximation on every finite element $[t_n,t_{n+1}]$ by $\hat{x}_n(t;h_{n})$.
To approximate the solution over the entire time interval $[0,T]$, we append the local approximations from each finite element:
\begin{align}\label{eq:continious_time_fesd_step}
	\hat{x}_h(t) &= \hat{x}_n(t;h_{n})\ \text{if}\ t\in [t_n,t_{n+1}],
\end{align}
where $h = \max_{n\in\{0,\ldots \NFE-1\}} h_n$. 
The subscript $h$ for a solution approximation $\hat{x}_h(t)$ indicates that we get different solution approximations on the entire interval $[0,T]$ as the maximum step size changes, and not that the entire approximation is parameterized by a single step size $h$. 
As the maximum step size shrinks, we expect $\hat{x}_h(t)$ to converge to an exact solution $x(t)$.
\color{black}
The set of all grid points is defined as $\mathcal{G} = \{t_0,\ldots, t_{\NFE}\}$.
We can also use this approach to construct continuous-time representations for the algebraic variables, which we denote by 
$\hat{\theta}_h, \hat{\alpha}_h$, $\hat{\lambda}_{h}^{\mathrm{n}}$ and $\hat{\lambda}_{h}^{\mathrm{p}}$.
The $n$-th switching point of the true solution is denoted by ${t}_{\mathrm{s},n}$ and one corresponding to a solution approximation by $\hat{t}_{\mathrm{s},n}$.
Similarly, the active sets (cf. \eqref{eq:active_set_defintion}) of the solution approximation are denoted by 
$\I(\hat{x}_h(t)) = \hat{\I}_n,\ t\in(\hat{t}_{\mathrm{s},n},\hat{t}_{\mathrm{s},n+1})$ and the active set at switching point $\hat{t}_{\mathrm{s},n}$ by $\I(\hat{x}_h(\hat{t}_{\mathrm{s},n})) = \hat{\I}_n^0$.
\begin{assumption}\label{ass:active_sets_steps}
	Let $\I^{0}_{n}$ be the active set at switching point ${x}({t}_{\mathrm{s},n})$ of true solution and 
	$\hat{\I}^{0}_{n}$ the active set at the switching point $\hat{x}(\hat{t}_{\mathrm{s},n})$ a solution approximation.
	If $\hat{t}_{\mathrm{s},n}$ is sufficiently close to ${t}_{\mathrm{s},n}$ and $\I^{0}_{n} = \hat{\I}^{0}_{n}$, then
	$\I_{n+1} = \hat{\I}_{n+1}$.
	Furthermore, if there are several possible new active sets, they are identical for both the true solution and its approximation.
\end{assumption}
\color{black}
This assumption requires that a given solution approximation and the corresponding true solution predict the same active sets. 
In other words, the solution approximation enters the same region or sliding mode as the true solution after a switching event.
In Stewart's reformulation, this statement can be directly proved using an auxiliary linear complementarity problem constructed with the problem data~\cite{Stewart1990a}.
In the case of DCS \eqref{eq:step_dcs}, such auxiliary problems are not available and the property is assumed directly.
The requirement that a sufficiently good solution approximation predicts the same active set as the true solution is needed to prove the high accuracy convergence of a switch detection method such as FESD, cf. \cite[Theorem 4.3 ]{Stewart1990a} and \cite[Theorem 16]{Nurkanovic2024a}.
\color{black}
\subsection{Solutions of the FESD problem are locally isolated}\label{sec:fesd_step_isolated_solutions}
In this section, we present a theorem on the regularity of solutions to the FESD problem~\eqref{eq:fesd_compact_step_representation}.
For brevity, we focus on the case where $c_{\Nstg} = 1$. 
For the reader's convince we restate the FESD problem
\begin{align}\label{eq:fesd_equation_step}
	&G_{\fesd}(\mathbf{Z},\mathbf{h},s_0,q, T)\! \coloneqq\!\!
	\begin{bmatrix}
	{G}_{\mathrm{std}}(\mathbf{Z},\mathbf{h},s_0,q,T)\\
	G_{\mathrm{cross}}(\mathbf{A},\mathbf{\Lambda^{\mathrm{p}}},\mathbf{\Lambda^{\mathrm{n}}})\\
	G_{\mathrm{eq}}(\mathbf{h},\mathbf{A},\mathbf{\Lambda^{\mathrm{p}}},\mathbf{\Lambda^{\mathrm{n}}})\\
	\sum_{n=0}^{\NFE-1} h_n - T			
\end{bmatrix}.
\end{align}
Furthermore, we provide a summary of the dimensions of all key functions and variables:
\begin{itemize}
	\item Degrees of freedom: $\mathbf{Z} =(\textbf{x},\mathbf{V},\mathbf{\Theta},\mathbf{A},\mathbf{\Lambda}^{\mathrm{p}},\mathbf{\Lambda}^{\mathrm{n}})$ and $\mathbf{h}$.
	\item Total number of degrees of freedom: $n_{\mathbf{Z}}+\NFE$, where $n_{\mathbf{Z}}= (\NFE+1)n_x + \NFE\Nstg n_x + n_{\theta}+3n_{\alpha}$.
	\item Dimension of $\mathbf{\Theta}$: $n_{{\theta}}= \NFE\Nstg\Nsys$.
	\item Dimension of $\mathbf{A},\mathbf{\Lambda}^{\mathrm{p}}, \mathbf{\Lambda}^{\mathrm{n}}$: $n_{\alpha} = \NFE\Nstg n_c$
	\item Parameters: $(s_0,q,T) \in \R^{n_x+n_u+1}$, 	
	\item Standard RK equations: ${G}_{\mathrm{std}}: \R^{n_{\mathbf{Z}}} \times \R^{\NFE} \times \R^{n_x} \times \R^{n_u} \times \R \rightarrow \R^{n_{\mathbf{Z}}}$, 
	\item Cross complementarity: $G_{\mathrm{cross}}:  \R^{n_{\alpha}} \times \R^{n_{\alpha}} \times \R^{n_{\alpha}} \rightarrow \R^{\NFE-1}$,
	\item Step equilibration: $G_{\mathrm{eq}}: \R^{n_{\NFE}} \times  \R^{n_{\alpha}}  \times \R^{n_{\alpha}} \times \R^{n_{\alpha}}   \rightarrow \R^{\NFE-1}$ and
	\item FESD equations: $G_{\fesd}: \R^{n_{\mathbf{Z}}} \times \R^{\NFE} \times \R^{n_x} \times \R^{n_u} \times \R \to \R^{n_{\mathbf{Z}} + 2\NFE-1}$.
\end{itemize}
The vectors $s_0 \in \R^{n_x}$, $q \in \R^{n_u}$ and $T\in \R$ are given parameters.
As a result, the system of equations \eqref{eq:fesd_equation_step} has a total of $n_{\mathbf{Z}} + \NFE$ unknowns and $n_{\mathbf{Z}} + 2\NFE-1$ equations, which means it is always over-determined.

\begin{theorem}\label{th:fesd_locally_unique_solution_step}
	Suppose that Assumptions \ref{ass:irk_scheme_step}, \ref{ass:solution_existence_fesd_step} and \ref{ass:regularity_step} hold. 
	Let ${s}_0$, ${q}_0$  and ${T}$ be some fixed parameters such that $G_\fesd(\mathbf{Z}^*,\mathbf{h}^*,{s}_0,{q},{T}) =0$ in \eqref{eq:fesd_equation_step}.
	Let $P^*\subseteq \R^{n_x}\times \R^{n_u} \times \R $ be the set of all parameters $(\hat{s}_0,\hat{q},\hat{T})$ such that $\mathbf{Z} \in \R^{n_\mathbf{Z}}$, which is the solution of $G_\fesd(\mathbf{Z},\mathbf{h},\hat{s}_0,\hat{q},\hat{T})=0$, has the same active set as $\mathbf{Z}^*$. 
	Additionally, suppose that $G_\fesd(\cdot)$ is continuously differentiable in $s_0,q$ and $T$ for all $(s_0,q,T) \in P^*$. 
	Then there exists a neighborhood ${P} \subseteq P^*$ of $({s}_0,{q}_0,{T})$ such that there exist continuously differentiable single valued functions  $\mathbf{Z}^{*}: {P} \to  \R^{n_{\mathbf{Z}}}$ and $\mathbf{h}^*: {P} \to \R^{\NFE}$.
\end{theorem}
\textit{Proof.} The proof follows similar lines as the proof of \cite[Theorem 14]{Nurkanovic2024a} and we omit it for brevity. \qed 

\color{black}
The main steps in the proof are as follows.
One observes that if there is an active-set change (i.e., a switch) between two adjacent finite elements, then the cross-complementarity conditions are implied by the standard complementarity conditions and are thus redundant and can be removed. 
The step size is locally determined by the step equlibration condition. 
Similarly, if there is a switch, then the cross complementarity conditions are binding and the corresponding step equilibration condition can be removed. 
In summary, the FESD problem \eqref{eq:fesd_equation_step} reduces to a square systems with $n_{\mathbf{Z}}+ \NFE$ equations and unknowns. 
Under the given assumptions, one can prove that the Jacobian of \eqref{eq:fesd_equation_step} is of full rank and apply the implicit function theorem to obtain the above result. 
\color{black}
\subsection{Convergence of the FESD method}
We proceed by stating the results of the convergence of the FESD method. 
We show that under suitable assumptions, the sequence of approximations $\hat{x}_h(\cdot)$ generated by the FESD method converges to a solution of \eqref{eq:pss}.
In particular, the FESD method has the same order as the underlying RK method for smooth ODE.

\begin{theorem}\label{th:integration_order_step}
	Let $x(\cdot)$ be a solution of \eqref{eq:pss} with finitely many active set changes for $t\in [0,T]$ with $x(0) =x_0$. 
	Suppose the following is true:
	\begin{enumerate}[(a)]
		\item The Assumptions \ref{ass:solution_existence_step} and \ref{ass:active_sets_steps} are satisfied. 
		\item The Assumptions \ref{ass:irk_scheme_step}, \ref{ass:solution_existence_fesd_step} and \ref{ass:regularity_step} hold for the FESD problem \eqref{eq:fesd_compact_step_representation}.
	\end{enumerate}
	Then $x(\cdot)$ is a limit point of the sequence of approximations $\hat{x}_h(\cdot)$, defined in Eq.~\eqref{eq:continious_time_fesd_step} as $h\downarrow 0$.
	Moreover, for sufficiently small ${h}>0$, the solution of \eqref{eq:fesd_compact_step_representation} generates a solution approximation $\hat{x}_h(t)$ on $[0,T]$ such that:
	\begin{subequations}\label{eq:fesd_convergence_step}
		\begin{align}
			|\tsnhat- \tsn | &= O({h}^p)\ \text{for every } n\in\{0,\dots,\Nswitch\}, \label{eq:fesd_convergence_t}\\
			\| \hat{x}_h(t) - x(t) \| &= O({h}^p), \ \text{for all } t \in\G.  \label{eq:fesd_convergence_x}
		\end{align}
	\end{subequations}
\end{theorem}
\textit{Proof.} 
The proof follows similar lines as the proof of \cite[Theorem 16]{Nurkanovic2024a}. 
The primary distinction lies in the prediction of new active sets.
In \cite[Theorem 16]{Nurkanovic2024a}, we use \cite[Assumption 8]{Nurkanovic2024a} to be able to apply \cite[Lemma A.2]{Stewart1990a} and demonstrate that both the approximation and the exact solution share the same active set in the vicinity of a switching point. 
In this theorem, the assertion emerges directly from Assumption \ref{ass:active_sets_steps}.
\qed

\color{black}

The proof of \cite[Theorem 16]{Nurkanovic2024a} is inspired by the proof of \cite[Theorem 4.3]{Stewart1990a} and is quite involved. The main idea is to consider intervals $[\tsn,\tsnn]$ with fixed active sets $\I_n$ where the dynamics are locally smooth, and we recover the accuracy of the underlying RK method. At a switching point, $\tsnhat$, respectively $\tsn$, one has to prove \eqref{eq:fesd_convergence_t} and that the solution approximation can continue to evolve with the same active set $\I_{n+1}$ as the true solution.
Then the argument can be used inductively.

Some distinction must be made for the case where the approximation switches before or after the true solution to obtain the results. 
Once \eqref{eq:fesd_convergence_t} is proved, \eqref{eq:fesd_convergence_x} can be obtained from some algeraic manipulations and the Lipschitz continuity of the local dynamics and its solution. 
The error is dominated by the maximum step size $h$; hence, it is used in the error estimate.
\color{black}
\subsection{Convergence of discrete-time sensitivities} \label{sec:sensitivity_convergence_step}
This section concludes by demonstrating that the numerical sensitivities (cf. Section \ref{sec:intro} for a definition) obtained using the FESD method for the step reformulation converge to the correct values with a high order of accuracy. 
We remind the reader that numerical sensitivities of standard time-stepping methods, e.g. \eqref{eq:dcs_irk_step_reformulation}, do not converge to the correct values, no matter how small the step size becomes~\cite{Stewart2010}.
The convergence of the sensitivities is crucial for the success of direct optimal control methods~\cite{Nurkanovic2020}.

Before stating the result, we assume the time derivatives of the solution approximation $\hat{x}_h$ converge accordingly. 
This assumption extends Assumption \ref{ass:irk_scheme_step} and allows us to consider a wide range of RK methods.
\begin{assumption}(RK derivatives) \label{ass:irk_scheme_derivative_step}
Regard the RK methods from Assumption \ref{ass:irk_scheme_step} applied to the differential algebraic equations \eqref{eq:step_dae}. 
The derivatives of the numerical approximation for the same RK method converge with order $1 \leq q\leq p$, i.e., $\| \dot{\hat{x}}_h(t) - \dot{x}(t) \| = O(h^q),\; t\in \G$.
\end{assumption}
We remind the reader that for collocation-based implicit RK methods for ODE in general it holds that $q = p-1$ \cite[Theorem 7.10]{Hairer1993}. 
\begin{theorem}[Convergence to exact sensitivities]\label{th:exact_sensitivites_step}
Suppose the assumptions of Theorem \ref{th:integration_order_step} and Assumption \ref{ass:irk_scheme_derivative_step} hold.
Assume that a single active-set change happens at time $t_{\mathrm{s},n}$, i.e., $||\I_n| - |\I_{n+1}|| \leq 1, n \in\{0,\dots, \Nswitch\}$.
 Then for $h\downarrow 0$ it holds that $\frac{\partial \hat{x}_h(t,x_0)}{\partial x_0}  \to  \frac{\partial {x}(t,x_0)}{\partial x_0}$ with the convergence rate
\begin{align}
	\| \frac{\partial \hat{x}_h(t,x_0)}{\partial x_0} - \frac{\partial {x}(t,x_0)}{\partial x_0}\| = O(h^{q}), \text{ for all } t\in \G.
\end{align}
\end{theorem}
\textit{Proof.} The proof is essentially the same as the proof of \cite[Theorem 18]{Nurkanovic2024a}, one has only to replace the local switching functions $\psi_{i,j}(x)$ by an appropriate switching function $\switchfun_k(x)$.\qed

\color{black}
The key to this proof is the fact that the cross complementarity conditions imply $\switchfun_j(\hat{x}_h(\hatts)) = 0$ at a switching point $\hatts$, cf. Section \ref{sec:fesd_cross_comp_step_reformulation}.
Most of the proof then applies the chain rule and the implicit function theorem to obtain the similar expression for the discrete-time numerical sensitivities as for the continuous-time case. 
\color{black}


\section{Efficient modeling with set-valued Heaviside step functions}
\label{sec:compact_rep_and_lifting}
In this section, we show how to efficiently represent common geometries of the PSS regions with the use of Heaviside step functions.
This is useful for reducing the complexity of the modeling process.
Moreover, we introduce a lifting algorithm, which makes the multi-affine expressions for $\theta_i$ in \eqref{eq:step_dcs_theta} ``less nonlinear'' with the help of auxiliary variables. 

\subsection{Overview of expressions for $\theta$ via Heaviside step functions}
We regard the following two sets:
$A = \{x \mid \switchfun_A(x) > 0\}, B = \{x \mid \switchfun_B(x) > 0 \}$, and let $\alpha_A \in \gamma(\switchfun_A(x))$ and $\alpha_B \in \gamma(\switchfun_B(x))$. 
Note that these set do not have to correspond to the base regions in Definition \ref{def:basis_sets}, but we can use them to define such regions.
Table \ref{tab:step_summary} provides an overview of how the elementary algebraic expressions for the multipliers $\theta_i$ are related to the geometric definition of a region $R_i$.
More complicated expressions can be obtained by combining the ones listed in Table~\ref{tab:step_summary}.

\begin{table}[t]
	\centering
	\caption{Expressions of $\theta_i$ for different definitions of $R_i$.}
	\begin{tabular}{@{}lll@{}}
		\specialrule{.15em}{1em}{0em} 
		\text{Definition of} $R_i$ & \text{Expression for} $\theta_i$ & \text{Sketch} \\
		\hline
		$R_i = A$ & $\theta_i =  \alpha_A$  & 
		\begin{minipage}{.35\textwidth}
			\vspace{0.05cm}
			\includegraphics[width=0.6\linewidth]{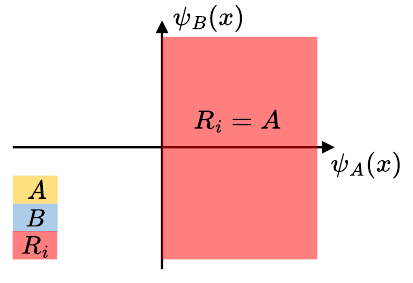}
		\end{minipage} 
		\\
		$R_i = A\cup B $& $\theta_i = \alpha_A+ \alpha_B$ & 			 \begin{minipage}{.35\textwidth}
			\vspace{0.02cm}
			\includegraphics[width=0.6\linewidth]{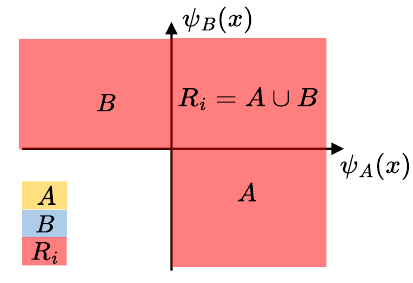}
		\end{minipage} 
		\\
		$R_i = A\cap B$ & $\theta_i =  \alpha_A  \alpha_B$& 			 \begin{minipage}{.35\textwidth}
			\vspace{0.02cm}
			\includegraphics[width=0.6\linewidth]{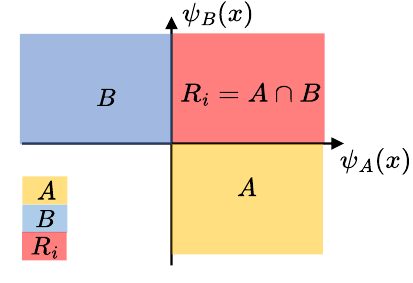}
		\end{minipage} 
		\\
		$R_i = \mathrm{int}(\R^{n_x} \setminus A) = \{x \mid \switchfun_A(x)<0\}$ & $\theta_i =1-\alpha_A$& 			 \begin{minipage}{.35\textwidth}
			\vspace{0.02cm}
			\includegraphics[width=0.6\linewidth]{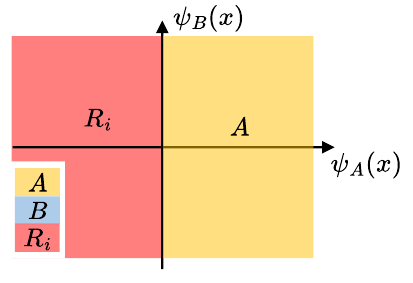}
		\end{minipage} 
		\\
		$R_i = A \setminus B$ & $\theta_i =\alpha_A- \alpha_B$ 			 
		&
		\begin{minipage}{.35\textwidth}
			\vspace{0.02cm}
			\includegraphics[width=0.6\linewidth]{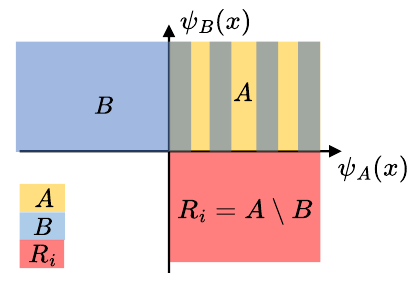}
		\end{minipage} 
		\\
		\specialrule{.15em}{0em}{0em} 
	\end{tabular}
	\label{tab:step_summary}
\end{table}

\begin{remark}[Sum of Filippov systems]
	In practice, one often encounters DIs that arise from the sum of several Filippov systems.
	This occurs, for example, if we have multiple surfaces with friction, or multiple objects touching the same frictional surface \cite{Stewart1996}.
	All developments from this paper can be extended to this case and are implemented in~\nosnoc.
	For the sake of brevity, we omit the corresponding equations and refer the reader to~\cite[Section 2.3]{Nurkanovic2024a}.
\end{remark}

\subsection{Representing unions of sets}\label{sec:unions_of_sets}
The regions $R_i$ may be given as unions of base sets $\tilde{R}_j$. 
Consequently, the number of multipliers $\theta \in \R^{\Nsys}$ decreases and it holds that $\Nsys < 2^{\nswitchfun}$.
For example, in the extreme case, we may have $R_1 = \tilde{R}_i$ and $R_2 = \cup_{j=1,j\neq i}^{2^{\nswitchfun}} \tilde{R}_j$, for some $i$, which significantly reduces the number of variables, since we are left with only two regions, i.e., $\Nsys =2$.
We illustrate such a case with by a simple example and discuss the general case in the sequel.
\begin{example}[Union of sets]\label{ex:union_of_sets}
	We regard an example with two scalar switching functions $\switchfun_1(x)$ and $\switchfun_2(x)$, with the basis regions
	$\tilde{R}_1 = \{ x\in \R^{n_x} \mid \switchfun_1(x) > 0, \switchfun_2(x)>0\},\; 
	\tilde{R}_2 = \{ x\in \R^{n_x} \mid \switchfun_1(x) > 0, \switchfun_2(x)<0\},\;
	\tilde{R}_3 = \{ x\in \R^{n_x} \mid \switchfun_1(x) < 0, \switchfun_2(x)>0\}$ and
	$\tilde{R}_4 = \{ x\in \R^{n_x} \mid \switchfun_1(x) < 0, \switchfun_2(x)<0\}$.
	The PSS is defined by the two regions:
	\begin{align*}
		R_1 &= \cup_{i=1}^{3} \tilde{R}_i = \{ x\in \R^{n_x} \mid \switchfun_1(x)>0\} \cup\{ x\in \R^{n_x} \mid \switchfun_1(x)\leq 0,\, \switchfun_2(x)>0 \},\\
		R_2 &= \tilde{R}_4= \{ x\in \R^{n_x} \mid \switchfun_1(x)< 0,\, \switchfun_2(x)<0 \}.
	\end{align*}
		\color{black}
	According to the second row of Table \ref{tab:step_summary}, if a region $R_i$ consists of the union of two sets, the expression for $\theta_i$ is the sum of two corresponding indicators.
	In the current example, $R_1$ consists of the union of three base regions, so the expressions for $\theta_1$ must be a sum of three terms corresponding to the three base sets. 
	On the other hand, all base sets are constructed from the intersection of two sets defined by $\switchfun_1(x)$ and $\switchfun_2(x)$, see the third row in Table \ref{tab:step_summary}.
	Thus, the products in the expressions refer to the intersection of the sets.
	\color{black}
	Therefore, the related Filippov system reads as $\dot{x} = \theta_1 f_1(x) + \theta_2 f_2(x)$, where
	\begin{align*}
		\theta_1 & = \alpha_1 \alpha_2+ \alpha_1 (1- \alpha_2) + (1-\alpha_1) \alpha_2  = \alpha_1+(1-\alpha_1) \alpha_2 ,\\
		\theta_2  &=(1-\alpha_1) (1-\alpha_2).
	\end{align*}
	By direct calculation, we verify that $\theta_1,\theta_2 \geq 0$ and $\theta_1 + \theta_2 = 1$.
	In contrast to the previous examples, the union of sets introduces a sum in the expressions for $\theta_1$.
\end{example}
We generalize the reasoning above as follows.
Given $\Nsys$ total regions and the  matrix $F(x,u) \in \R^{n_x} \times \R^{n_F}$, where $n_F = 2^{\nswitchfun}$ is the number of possibly repeating columns.
We define the index sets $\mathcal{R}_k = \{ i \mid F_{\bullet,i}(x,u) = f_k(x,u) \}$, i.e., \textcolor{black}{the set if the indices of the columns of $F(x,u)$ equal to $f_k(x,u)$}.
Note that if we have no unions, then $n_F = \Nsys$ and $\mathcal{R}_k = \{k\}$ for all $k \in \mathcal{J}$.
Using these definitions, the expression for $\theta_i$ in \eqref{eq:theta_via_step} reduces to:
\begin{align*}
	\theta_i &= \sum_{k \in \mathcal{R}_i}  \prod_{j\in \C}  \frac{1-{S}_{k,j}}{2}+{S}_{k,j}\alpha_j.
\end{align*}

\subsection{A lifting algorithm for the multi-affine terms}\label{sec:step_lifting_algortihm}
\color{black}
It can be seen from \eqref{eq:theta_via_step}, the expression of $\theta_i$ consist of a product of $\nswitchfun$ affine terms (of the form of $\alpha_j$ and $1-\alpha_k$).
If $\nswitchfun$ is large, then this expression is very nonlinear.
\color{black}
To reduce the nonlinearity, we introduce auxiliary lifting variables as in the \textit{lifted} Newton's method \cite{Albersmeyer2010}, which iterates on a larger but less nonlinear problem.
Whenever there are more than two terms in the multi-affine expression for $\theta_i$, we introduce lifting variables 
$\beta_k$ and derive an equivalent formulation, which has only bilinear terms.
\color{black}
Now, instead of having $\Nsys$ multi-affine expressions, we have $\Nsys+n_\beta$ expressions, but which are less nonlinear.
\color{black}
We exploit the structure of the matrix ${S}$ and derive an easy-to-implement algorithm that automates the lifting procedure.
To give an idea of the final results we aim to obtain, we illustrate the lifting procedure with an example.

\begin{example}\label{ex:step_lifting}
	Regard a PSS with $\nswitchfun=3$ switching functions and $\Nsys = 8$ modes, i.e., the PSS regions match the basis sets, $R_i = \tilde{R}_i, \ i = \{1,\ldots,8\}$.
	The matrix ${S}\in \R^{8\times3}$, and the expression for the multipliers $\theta \in \R^8$ read as
	\begin{align*}
		{S} &=
		\begin{bmatrix}
			1 & 1 & 1\\
			1 & 1 & -1\\
			1 & -1 & 1\\
			1 & -1 & -1\\
			-1 & 1 & 1\\
			-1 & 1 & -1\\
			-1 & -1 & 1\\
			-1 & -1 & -1
		\end{bmatrix},\
		G_{\mathrm{F}}(\theta,\alpha) =
		\begin{bmatrix}
			\theta_1 \shortminus \alpha_1\alpha_2\alpha_3\\
			\theta_2 \shortminus \alpha_1\alpha_2	(1\shortminus\alpha_3)\\
			\theta_3 \shortminus \alpha_1(1\shortminus\alpha_2)\alpha_3\\
			\theta_4 \shortminus \alpha_1(1\shortminus\alpha_2)	(1\shortminus\alpha_3)\\
			\theta_5 \shortminus (1\shortminus\alpha_1)\alpha_2\alpha_3\\
			\theta_6 \shortminus (1\shortminus\alpha_1)\alpha_2	(1\shortminus\alpha_3)\\
			\theta_7 \shortminus (1\shortminus\alpha_1)(1\shortminus\alpha_2)\alpha_3\\
			\theta_8 \shortminus (1\shortminus\alpha_1)(1\shortminus\alpha_2)	(1\shortminus\alpha_3)
		\end{bmatrix}  =0.
	\end{align*}
	
	We can introduce the lifting variable	$\beta \in \R^4$ and obtain
	\begin{align*}
		&G_{\beta}(\alpha,\beta) =
		\begin{bmatrix}
			\beta_1 -  \alpha_1\alpha_2\\
			\beta_2 -  \alpha_1(1\shortminus\alpha_2)\\
			\beta_3 -  (1\shortminus\alpha_1)\alpha_2\\
			\beta_4 - (1\shortminus\alpha_1)(1\shortminus\alpha_2)
		\end{bmatrix}=0,\
		G_{\theta}(\theta,\alpha,\beta)
		=\begin{bmatrix}
			\theta_1 \shortminus \beta_1\alpha_3\\
			\theta_2 \shortminus \beta_1	(1\shortminus\alpha_3)\\
			\theta_3 \shortminus \beta_2\alpha_3\\
			\theta_4 \shortminus \beta_2(1\shortminus\alpha_3)\\
			\theta_5 \shortminus \beta_3\alpha_3\\
			\theta_6 \shortminus \beta_3(1\shortminus\alpha_3)\\
			\theta_7 \shortminus \beta_4\alpha_3\\
			\theta_8 \shortminus \beta_4(1\shortminus\alpha_3)
		\end{bmatrix}=0.
	\end{align*}
	The equation $G_{\beta}(\alpha,\beta)=0$ relates the lifting variables $\beta_i$ with the variables $\alpha_j$, whereas $G_{\theta}(\theta,\alpha,\beta)$ provides expressions for $\theta_i$ via $\beta_i$ and the remaining $\alpha_j$.
	By replacing $G_{\mathrm{F}}(\theta,\alpha)= 0$ with $G_{\mathrm{lift}}(\theta, \alpha,\beta)  \coloneqq
	(G_{\beta}(\alpha,\beta),	G_{\theta}(\theta,\alpha,\beta) ) =0$, we obtain an equivalent system of equations that only consists of bilinear terms.
\end{example}

\begin{algorithm}\label{algo:lifting}
	\caption{Lifting algorithm for the step DCS \eqref{eq:step_dcs}}
	\label{alg:lifting_dcs}
	\begin{algorithmic}[1]
		\State \textbf{Input:} ${S}, n_{\mathrm{d}}$
		\State \textbf{Initialize:}
		$\tilde{S} \leftarrow {S}$,
		$k  \leftarrow 0$;
		$\tilde{\theta} \leftarrow e \in \R^{n_f}$,
		$G_{\theta}(\theta,\alpha,\beta) \leftarrow [\;]$,
		$G_{\beta}(\alpha,\beta) \leftarrow [\;]$.
		\For{$j=1:\nswitchfun$}
		\State $\tilde{N} \leftarrow \sum_{j=1}^{\nswitchfun} |\tilde{S}_{\bullet,j}|$
		\State $\ILift \leftarrow \{ i \mid \tilde{N}_i = j \}$
		\State
		$\tilde{\theta} \leftarrow
		\tilde{\theta} \cdot  \Big(\frac{e-\tilde{S}_{\bullet,j}}{2}+{S}_{\bullet,j}^{\mathrm{temp}}\cdot \alpha_j \Big)$
		\If {$\ILift \neq \emptyset $}
		\State
		$ G_{\theta}(\theta,\alpha,\beta) \leftarrow
		(G_{\theta}(\theta,\alpha,\beta),
		\theta_{k+\ILift}- \tilde{\theta}_{\ILift})$
		\State
		Remove entries of $\tilde{\theta}$ with index in $\ILift$
		\State
		Remove rows of $\tilde{S}$ with index in $\ILift$
		\State
		$k \leftarrow k+\max(\ILift)$
		\EndIf
		\If{ $j \in \{ n_{\mathrm{d}},\ldots,\nswitchfun-1\}$}
		\State
		$\{ \I_{\mathrm{red}},\I_{\mathrm{full}}\} =
		\texttt{unique}(\tilde{S}_{\bullet,\{1,\ldots,j\}})$,
		\State $\beta \leftarrow (\beta,\beta^{j})$ where  $\beta^j \in \R^{|\I_{\mathrm{red}}|}$
		\State $G_{\beta}(\alpha,\beta) \leftarrow
		(G_{\beta}(\alpha,\beta)
		\beta^j - \tilde{\theta}_{\I_{\mathrm{red}}})$
		\State $\tilde{\theta} \leftarrow \beta^j_{\I_{\mathrm{full}}}$
		\EndIf
		\EndFor
		\State  $G_{\mathrm{Lift}}(\theta,\alpha,\beta) \coloneqq (G_{\theta}(\theta,\alpha,\beta),G_{\beta}(\alpha,\beta))$
		\State \textbf{Output:} $G_{\mathrm{Lift}}(\theta,\alpha,\beta)$ , $\beta$
	\end{algorithmic}
\end{algorithm}

We proceed by outlining a general lifting algorithm.
The expressions for $\theta_i$ consist of the product of $\nswitchfun$ affine terms.
Our goal is to have at most $n_{\mathrm{d}}$ terms in the multi-affine expression for $\theta_i$. 
For example, if we pick $n_{\mathrm{d}} =2$, we have only bilinear expressions in the equations defining $\theta$ and $\beta$.
Thus, the parameter $\nswitchfun \geq n_{\mathrm{d}} \geq 2$, controls the number of terms in the multi-affine expressions and implicitly the number of new lifting variables $\beta \in \R^{n_\beta}$.
Given the matrix ${S}$, our goal is to automatically obtain the constraint $G_{\mathrm{lift}}(\theta,\alpha,\beta)=0$.

The algorithm outlined above can be implemented using a symbolic framework such as \texttt{CasADi} \cite{Andersson2019}.
We provide the pseudo code in Algorithm \ref{alg:lifting_dcs}, which introduces the lifting algebraic variables $\beta$ and new \textit{lifted} expressions for $\theta_i$, namely $G_{\mathrm{lift}}(\theta,\alpha,\beta)$.
Note that we make use of three helper variables, the matrix $\tilde{S}$ and the vectors $\tilde{\theta}$ and $\tilde{N}$.
The matrix $\tilde{S}$ is a submatrix of ${S}$, where we have removed the rows with index $i$, for which we already have a (lifted) expression for $\theta_i$.
The vector $\tilde{N}$, defined in line 4, collects the number of nonzero entries of every row $\tilde{S}$.
In other words, it keeps track of how many terms are in the initial expressions for $\theta_i$, that are not yet lifted.

The main loop iterates from $j = 1$ to $\nswitchfun$ and provides in every iteration the expressions for all $\theta_i$ that have exactly $j$ terms in their multi-affine expression.
The index set  $\ILift = \{ i \mid \tilde{N}_i = j \}$, defined in line 5, contains the indices of $\theta$, that have exactly $j$ entries in their corresponding multi-affine expression.
In line 6, we define the auxiliary variable $\tilde{\theta}$, which stores the intermediate expressions for $\theta$ with up to $j$ terms in the product.
The index $k$ stores the index of the last $\theta_k$ for which a lifted expression was derived.
For $j \leq n_{\mathrm{d}}$ the expressions for $\theta_i$ are unaltered.
This is treated in lines 7-11.

As soon as  $j > n_{\mathrm{d}}$, the algorithm introduces new lifting variables $\beta^j$ (line 15) and changes the expression for $\tilde{\theta}$ accordingly.
This is done in lines 13-17.
A key tool is the function \texttt{unique} in line 14.
It is available in \texttt{MATLAB} and in the \texttt{numpy} package in \texttt{python}.
It works as follows: given a matrix $A \in \R^{m \times n}$ it returns a matrix $\tilde{A} \in \R^{p \times n}$, with $p \leq m$.
This is the matrix constructed from $A$ by removing its repeating rows.
More importantly for our needs, it returns the index sets $\I_{\mathrm{red}}$ and $\I_{\mathrm{full}}$, with $|\I_{\mathrm{red}}| = p$  and $|\I_{\mathrm{full}}| = m$.
The index sets have the properties
$A =  \begin{bmatrix}
	\tilde{A}_{i,\bullet} \mid i \in \I_{\mathrm{full}}
\end{bmatrix} \in \R^{m \times n} $
and
$\tilde{A} =  \begin{bmatrix}
	{A}_{i,\bullet} \mid i \in \I_{\mathrm{red}} \in \R^{p \times n}
\end{bmatrix}$.

This enables the use of the same $\beta^j$ for several $\theta_i$ if they share the same terms in the corresponding multi-affine expressions, cf. lines 15-17.
After the loop is finished, the algorithm outputs $G_{\mathrm{Lift}}(\theta,\alpha,\beta)$ and $\beta$.
One can verify that Algorithm \ref{alg:lifting_dcs} produces the same output as Example~\ref{ex:step_lifting}.
It can be shown, that for a given $n_{\mathrm{d}} < \nswitchfun$, the total number of new lifting variables is $n_{\beta} = 2^{\nswitchfun}-2^{n_{\mathrm{d}}}$.

\subsection{Comparisons of Stewart's and the Heaviside step reformulation}\label{sec:stewart_vs_step_complexity}
We compare Stewart's reformulation \eqref{eq:stewart_dcs} and the Heaviside step reformulation \eqref{eq:step_dcs} based on their total number of algebraic variables, complementarity constraints, and equality constraints for a given number of switching functions $\nswitchfun$.
The total number of regions (and multipliers $\theta_i$) in Stewart's reformulation is always $\Nsys = 2^{\nswitchfun}$.
\color{black}
In contrast to the Heaviside step reformulation, we cannot reduce the number of variables if the regions $R_i$ are defined as unions of the base regions $\tilde{R}_i$.
\color{black}
On the other hand, in the Heaviside step reformulation, depending on the geometry of the regions $R_i$, $\Nsys$ is an integer in $[2,2^{\nswitchfun}]$. 
In the step reformulation, we may introduce $n_{\beta}$ lifting variables to reduce the nonlinearity.
If $n_{\mathrm{d}} > \nswitchfun$, this leads to $n_{\beta} = 2^{\nswitchfun}-2^{n_{\mathrm{d}}}$ additional lifting variables and equations.

We compare now the number of algebraic variables. 
In Stewart's reformulation, we have $\lambda \in \R^{2^{\nswitchfun}}$ and $\mu \in \R$.
In the Heaviside step reformulation, we have $\alpha,\lambdap,\lambdan \in \R^{\nswitchfun}$.
Thus, the total number of algebraic variables in the former case is $n_{\mathrm{alg}}^{\mathrm{S}} = 2 \cdot 2^{\nswitchfun}+1$, and in the later case
$n_{\mathrm{alg}}^{\mathrm{H}} = \Nsys + 3 \nswitchfun+n_{\beta}$.
The number of complementarity constraints $n_{\mathrm{comp}}$ in Stewart's case is $n_{\mathrm{comp}}^{\mathrm{S}} = 2^{\nswitchfun}$, and in the Heaviside step case $n_{\mathrm{comp}}^{\mathrm{H}} = 2\nswitchfun$, i.e., we have exponential versus linear complexity.
Finally, in Stewart's reformulation, we have in total $n_{\mathrm{eq}}^{\mathrm{S}} = 2^{\nswitchfun} + 1$ equality constraints ($g_i(x) = \lambda_i - \mu$ and $e^\top \theta=1$).
In the Heaviside step case, there are $n_{\mathrm{eq}}^{\mathrm{H}} = \Nsys + n_{\beta} + \nswitchfun$ equality constraints, for the definitions of $\theta_i$, $\beta_i$ and the constraints $\switchfun_i(x) = \lambdap_i - \lambdan_i$, respectively.
The numbers of variables and constraints are summarized in Table~\ref{tab:stewart_vs_step}.

\begin{table}[t]	
	\centering
	\caption{Comparison of the problem sizes in Stewart's and the step reformulation for a fixed $\nswitchfun$.}
	\centering
	\begin{tabular}{@{}l|l|ll|ll@{}}	
			\specialrule{.15em}{0em}{0.0em} 
			\textbf{Ref.} & $n_{f}$  &$n_{\beta}$ &$n_{\mathrm{alg}}$ &$n_{\mathrm{comp}}$ &$n_{\mathrm{eq}}$\\
			\hline
			Stewart &  $2^{\nswitchfun}$    & $0$ & $2\cdot 2^{\nswitchfun}\!+\!1$ & $2^{\nswitchfun}$ & $2^{\nswitchfun}\!+\!1$ \\
			Step    & $[2,2^{\nswitchfun}]$ & $\!\begin{cases}2^{\nswitchfun}\!-\!2^{n_{\mathrm{d}}},\ n_{\mathrm{d}}\! \leq\! \nswitchfun\\	0 ,\ n_{\mathrm{d}} > \nswitchfun	\end{cases}$ &$\Nsys\!+\!3{\nswitchfun}\!+\!n_{\beta}$ & $2{\nswitchfun}$ & ${\nswitchfun}\!+\!n_{\beta}\!+\!\Nsys$\\
			\specialrule{.15em}{0em}{0.0em} 
		\end{tabular}
		\label{tab:stewart_vs_step}
	\end{table}
	\begin{figure}[t]
		\centering
		{\includegraphics[width= 0.98\columnwidth]{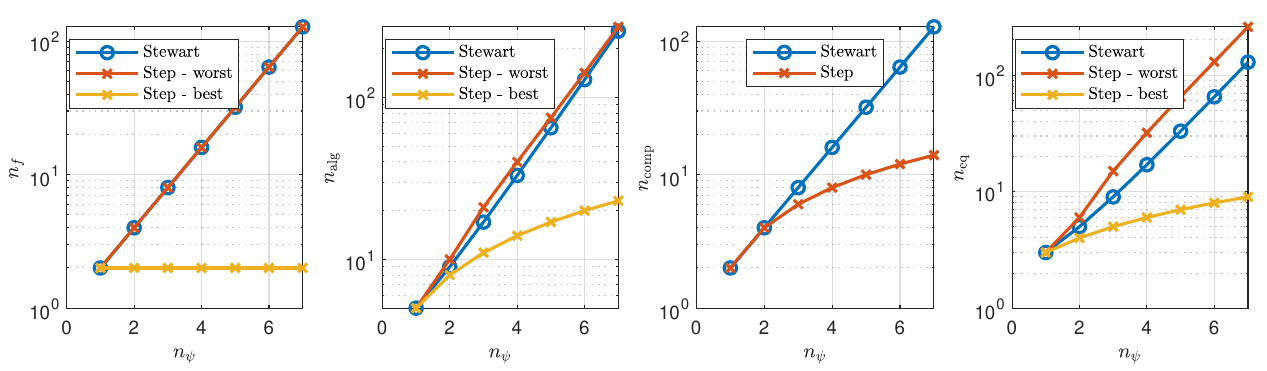}}
		\caption{Comparison of the complexities of Stewart's and the step reformulation.
		}
		\label{fig:stewart_vs_step}
	\end{figure}
	Figure \ref{fig:stewart_vs_step} illustrates the different quantities for several $\nswitchfun$. 
	We plot for the Heaviside step reformulation two extreme scenarios:
	\begin{enumerate}
		\item Worst complexity case - every basis set defines a PSS region, $n_f = 2^{\nswitchfun}$, we lift to have only bilinear terms, i.e., $n_{\mathrm{d}} = 2$ (maximizes the number of lifting variables) - red line in the plots.
		\item Best complexity case - no lifting and only two regions ($n_f = 2$) for all $\nswitchfun$ -yellow line in the plots. 
	\end{enumerate}
	Note that in both cases the Heaviside step reformulation has the same number of complementarity constraints.
	For smaller values of $\nswitchfun$, both reformulations have similar complexity. 
	For a large number of switching functions, the step reformulation has fewer variables. 
	However, if there is no lifting, the problem can become very nonlinear for large $\nswitchfun$.
	\color{black}
	An extensive numerical comparison of the use of the two approaches in simulation and optimal control problems has been done in \cite[Section 5.8]{Pozharskiy2023}. 
	It turns out that statistically both formulations have similar performance, but on specific examples one can outperform the other. 
	In addition, the Heaviside step reformulation provides more modeling flexibility.
	\color{black}
\section{Numerical examples}\label{sec:numerical_examples}
\color{black}
In this section, we consider two numerical examples. 
First, we consider a numerical simulation example of a gene regulatory network and empirically demonstrate that the FESD method has higher order integration accuracy.
Second, we consider an optimal control example of a planar two-link monoped that must reach a certain goal in the horizontal direction. 
We compare the Heaviside step reformulation with the Stewart reformulation in terms of total computational time and cost per iteration.
\color{black}
\subsection{Gene regulatory networks}\label{sec:numerical_simulation_gen_reg}
In the simulation experiment we consider a gene regulatory network model, which is described by a DI of the form \eqref{eq:di_ap}.
The model is not a Filippov system, but a more general DI. 
More numerical examples with FESD for the step reformulation can be found in \cite{Nurkanovic2023a}.
In this reference, we have generated integration order plots for a Filippov system and confirmed empirically the results of Theorem~\ref{th:integration_order_step}.

\begin{figure}
	\centering
	\includegraphics[width=0.85\columnwidth]{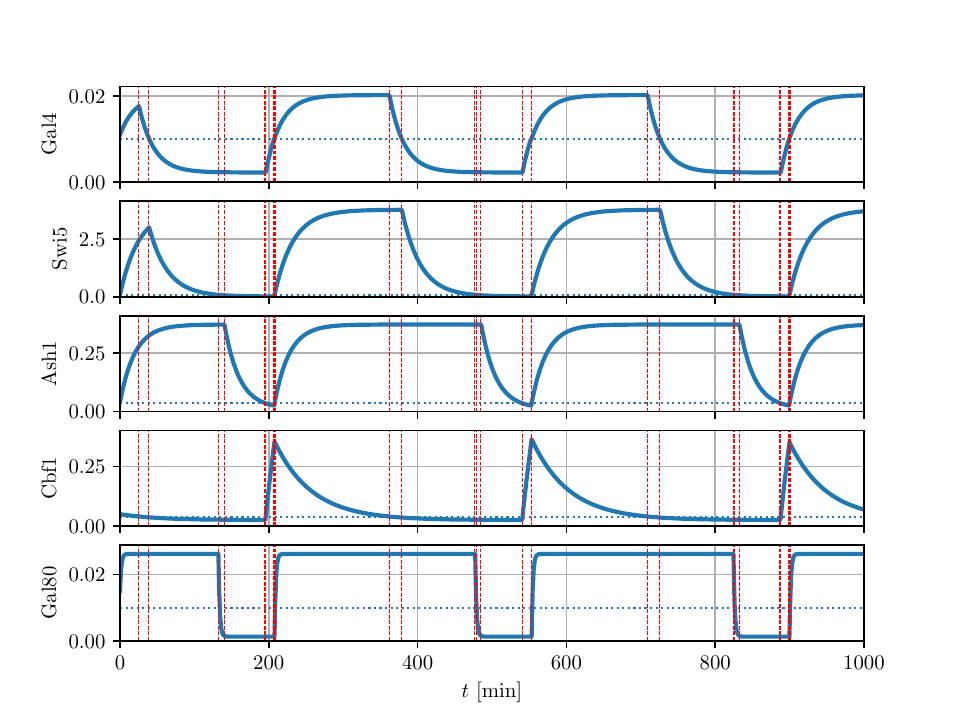}
	\caption{State trajectory of the IRMA example with red dotted vertical lines indicating the switches.}
	\label{fig:irma_states}
\end{figure}

We regard the IRMA example, a synthetic network composed of five genes, originally proposed in~\cite{Cantone2009}.
This example is inspired by~\cite{Acary2014}, which includes more examples with Heaviside step functions that are implemented in \nosnoc~\cite{nosnoc_py}.
\paragraph*{IRMA model}
In Figure~\ref{fig:irma_states}, we reproduce the state trajectories from \cite[Fig. 11]{Acary2014}.
Additionally, the vertical lines show the switching times of the selection variables $\alpha_i$.
The states of this system are the protein concentrations of Gal4, Swi5, Ash1, Cbf1, and Gal80, which are denoted by $x_1,\dots, x_5$.
There are seven switching functions, which are defined by the states crossing certain thresholds, which are plotted as horizontal lines in Figure~\ref{fig:irma_states}.
Specifically, the switching functions are
\begin{align*}
&\switchfun(x) = (x_1 - 0.01,  x_2 - 0.01,  x_2 - 0.06,  x_2 - 0.08,  x_3 - 0.035,  x_4 - 0.04,  x_5 - 0.01). 
\end{align*}
Note that Swi5 ($x_2$) has three threshold values.
The continuous-time dynamics of the system are given by
\begin{subequations}
\begin{align}
  \dot{x}_1 & \in - p_1 x_1 + \kappa_1^1 + \kappa_1^2 \gamma(\switchfun_6(x)), \\
  \dot{x}_2 & \in - p_2 x_2 + \kappa_2^1 + \kappa_2^2 \gamma(\switchfun_1(x)) (1-u) \gamma(\switchfun_7(x)), \\
  \dot{x}_3 & \in - p_3 x_3 + \kappa_3^1 + \kappa_3^2 \gamma(\switchfun_3(x)), \\
  \dot{x}_4 & \in - p_4 x_4 + \kappa_4^1 + \kappa_4^2 \gamma(\switchfun_2(x)) (1-\gamma(\switchfun_5(x))), \\
  \dot{x}_5 & \in - p_5 x_5 + \kappa_5^1 + \kappa_5^2 \gamma(\switchfun_4(x)).
\end{align}
\end{subequations}
Note that $u \in \{0, 1\}$ is an external input, which is set to $u = 1$ in the scenario considered here.
The initial state is given by $x_0 = (0.011, 0.09, 0.04, 0.05, 0.015)$.
The parameter values are given as
\begin{align*}
  p &= (0.05, 0.04, 0.05, 0.02, 0.6), \\
  \kappa^1 &= (1.1 \cdot 10^{-4}, 3 \cdot 10^{-4}, 6 \cdot 10^{-4}, 5 \cdot 10^{-4}, 7.5 \cdot 10^{-4})\\
  \kappa^2 &= (9 \cdot 10^{-4}, 0.15, 0.018, 0.03, 0.015).
\end{align*}
%

\paragraph*{Integration order experiment}
\begin{figure}
    \includegraphics[width=0.49\columnwidth]{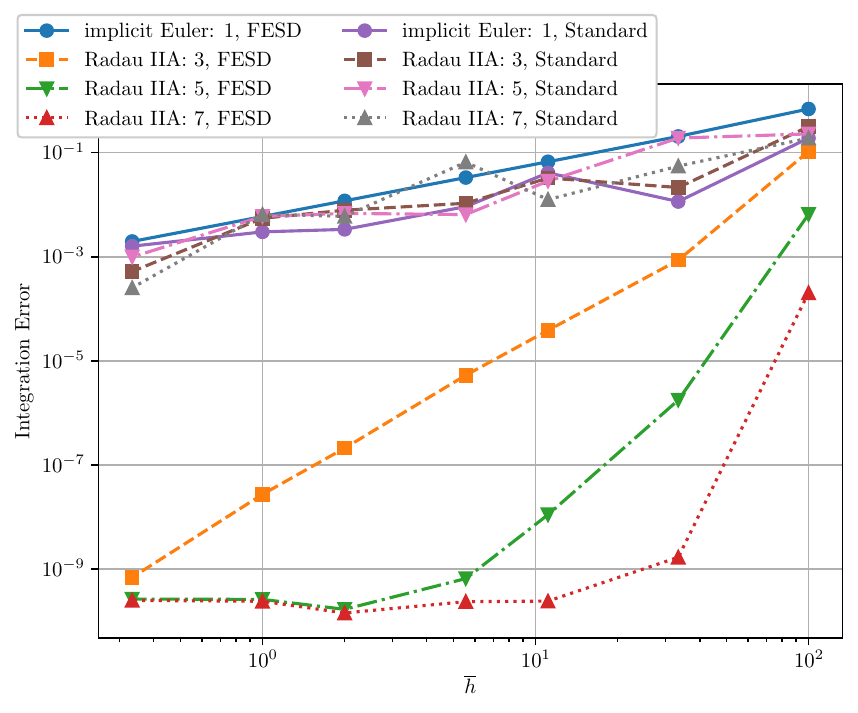}
    \includegraphics[width=0.49\columnwidth]{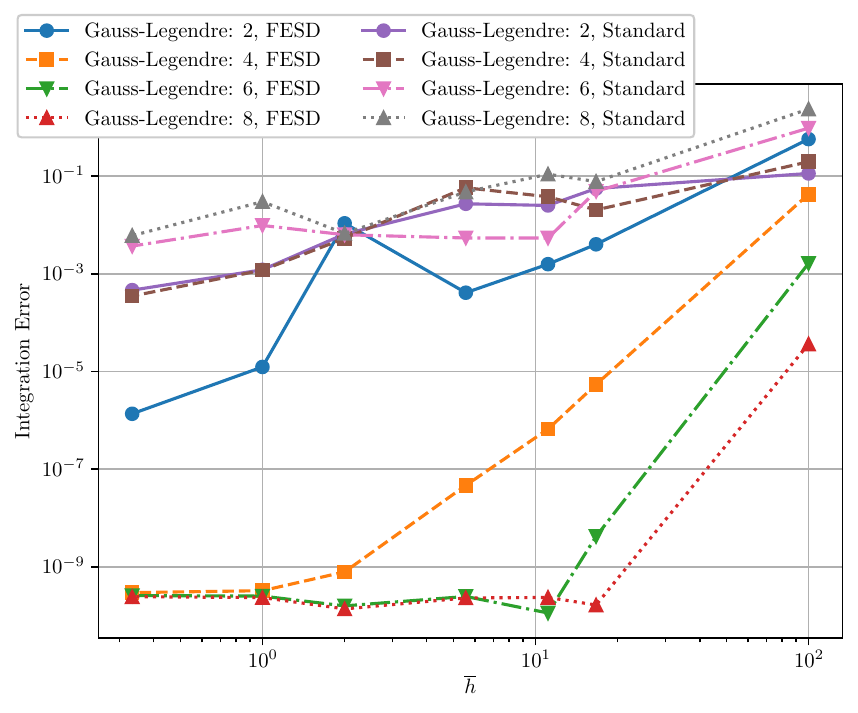}
    \caption{Accuracy vs. step size: Simulation of first 100 minutes of the trajectory in Figure~\ref{fig:irma_states} with different RK schemes and step sizes.}
    \label{fig:irma_bm}
\end{figure}
To showcase that the FESD step reformulation preserves the integration order of the underlying Runge-Kutta method, we simulate the first 100 minutes of the trajectory depicted in Figure~\ref{fig:irma_states}.
This time interval contains exactly two switches.
We use $ \NFE = 3 $, to capture those switches even with a single FESD step. We plot the integration error over the average step size $\bar{h} = {T_{\mathrm{sim}}}/{\Nsim\NFE}$, where $\Nsim$ is the total number of simulation steps.
The results are summarized in Figure~\ref{fig:irma_bm} for FESD with underlying Radau IIA of orders 1, 3, 5, and 7, and Gauss-Legendre methods of orders 2, 4, 6, and 8.
It can be seen that the FESD step formulation preserves the integration order of the underlying Runge-Kutta method, while with the standard time-stepping approach, using methods that typically have higher-order integration accuracy degrade to order one without FESD.
This experiment shows, that if a feasible solution is found the integration order of the underlying Runge-Kutta method is preserved.
The implementation of the example is publicly available\footnote{\url{https://github.com/FreyJo/nosnoc_py/blob/main/examples/Acary2014/irma_integration_order_experiment.py}}.



\subsection{Optimal control example with state jumps}\label{sec:robot_ocp}
Now we further investigate the use of FESD for the Heaviside step reformulation by applying it to an optimal control problem. 
We regard the problem of synthesizing dynamic motions of the two-link \textit{Capler} robot with state jumps and friction~\cite{Carius2018}. 
Systems with state jumps do not directly fit the form of ODEs with set-valued Heaviside step functions~\eqref{eq:di_ap}.
However, using the time-freezing reformulation we can transform the system with state jumps into a PSS of the form of~\eqref{eq:pss}~\cite{Nurkanovic2021} and use the methods developed in this paper, 
\color{black}
see also~Remark~ \ref{rem:time_freezing}.
\color{black}
The monoped's configuration is described by four degrees of freedom $q = (q_x,q_z,\phi_\mathrm{knee},\phi_\mathrm{hip})$, where $(q_x,q_z)$ are the coordinates of the monoped's base at the hip and $\phi_\mathrm{knee}$, $\phi_\mathrm{hip}$ are the angles of the hip and knee.
The robot is actuated by two direct-drive motors at the hip and knee joints. 
The control variables are the torques of these motors $u(t)=(u_{\mathrm{knee}}(t),u_{\mathrm{hip}}(t))$.
Denote by $p_{\mathrm{foot}}(q) = (p_{\mathrm{foot},x}(q),p_{\mathrm{foot},z}(q))$ and $p_{\mathrm{knee}}(q) = (p_{\mathrm{knee},x}(q),p_{\mathrm{knee},z}(q))$ the kinematic position of the robot's foot and knee, respectively. 
We model a single contact point, the tip of the robot's foot touching the ground, which is expressed via the unilateral constraint 
$$f_c(q) = p_{\mathrm{foot},z}(q) \geq 0.$$ 
Moreover, we denote the expression summarizing the Coriolis, control, and gravitational forces by $f_v(q,u)$, the inertia matrix by $M(q)$, the normal contact Jacobian by $J_\mathrm{n}(q)\in \R^{n_q}$, and the contact tangent by $J_\mathrm{t}(q)\in \R^{n_q}$.
The detailed derivation of all mentioned functions, i.e., the model equations, kinematic expressions, and all parameters for the robot can be found in~\cite[Appendix A]{Gehring2011}.

After the time-freezing reformulation, the PSS state consists of $x(\tau)=(q(\tau),v(\tau),t(\tau)) \in \R^9$, with $q(\tau)$ being the position, $v(\tau)$ the velocity, $t(\tau)$ is a clock state needed for the time-freezing reformulation and $\tau$ is the time of the ODE, cf. \cite{Nurkanovic2021}.
For the time-freezing PSS, we have in total \secsub{three} switching functions: the gap function, as well as the normal and tangential contact velocities~\cite{Nurkanovic2021}:
\begin{align*}
	\psi(x) = (f_c(q), J_\mathrm{n}(q)^\top v, J_\mathrm{t}(q)^\top v).
\end{align*}
This allows a definition of eight base sets via the sign matrix (cf. Definition \ref{def:basis_sets}):
\begin{align*}
	{S} &=
	\begin{bmatrix}
		1 & 1 & 1\\
		1 & 1 & -1\\
		1 & -1 & 1\\
		1 & -1 & -1\\
		-1 & 1 & 1\\
		-1 & 1 & -1\\
		-1 & -1 & 1\\
		-1 & -1 & -1\\
	\end{bmatrix}
\end{align*}

However, the time-freezing PSS has only three regions, where the first one consists of the union of the first six base sets, and the other two match the two remaining base sets, i.e.:
\begin{align*}
	R_1 &= \cup_{i=1}^{6} \tilde{R}_i = \{x \in \R^{n_x} \mid f_c(q)>0\}\cup \{x \in \R^{n_x} \mid f_c(q)<0,J_\mathrm{n}(q)^\top v>0\},\\
	R_2 &= \tilde{R}_7= \{x \in \R^{n_x} \mid f_c(q)<0,J_\mathrm{n}(q)^\top v<0,J_\mathrm{t}(q)^\top v>0\},\\
	R_3 &= \tilde{R}_8 = \{x \in \R^{n_x} \mid f_c(q)<0,J_\mathrm{n}(q)^\top v<0,J_\mathrm{t}(q)^\top v<0\}.
\end{align*}
In region $R_1$, we define the unconstrained (free flight) dynamics of the monoped, and in regions, $R_2$ and $R_3$, auxiliary ODEs that mimic state jumps in normal and tangential directions due to frictional impacts:
\begin{align*}
	&f_1(x,u)  = (q, M(q)^{-1}f_v(q,u) , 1),\\
	&f_2(x)  = (\mathbf{0}_{4,1},  M(q)^{-1}(J_{\mathrm{n}}(q) - J_{\mathrm{t}}(q)\mu)a_{\mathrm{n}}, 0), \\
	&f_3(x)  = (\mathbf{0}_{4,1}, M(q)^{-1}(J_{\mathrm{n}}(q) + J_{\mathrm{t}}(q)\mu)a_{\mathrm{n}} , 0).
\end{align*}

\begin{figure}[t]
	\centering
	\centering
	\vspace{0.2cm}
	{\includegraphics[scale=0.28]{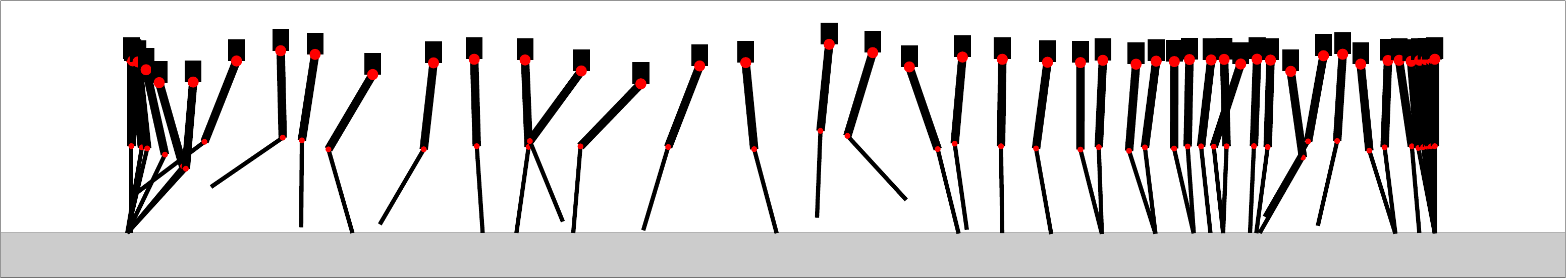}}
	\caption{Illustration of several frames of the solution of the discretized OCP.}
	\label{fig:solution_ocp}
	\centering
\end{figure}
\begin{figure}[t]
	\centering
	\centering
	\vspace{0.2cm}
	{\includegraphics[scale=0.38]{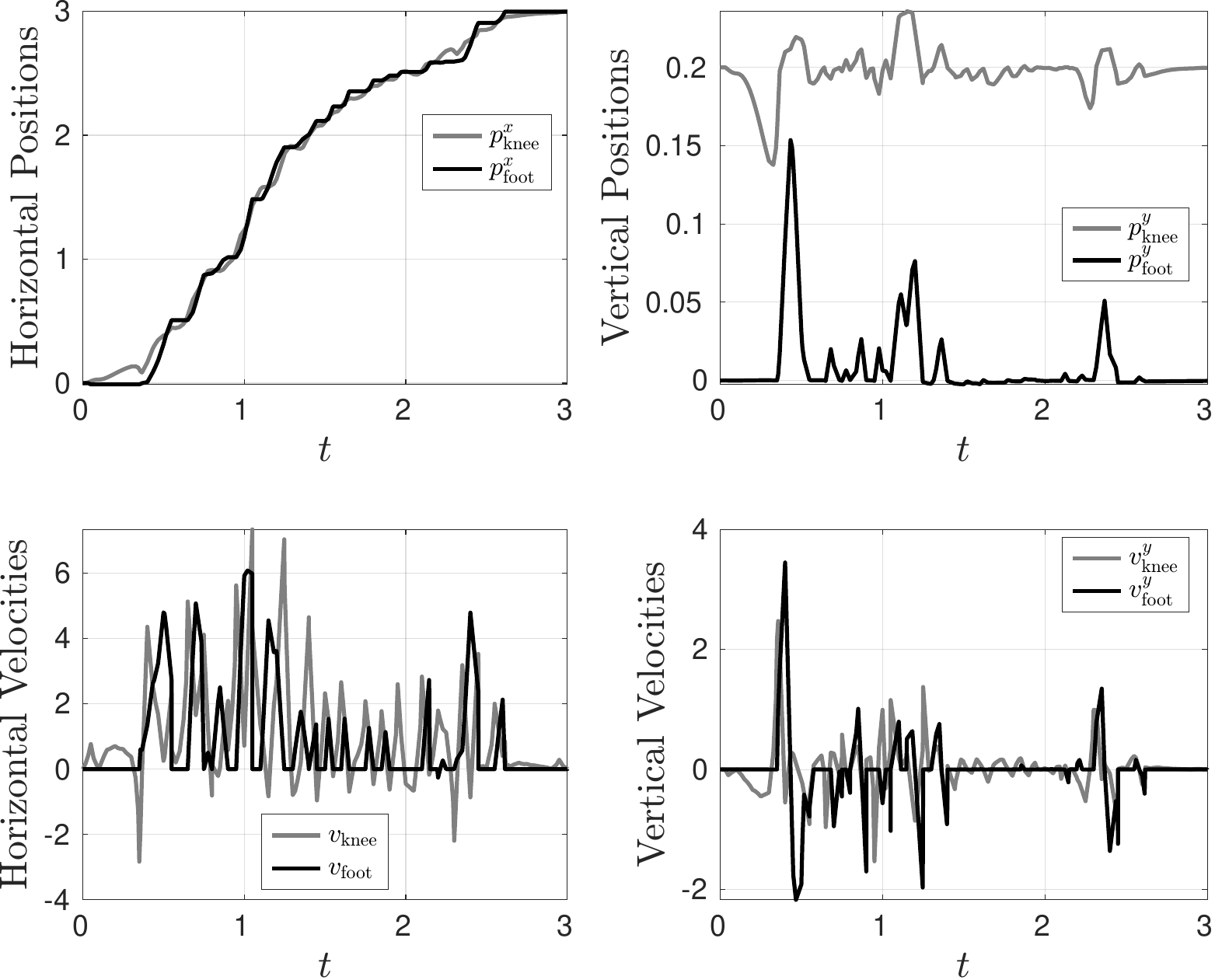}}
	\caption{The monoped's vertical and horizontal positions (top plots), and velocities (bottom plots).}
	\label{fig:kinematics_ocp}
\end{figure}

Here, $a_{\mathrm{n}}=200$ is the auxiliary ODE's constant~\cite{Nurkanovic2021} and the coefficient of friction is $\mu = 0.8$.
Observe that the clock state dynamics are $\frac{\dd t}{\dd \tau} = 1$ in $R_1$, and $\frac{\dd t}{\dd \tau} = 0$ in $R_2$ and $R_3$.
A solution trajectory of a PSS is continuous in time. 
By taking the pieces of the trajectory where $\frac{\dd t}{\dd \tau}>0$, we recover the solution of the original system with state jumps~\cite{Nurkanovic2021}.

It follows from the discussion in Section \ref{sec:compact_rep_and_lifting}, that with the Heaviside step reformulation, we have $\theta \in \R^{3}$.
Since we have $\nswitchfun =3$, the step DCS \eqref{eq:step_dcs} in this example has three complementarity constraints and six equality constraints.
For Stewart's reformulation, we cannot exploit the union of basis sets but must define eight regions, which are equal to the basis sets.
\secsub{The first six regions are equipped with $f_1(x,u)$ and the remaining two with $f_2(x)$ and $f_3(x)$ respectively.}
In this case, we have $\theta \in \R^8$, eight complementarity constraints and nine equality constraints in the DCS~\eqref{eq:stewart_dcs}.
Using the Heaviside step reformulation for this problem allows us to reduce the number of regions we need to define from eight to only three.

\begin{figure}[t]
	\centering
	\centering
	{\includegraphics[scale=0.27]{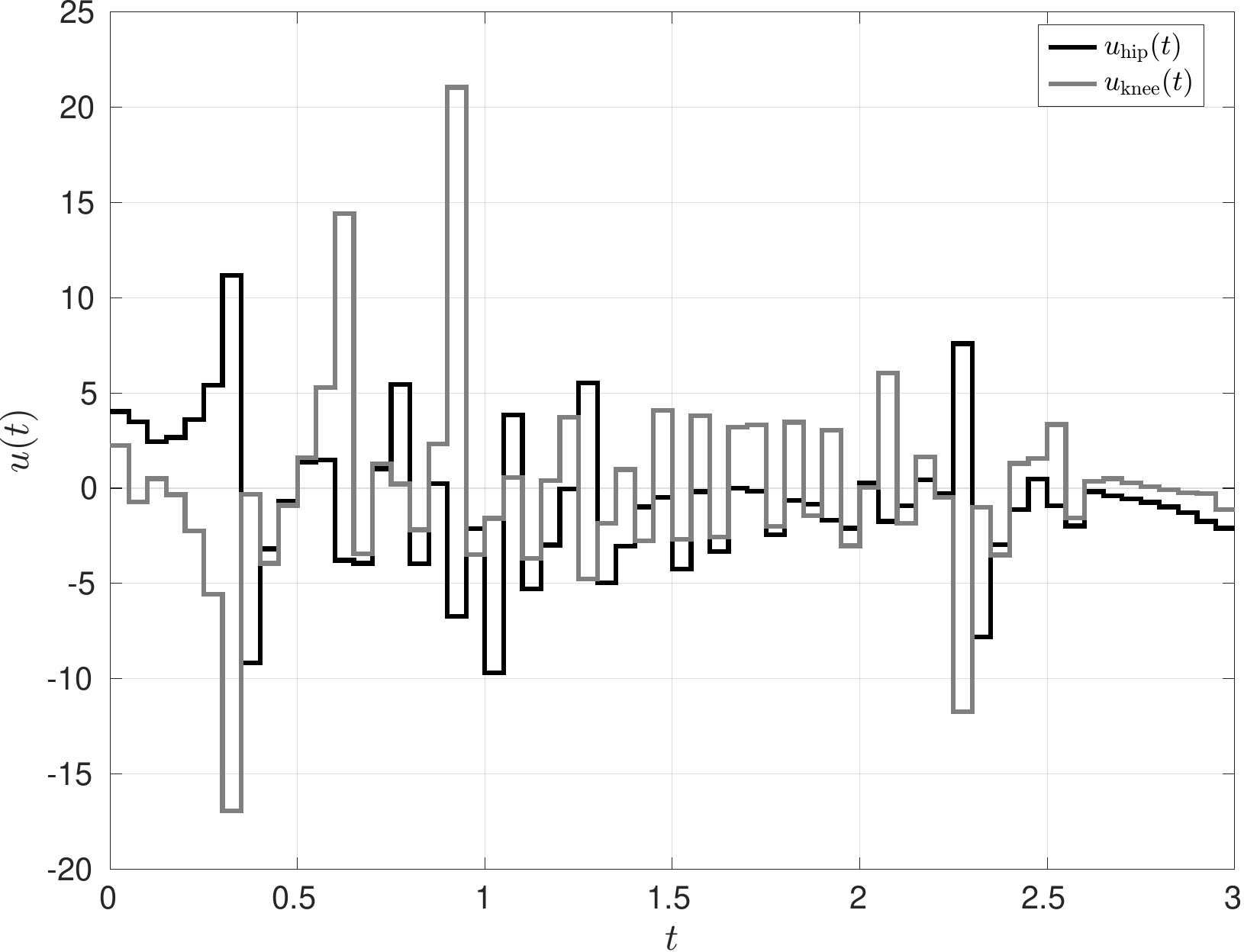}}
	\caption{The optimal controls for the example of $N = 60$ control stages.}
	\vspace{-0.5cm}
	\label{fig:monoped_ocp_controls}
	\centering
\end{figure}

To compare the performances of the two reformulations, we run an experiment in which the robot reaches a target position $q_{\mathrm{target}}=(3,0.4,0,0)$ with zero velocity $v_{\mathrm{target}} =\mathbf{0}_{4,1}$ in $T= 3.0$ seconds.
The initial state is $x(0)  = (0,0.4,0,0,0,0,0,0,0)$.
Given a reference $x^{\mathrm{ref}}(t)$, which is a spline interpolation between the initial and final position, and includes two jumps, we define the least-squares objective with the running and terminal costs
\begin{align*}
	L(x(\tau),u(\tau)) &=(x(\tau)-x^{\mathrm{ref}}(\tau))^\top Q (x(\tau)-x^{\mathrm{ref}}(\tau)) + \rho_u u(\tau)^\top u(\tau),\\
	R(x(T)) &= (x(T)-x^{\mathrm{ref}}(T)^\top Q_T (x(T)-x^{\mathrm{ref}}(T)),
\end{align*}
where the weight matrices are
$Q =  \mathrm{diag}(1, 1, 10, 1, 10^{-6}, 10^{-6}, 10^{-6}, 10^{-6}, 0)$,
$\rho_u = 0.01$, and $Q_T = \mathrm{diag}(10^3, 10^3, 10^3, 10^3, 10, 10, 10, 10, 0)$.
We define the state and control bound constraints:
	\begin{align*}
		& x_{\mathrm{lb}} \leq  x(\tau) \leq x_{\mathrm{ub}},\\
		&u_\mathrm{lb}\leq u(\tau) \leq u_{\mathrm{ub}},
	\end{align*}
where $x_{\mathrm{ub}} = (3.5, 10, \pi, \pi, 100, 100, 100, 100,\infty)$,
$x_{\mathrm{lb}} = (-0.5, 0, -\pi, -\pi, -100, -100, -100, -100, -\infty)$,
$u_{\mathrm{ub}} = (100,100)$, and $u_{\mathrm{lb}} =- u_{\mathrm{ub}}$.

Collecting all the above, we can define an OCP of the form of \eqref{eq:ocp}, which we discretized with the FESD Radau IIA scheme of order 3 ($\Nstg = 2$), with $\NFE =3$ finite elements on every control interval.
This OCP is discretized and solved with \nosnoc\ in a homotopy loop with \texttt{IPOPT}~\cite{Waechter2006}.
Figure~\ref{fig:solution_ocp} illustrates several frames of an example solution ($N=60$).
Figure~\ref{fig:kinematics_ocp} shows the relevant vertical and horizontal positions and velocities, and Figure \ref{fig:monoped_ocp_controls} shows the optimal controls.
\color{black}
From these plots we can see that the robots make several small jumps, mostly landing with low vertical velocities to reduce energy dissipation, and successfully reaches the target position with zero terminal velocity.
\color{black}

\begin{figure}[t]
	\centering 
	\includegraphics[width=0.45\linewidth]{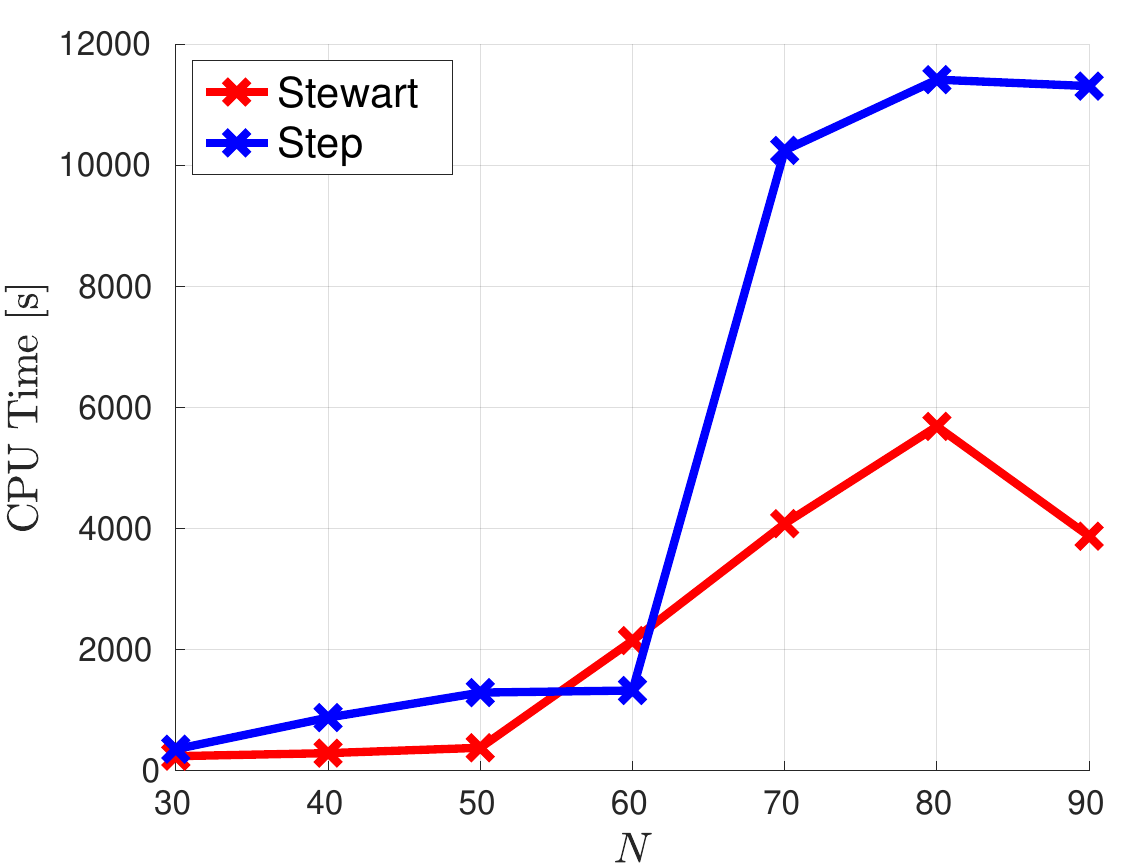}
	\includegraphics[width=0.45\linewidth]{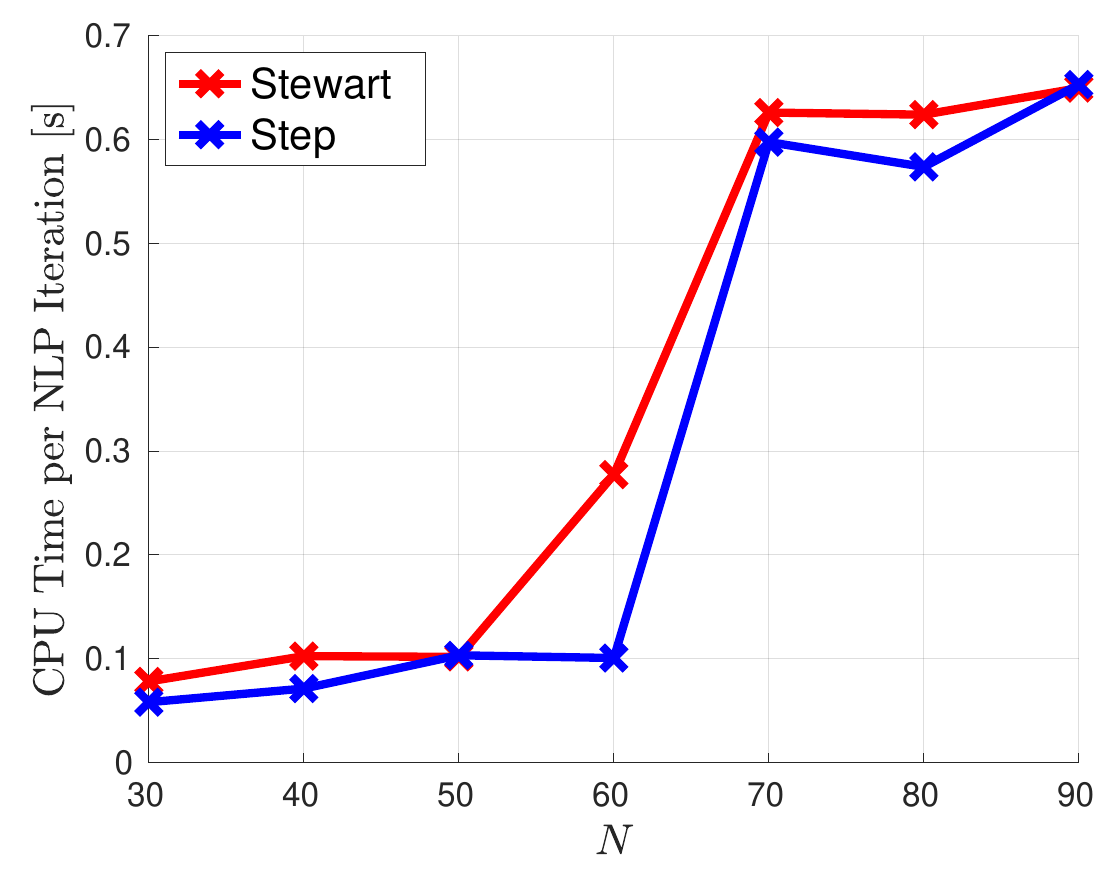}
	\caption{The total CPU time of the homotopy loop for different $N$ (left plot) and the CPU time per NLP solver iteration (right plot).}
	\label{fig:ocp_comparisson}
	\vspace{-0.6cm}
\end{figure}

Next, we solve this OCP for different values $N$ (number of control intervals) from 30 to 90 in increments of 10 and compare FESD for the Heaviside step reformulation (this paper) to FESD derived for Stewart's reformulation~\cite{Nurkanovic2024a}. 
We compare the  CPU time per NLP iteration and total CPU time for both approaches in Figure~\ref{fig:ocp_comparisson}.
As expected, due to the smaller number of variables and constraints, the Heaviside step reformulation leads to faster NLP iterations than the Stewart reformulation.
In our experiments, we try different linear solvers from the \texttt{HSL} library \cite{HSL}.
The linear solver choice has a particularly strong influence on performance for both reformulations. 
After evaluation, we select the best linear solver for each one, which is \texttt{ma27} for Stewart reformulation and \texttt{ma97} for the Heaviside step reformulation.
We can see the linear algebra costs increase significantly, for $N > 60$, as the linear solvers need more time to factorize larger systems.
The total computation time is influenced by several factors: homotopy algorithm, initialization, NLP solver performance, linear algebra solver, etc., and as such shows a less clear trend. 
In this example, Stewart's reformulation achieves a consistently lower computation times despite a higher cost per iteration.  
On the other hand, in our experiments in \cite{Nurkanovic2023a}, we used a similar setting and the Heaviside step reformulation led to lower total computational times.
We can conclude that which reformulation is the better depends strongly on the example and the chosen optimization algorithm parameters.
\color{black}
Further comparisons of the two approaches can be found in \cite[Section 5.2]{Pozharskiy2023}.
\color{black}


\section{Summary}\label{sec:conclusion}
This paper extends the Finite Elements with Switch Detection (FEDS) method to nonsmooth dynamical systems with set-valued Heaviside step functions.
The step functions allow to encode logical relations within a dynamical system. 
These systems cam be equivalent to Filippov systems, and this paper focuses on this case.
The set-valuedness of the Heaviside step functions leads to a differential inclusion. 
However, by using the optimality condition of a parametric linear program whose solution map corresponds to the step function, one obtains an equivalent Dynamic Complementarity System (DCS). 
Now the nonsmoothness and combinatorial structure is expressed algebraically and is thus suitable for numerical computations.
In the derivation of FESD, we exploit the continuity of the Lagrange multipliers in the DCS, which in turns enables the accurate detection the nonsmooth transition in time. 
This is also necessary for the correct computation of sensitivities of the discretized nonsmooth system.
We show how to apply the FESD method to discretize optimal control problems subject to nonsmooth systems with set-valued step functions and provide some convergence results of the method in simulation problems.
Compared to Stewart's reformulation used in~\cite{Nurkanovic2024a}, the reformulation considered here allows for a more compact reformulation of the same system, but with a smaller number of algebraic variables.
The comparison on an optimal control example shows that the new approach has a lower cost per iteration.
In summary, we extend the FESD method to a problem class that allows more modeling flexibility than in our initial study in \cite{Nurkanovic2024a}.
Our claims are verified by numerical examples.
An implementation of the new method is provided in the open-source software package~\nosnoc~\cite{Nurkanovic2022b}.
\color{black}


\end{document}